\documentclass[11pt]{amsart}
\usepackage[margin= 1in]{geometry}
\usepackage{graphicx,dsfont}
\usepackage{amssymb}
\usepackage{mathtools}
\usepackage{hyperref}
\usepackage{bm}
\usepackage{soul}

\usepackage{xcolor}
\usepackage[title]{appendix}
\usepackage{amsmath, mathrsfs, stmaryrd, color,bbm}

\let\oldmathscr\mathscr

\DeclareFontFamily{U}{BOONDOX-calo}{\skewchar\font=45 }
\DeclareFontShape{U}{BOONDOX-calo}{m}{n}{
  <-> s*[1.05] BOONDOX-r-calo}{}
\DeclareFontShape{U}{BOONDOX-calo}{b}{n}{
  <-> s*[1.05] BOONDOX-b-calo}{}
\DeclareMathAlphabet{\mathcalboondox}{U}{BOONDOX-calo}{m}{n}

\SetMathAlphabet{\mathcalboondox}{normal}{U}{BOONDOX-calo}{m}{n}
\DeclareMathAlphabet{\mathbcalboondox}{U}{BOONDOX-calo}{m}{n}

\usepackage{float}
\usepackage{verbatim}
\usepackage{comment}
\usepackage[utf8]{inputenc}
\usepackage[english]{babel}
\usepackage{enumitem}
\usepackage{booktabs}

\usepackage{csquotes}
\usepackage[style=alphabetic, giveninits=true, maxbibnames=4, doi=false,isbn=false,url=false]{biblatex}
\renewbibmacro{in:}{}
\DeclareFieldFormat[article]{title}{#1} 
\DeclareFieldFormat[article]{pages}{#1} 
\AtBeginBibliography{\small}
\bibliography{KPZsheet.bib} 

\newtheorem{theorem}{Theorem}

\newtheorem{claim}[theorem]{Claim}
\newtheorem{definition}[theorem]{Definition}
\newtheorem{proposition}[theorem]{Proposition}
\newtheorem{corollary}[theorem]{Corollary}
\newtheorem{conjecture}[theorem]{Conjecture}
\newtheorem{lemma}[theorem]{Lemma}

\newcommand{\N}{\mathbb{N}}

\newcommand{\whitenoise}{\ensuremath{\oldmathscr{\dot{W}}}}

\numberwithin{theorem}{section}
\numberwithin{equation}{section}

\begin{document}
\title[]{The KPZ equation and the directed landscape}
\author{Xuan Wu}
\address{Department of Mathematics, University of Illinois Urbana-Champaign, Altgeld Hall 
1409 W Green St, Urbana, IL 61801}
\email{xuanzw@illinois.edu}

\begin{abstract}
This paper proves the convergence of the narrow wedge solutions of the KPZ equation to the Airy sheet and the directed landscape in the locally uniform topology. This is the first convergence result to the Airy sheet and the directed landscape established for a positive temperature model. We also give an independent proof for the convergence of the KPZ equation to the KPZ fixed point for general initial conditions in the locally uniform topology. Together with the directed landscape convergence, we show the joint convergence to the KPZ fixed point for multiple initial conditions.
\end{abstract}

\vspace{1.5cm}

\maketitle
\section{Introduction}
\subsection{Kardar-Parisi-Zhang equation}
The KPZ equation was introduced in 1986 by Kardar, Parisi and Zhang \cite{KPZ} as a model for random interface growth and it describes the evolution of a function $\mathcal{H}(t, y)$ recording the height of an interface at time $t$ above position $y$. The KPZ equation is written formally as a stochastic partial differential equation (SPDE), 
\begin{equation}\label{eq:KPZ}
\partial_t\mathcal{H} = \frac{1}{2}\partial_y^2 \mathcal{H} + \frac{1}{2}(\partial_y\mathcal{H})^2 + \whitenoise,
\end{equation}
where $\whitenoise$ is a space-time white noise on $\mathbb{R}^2$ (for mathematical background or literature review, see \cite{Cor, QS15} for instance). The KPZ equation \eqref{eq:KPZ} is a canonical member of the associated KPZ universality class. Recently, there have been significant advancements in understanding the large-time behavior of solutions to the KPZ equation. Quastel-Sarkar \cite{QS} and Vir\'{a}g \cite{V} independently proved that, for a given initial function, the solution to \eqref{eq:KPZ} converges to a universal Markov process known as the KPZ fixed point \cite{MQR}.

In this paper we prove the convergence of the KPZ equation to the directed landscape \cite{DOV} in the locally uniform topology. This allows us to provide an independent proof for the convergence of the KPZ equation to the KPZ fixed point for general initial conditions. Moreover, the KPZ fixed point itself does not fully capture the large-time dynamics of the KPZ equation. A natural subsequent question arises: if we have two (or more) solutions with different initial functions and starting at different times, how can we describe the joint law of these solutions at a large time? It should be noted that these solutions are correlated through the same white noise in \eqref{eq:KPZ}. We also establish the joint convergence to the KPZ fixed point in the locally uniform topology for multiple initial conditions.

Before delving into our main results in Subsection~\ref{sec:mainresult}, we introduce related objects in the following two subsections.

\subsection{The narrow wedge solution}

The KPZ equation is related to the stochastic heat equation (SHE) with multiplicative noise through the Hopf–Cole transformation. Denote $\mathcal{Z}(t, y)$ as the solution to the following SHE, 
\begin{equation}\label{eq:SHE}
\partial_t \mathcal{Z} = \frac{1}{2}\partial_y^2 \mathcal{Z} + \mathcal{Z} \whitenoise.
\end{equation} 
The Hopf-Cole solution to the KPZ equation \eqref{eq:KPZ} is defined by taking $\mathcal{H}(t,y)=\log \mathcal{Z}(t,y).$ The fundamental solutions to the SHE \eqref{eq:SHE} are of great importance. For fixed $(s,x)\in\mathbb{R}^2$, we denote by $\mathcal{Z}(s,x;t,y)$ the solution to \eqref{eq:SHE} on $t>s$, $y\in\mathbb{R}$ with the delta initial condition $\mathcal{Z}(s,x;s,y)=\delta(y-x)$. We take a logarithm and define the \textbf{narrow wedge solutions} to \eqref{eq:KPZ}:
\begin{align}\label{equ:0309}
\mathcal{H}(s,x;t,y)\coloneqq \log \mathcal{Z}(s,x;t,y).
\end{align}
Let us note that in the literature, the narrow wedge solution often stands for the $2$-variable random function $\mathcal{H}(0,0;t,y)$. In contrast, in this paper the narrow wedge solution is a $4$-variable random function.

Let $\mathbb{R}^4_+\coloneqq \{(s,x,t,y)\in\mathbb{R}^4\, |\, s<t\}$. It is recently proved in \cite{AJRS} that all $\mathcal{Z}(s,x;t,y),$  $(s,x;t,y)\in\mathbb{R}^4_+$ can be coupled in one probability space to form a process on $\mathbb{R}_+^4$ with many desired features. In the following proposition, we collect some of the results in \cite[Theorem 2.2, Proposition 2.3 and Lemma 3.12]{AJRS}. We formulate them in terms of narrow wedge solutions which are more suitable for our purpose. 
\begin{proposition}[\cite{AJRS}]\label{thm:forcoupling}
There exists a coupling of  $$\{\mathcal{H}(s,x;t,y),\ (s,x;t,y)\in\mathbb{R}^4_+ \}$$ with the following properties. 
\begin{enumerate}
\item $\mathcal{H}(t,x;s,y)$ is a random continuous function on $\mathbb{R}^4_+$ .
\item Almost surely for all $(s,x,t,y)\in\mathbb{R}^4_+$ and $r\in (s,t)$, it holds that
\begin{align}\label{equ:KPZ_linear}
\exp\big(\mathcal{H}(s,x;t,y)\big)=\int_{-\infty}^\infty \exp\big( \mathcal{H}(s,x;r,z)+\mathcal{H}(r,z;t,y)\big)\, dz. 
\end{align}
\item For any fixed $(r,u)\in\mathbb{R}^2$, $\mathcal{H}(s+r,x+u;t+r,y+u)\overset{d}{=} \mathcal{H}(s,x;t,y).$ $\mathcal{H}(-t,y;-s,x) \overset{d}{=} \mathcal{H}(s,x;t,y).$ Moreover, for any $\nu\in\mathbb{R}$,
\begin{align*}
\mathcal{H}(s ,x+\nu s ;t ,y+\nu t)\overset{d}{=} \mathcal{H}(s,x;t,y)-\nu (y-x)-2^{-1}\nu^2(t-s).
\end{align*}
  
\item Fix finitely many disjoint open intervals $\{(s_j,t_j)\}_{j=1}^m$. The random functions $ \mathcal{H}(s_j,\cdot ;t_j,\cdot )$ are independent.
\end{enumerate}  
\end{proposition}

Let us define \textbf{KPZ sheets}, which will play an important role in the paper.  For $t=T>s=0$ fixed, the marginal $\mathcal{H}(0,x;T,y)$, viewed as a random continuous function on $\mathbb{R}^2$, is called a KPZ sheet. We denote it by
\begin{equation}\label{def:KPZ_sheet}
\mathbcalboondox{h}^T(x,y)\coloneqq \mathcal{H}(0,x;T,y).
\end{equation}

\subsection{The Airy line ensemble, Airy sheets and the directed landscape} 
In this subsection, we introduce several central objects in the KPZ universality class: the Airy line ensemble, Airy sheets and the directed landscape.
\begin{definition}\label{def:Airy}
The \textbf{stationary Airy line ensemble} $\tilde{\mathcal{A}} = \{\tilde{\mathcal{A}}_1 > \tilde{\mathcal{A}}_2 > \cdots \}$ is a collection of countable many random functions indexed by $\mathbb{N}$. The law of $\tilde{\mathcal{A}}$ is uniquely determined by its determinantal structure. More precisely, for any finite set $I=\{u_1, \cdots, u_k\}\subset\mathbb{R}$, the point process on $I\times\mathbb{R}$ given by $\{(s,\tilde{\mathcal{A}}_i(s)):i\in\mathbb{N},s\in I\}$ is a determinantal point process with kernel given by  
\begin{equation*}
K(s_1,x_1;s_2,x_2)=\left\{ \begin{array}{cc}
\int_0^{\infty} e^{-z(s_1-s_2)}\textup{Ai}(x_1+z)\textup{Ai}(x_2+z)dz & \textrm{if} \quad s_1\geq s_2,\\
-\int_{-\infty}^0 e^{-z(s_1-s_2)}\textup{Ai}(x_1+z)\textup{Ai}(x_2+z)dz  & \textrm{if} \quad s_1< s_2,
\end{array} \right.
\end{equation*}
where $\textup{Ai}$ is the Airy function. The \textbf{parabolic Airy line ensemble} $ {\mathcal{A}}=\{ {\mathcal{A}}_1> {\mathcal{A}}_2>\dots\}$ is defined by
\begin{equation}
 {\mathcal{A}}_i(x):= \tilde{\mathcal{A}}_i( x)-x^2.
\end{equation}
\end{definition}
The finite-dimensional marginal of the stationary Airy line ensemble was introduced by Pr\"{a}hofer and Spohn \cite{PS} in which it was called the ``multi-line Airy process." Later, Corwin and Hammond \cite{CH14} showed that $\mathcal{A}$ can be realized as a random continuous function on $\mathbb{N}\times\mathbb{R}$ through the Brownian Gibbs property. The first indexed random function, $\tilde{\mathcal{A}}_1$, is of particular importance and is known as the $\textup{Airy}_2$ process. 

In the monumental work \cite{DOV}, Dauvergne, Ortmann and Vir\'{a}g constructed Airy sheets and the directed landscape via the parabolic Airy line ensemble. The directed landscape can be viewed as ``fundamental solutions" to the KPZ fixed point and Airy sheets are fixed time marginals of the directed landscape. We follow the presentation in \cite{DOV}  and define Airy sheets and the directed landscape through their characterization properties. We remark that it was proved in \cite{DOV} that those properties uniquely determine the laws of Airy sheets and the directed landscape respectively. 
\begin{definition}\label{def:Airy_sheet}
The \textbf{Airy sheet} $\mathcal{S}(x,y)$ is a $C(\mathbb{R}^2,\mathbb{R})$-valued random variable which can be coupled with the parabolic Airy line ensemble $\mathcal{A}$ with the following properties. 
\begin{enumerate}
\item $\mathcal{S}(\cdot+t ,\cdot+t )$ has the same distribution as $\mathcal{S}(\cdot  ,\cdot  ).
$
\item $S(0,\cdot)=\mathcal{A}_1(\cdot)$. 
\item Almost surely for all $x>0$ and $y_1, y_2$ in $\mathbb{R}$, we have
\begin{equation}\label{equ:D}
\begin{split}
\lim_{k\to\infty} \mathcal{A}[(-2^{-1/2}k^{1/2}x^{-1/2} ,k)\xrightarrow{\infty}  (y_2,1)] -\mathcal{A}[(-2^{-1/2}&k^{1/2}x^{-1/2},k)\xrightarrow{\infty}  (y_1,1)]\\
 &=\mathcal{S}(x,y_2)-\mathcal{S}(x,y_1).
\end{split} 
\end{equation}
\end{enumerate}
Here $\mathcal{A}[(x,k)\xrightarrow{\infty} (y,1)]$ is the last passage time from $(x,k)$ to $(y,1)$ on the parabolic Airy line ensemble. We refer readers to \eqref{def:freeenergy} in Section~\ref{sec:polymer} for the detailed definition. For any $s>0$, $s\mathcal{S}(s^{-2}x,s^{-2}y)$ is called an \textbf{Airy sheet of scale} $\bm{s}$ .  
\end{definition}

\begin{definition}\label{def:landscape}
The \textbf{directed landscape} $\mathcal{L}(s,x;t,y)$ is a $C(\mathbb{R}^4_+,\mathbb{R})$-valued random variable with the following properties.
\begin{enumerate}
\item Given $s<t$, $\mathcal{L}(s,\cdot;t,\cdot)$ is distributed as an Airy sheet of scale $(t-s)^{1/3}$.
\item For any finite disjoint open intervals $\{(s_j,t_j)\}_{j=1}^m$, the random functions $ \mathcal{L}(s_j,\cdot;t_j,\cdot) $ are independent.
\item Almost surely for all $s<r<t$ and $x,y\in\mathbb{R}$, it holds that
\begin{align*}
\mathcal{L}(s,x;t,y)=\max_{z\in\mathbb{R}}\big( \mathcal{L}(s,x;r,z)+\mathcal{L}(r,z;t,y) \big).
\end{align*}
\end{enumerate}

\end{definition}
As a direct consequence of Definition~\ref{def:landscape}, we have the following description on the marginal law of $\mathcal{L}$ when the time variables are restricted on a finite set.
\begin{corollary}\label{cor:landscape_finite}
Fix a finite set $\{t_1<t_2<\dots <t_m\}$. Then the restriction of the directed landscape, $\{\mathcal{L}( t_i,\cdot ; t_j,\cdot )\}$ is uniquely characterized by the following properties.
\begin{enumerate}
\item  For all $1\leq i<j\leq m$, $\mathcal{L}( t_i,\cdot ; t_j,\cdot )$ is distributed as an Airy sheet of scale $ (t_j-t_i)^{1/3}$.
\item  $\{\mathcal{L}( t_i,\cdot ; t_{i+1},\cdot )\}_{i=1}^{m-1}$ are independent.
\item  Almost surely for all $x,y\in\mathbb{R}$ and $1\leq i<j<k\leq m $,
\begin{align*}
\mathcal{L}( t_i,x ; t_k,y )=\max_{z\in\mathbb{R}}\left( \mathcal{L}( t_i,x ; t_j,z )+\mathcal{L}( t_j,z ; t_k,y ) \right).
\end{align*}
\end{enumerate}   
\end{corollary}

\subsection{Main results}\label{sec:mainresult}

In this subsection, we perform the $1:2:3$ scaling to $\mathcal{H}(s,x;t,y)$ and state our main results concerning the large-time asymptotics. For $T>0$, the \textbf{scaled narrow wedge solutions} are given by
\begin{align}\label{def:KPZ_4_scaled}
\mathfrak{H}^T( s,x ; t,y ):=2^{1/3}T^{-1/3}\mathcal{H}( Ts,2^{1/3}T^{2/3}x  ; Tt, 2^{1/3}T^{2/3}y )+(t-s)2^{1/3}T^{ 2/3}/24.
\end{align}
For $t=1$ and $s=0$ fixed, we call the marginal $\mathfrak{H}^T( 0,x ; 1,y )$ the \textbf{scaled KPZ sheet} and denote it by
\begin{equation}\label{def:KPZ_sheet_rescaled0}
\mathfrak{h}^T(x,y)\coloneqq \mathfrak{H}^T( 0,x ; 1,y ).
\end{equation}
Note that from \eqref{def:KPZ_sheet}, $\mathfrak{h}^T $ can be expressed in terms of the KPZ sheet $\mathbcalboondox{h}^T$ as
\begin{equation}\label{def:KPZ_sheet_rescaled}
\mathfrak{h}^T(x,y)\coloneqq 2^{1/3}T^{-1/3}\mathbcalboondox{h}^T(2^{1/3}T^{2/3}x,2^{1/3}T^{2/3}y)+ {2^{1/3}T^{2/3}}/{24}.
\end{equation}

It is conjectured that the KPZ fixed point describes the large-time behavior of solutions to the KPZ equation \eqref{eq:KPZ}. In \cite{ACQ}, Amir, Corwin and Quastel gave strong evidence for this conjecture and proved that $\mathfrak{H}^T(0,0;1,0)$ converges to the Tracy-Widom law. Equivalently, using the notation introduced above, they showed that $\mathfrak{H}^T(0,0;1,0)$ converges in distribution to $\mathcal{L}(0,0;1,0)$. See contemporary physics work in \cite{Calabrese_2010,Dotsenko_2010,SS}. Recently, a breakthrough was made by two groups, Quastel-Sarkar \cite{QS} and Vir\'{a}g \cite{V}. The authors independently proved that $\mathfrak{H}^T(0,0;t,y)$, as a random function on $\mathbb{R}$, converges in distribution to $\mathcal{L}(0,0;t,y)$. 

In this paper, we establish the convergence of $\mathfrak{H}^T(s,x;t,y)$ to $\mathcal{L}(s,x;t,y)$ as a four-parameter process in the locally uniform topology. This allows us to provide an independent proof for the convergence of the KPZ equation to the KPZ fixed point for general initial conditions. We also establish the joint convergence to the KPZ fixed point in the locally uniform topology for multiple initial conditions.

Before stating our main results, we note that for a topological space $\mathcal{T}$, we always equip $C(\mathcal{T},\mathbb{R})$, the collection of continuous functions on $\mathcal{T}$, with the topology of uniform convergence on compact subsets.

\begin{theorem}\label{thm:KPZtoLandscape}
The scaled narrow wedge solutions $ \mathfrak{H}^T$ converge in distribution to the directed landscape $ \mathcal{L} $ as $T$ goes to infinity. Here $\mathfrak{H}^T $ and $\mathcal{L} $ are viewed as $C(\mathbb{R}^4_+,\mathbb{R})$-valued random variables.
\end{theorem}


{\color{black}The directed landscape convergence has been proved through line ensembles for integrable last passage percolation models \cite{DOV, dauvergne2022} in the locally uniform topology. \cite{aggarwal2024} proved (finite dimensional in time and locally uniform in space) convergence to the directed landscape for the colored ASEP, and the stochastic six-vertex models. Recently, \cite{dauvergne2025} provided a different approach to proving the convergence to the directed landscape in the sense of finite dimensional distributions both for time and space based on a new characterization of the directed landscape.}

As a crucial middle step of proving Theorem \ref{thm:KPZtoLandscape}, we show the convergence of the scaled KPZ sheet to the Airy sheet.

\begin{theorem}\label{thm:KPZtoAiry_sheet}
The scaled KPZ sheets $\mathfrak{h}^T $ converge in distribution to the Airy sheet $\mathcal{S}$ as $T$ goes to infinity. Here $\mathfrak{h}^T$ and $\mathcal{S}$ are viewed as $C(\mathbb{R}^2,\mathbb{R})$-valued random variables.  
\end{theorem}

Up until now, we have focused on narrow wedge solutions. In the following, we examine solutions to the KPZ equation with general initial functions. For a given initial function, Quastel-Sarkar \cite{QS} and Vir\'{a}g \cite{V} have shown that such solutions converge to the KPZ fixed point. We will give an independent proof for the convergence of KPZ equation to the KPZ fixed point in the locally uniform topology based on the convergence of the KPZ equation to the directed landscape.

The KPZ fixed point is a Markov process and it is believed to govern the large-time asymptotics of models in the KPZ universality class, including the KPZ equation \eqref{eq:KPZ}. The KPZ fixed point was first constructed in \cite{MQR} and the authors provided the transition probability through the Fredholm determinant. In this paper, we focus on the following variational description proven in \cite[Corollary 4.2]{MR4190063}. For a continuous function $f(x)$ that satisfies $f(x)\leq C(1+ |x| )$ for some $C$, the KPZ fixed point with initial function $f(x)$ is a random process on $(t,y)\in (0,\infty)\times\mathbb{R}$ given by
\begin{align}\label{equ:infinity_convolution}
(f\, \overline{\otimes}\,  \mathcal{L})(t,y):=\sup_{x\in\mathbb{R}}\big( f(x)+\mathcal{L}(0,x;t,y) \big).
\end{align}
For $T>0$, define
\begin{align}\label{equ:T_convolution}
(f\otimes_T\mathfrak{H}^T)(t,y)\coloneqq 2^{1/3}T^{-1/3}\log\int\exp\left( 2^{-1/3}T^{1/3}\left( f(x)+\mathfrak{H}^T(0,x;t,y) \right) \right) dx.
\end{align}
It is straightforward to check that $(f\otimes_T\mathfrak{H}^T)(t,y)$ is the 1:2:3 scaled solution to the KPZ equation \eqref{eq:KPZ} with initial data $2^{-1/3}T^{1/3}f(2^{-1/3}T^{-2/3}x)$ at time $0$.

\begin{theorem}\label{thm:single_sol}
 $f\otimes_T\mathfrak{H}^T$ converges in distribution to $f\,\overline{\otimes}\,\mathcal{L}$ as $T$ goes to infinity. Here $f\otimes_T\mathfrak{H}^T$ and $f\,\overline{\otimes}\,\mathcal{L}$ are viewed as $C((0,\infty)\times\mathbb{R} ,\mathbb{R})$-valued random variables.
\end{theorem}

With Theorem~\ref{thm:KPZtoLandscape}, we can extend the convergence result in Theorem~\ref{thm:single_sol} to multiple initial conditions. Generalizing \eqref{equ:infinity_convolution} and \eqref{equ:T_convolution}, we define
\begin{align}\label{equ:infinity_convolution_s}
(f\, \overline{\otimes}\,  \mathcal{L})(s;t,y)\coloneqq & \sup_{x\in\mathbb{R}}\big( f(x)+\mathcal{L}(s,x;t,y) \big),
\end{align}
\begin{align}\label{equ:T_convolution_s}
(f\otimes_T\mathfrak{H}^T)(s;t,y)\coloneqq &2^{1/3}T^{-1/3}\log\int\exp\left( 2^{-1/3}T^{1/3}\left( f(x)+\mathfrak{H}^T(s,x;t,y) \right) \right) dx.
\end{align} 
It can be verified that $\left( f \otimes_T \mathfrak{H}^T \right)(s ;t,y)$ is the 1:2:3 scaled solution to the KPZ equation~\ref{eq:KPZ} with initial condition $2^{-1/3}T^{1/3}f (2^{-1/3}T^{-2/3}x)$ at time $Ts$. Fix $s_1,\dots , s_N\in\mathbb{R}$ and functions $f_1(x),\dots, f_N(x)$ that satisfy $f_i(x)\leq C(1+|x|)$ for some $C$, we have the following joint convergence. 

\begin{corollary}\label{cor:multi_sol}
$\left( f_i\otimes_T \mathfrak{H}^T \right)(s_i;t,y),\ 1\leq i\leq N$ jointly converge in distribution to $\left( f_i\, \overline{\otimes}\, \mathcal{L}  \right)(s_i;t,y)$,\ $1\leq i\leq N$. Here $\left( f_i\otimes_T \mathfrak{H}^T \right)(s_i;t,y) $ and 
$\left( f_i\, \overline{\otimes}\, \mathcal{L}  \right)(s_i;t,y)$ are viewed as $C((s_i,\infty)\times\mathbb{R},\mathbb{R})$-valued random variables. 
\end{corollary}

Building on Theorems~\ref{thm:KPZtoLandscape} and \ref{thm:KPZtoAiry_sheet}, one can establish further connections between objects related to the KPZ equation and the KPZ fixed point. For instance, in \cite[Theorem 1.6]{DZ}, Das and Zhu proved the one-point convergence of continuum directed random polymer paths to a geodesic in the directed landscape. Using Theorem~\ref{thm:KPZtoAiry_sheet}, the authors were able to upgrade this to a process-level convergence \cite[Theorem 1.7]{DZ}. 

\subsection{O'Connell Yor polymers and the KPZ line ensemble}  
For $x< y$ and natural numbers $\ell \geq m$, we denote by $\mathcal{Q}[(x,\ell )\to (y,m)]$ the collection of directed semi-discrete paths from $(x,\ell)$ to $(y,m)$. By considering the jump times, $\mathcal{Q}[(x,\ell )\to (y,m)]$ can be identified with the convex set $\{x\leq t_\ell\leq t_{\ell-1}\leq \dots \leq t_{m+1}\leq y\}$ in $\mathbb{R}^{\ell-m}$. We use $\pi$ to denote a path and write $d\pi$ for the Lebesgue measure on $\mathcal{Q}[(x,\ell )\to (y,m)]$. Let $f=\{f_1,f_2,\dots, f_n\}$ be $n$ continuous functions with $n\geq \ell$. For a path $\pi$ from $(x,\ell)$ to $(y,m)$ with jump times $t_\ell\leq t_{\ell-1}\leq \dots \leq t_{m+1}$, we assign the weight  $f(\pi)\coloneqq \sum_{j=m}^{\ell} f_j(t_j)-f_{j}(t_{j+1})$ with convention that $t_{\ell+1}=x$ and $t_{m}=y$. For $\beta>0$, the probability measure on paths proportional to $\exp\left(\beta {f(\pi)}\right)d\pi$ is called the $\bm{\beta}$\textbf{-polymer measure}. When $\beta$ goes to infinity, the $\beta$-polymer measures concentrate on paths that achieve the maximum of $f(\pi)$. Those maximum paths are called \textbf{geodesics}. The $\bm{\beta}$\textbf{-free energy} from $(x,\ell)$ to $(y,m)$ is given by
\begin{align*} 
f[(x,\ell)\xrightarrow{\beta} (y,m)]\coloneqq \beta^{-1} \log\int_{\mathcal{Q}[(x,\ell)\to (y,m)]} \exp\left(\beta {f(\pi)}\right)\, d\pi.
\end{align*}
When $\beta$ goes to infinity, we obtain the \textbf{last passage time} as
\begin{align*}
 f[(x,\ell)\xrightarrow{\infty} (y,m)]\coloneqq \max_{\pi\in \mathcal{Q}[(x,\ell)\to (y,m)]} f(\pi).
\end{align*}
We will mostly work with the case $\beta=1$. Hence we simply write polymer measure and $f[(x,\ell)\to (y,m)]$ for $1$-polymer measure and $f[(x,\ell)\xrightarrow{1} (y,m)]$ respectively.

Next, we discuss the O'Connell-Yor polymer model \cite{OY}. Let $B=\{B_1,B_2,\dots\}$ be a collection of i.i.d. standard two-sided Brownian motions. The $n$-th O'Connell-Yor free energy is a random continuous function on $(0,\infty)$ defined by $Y^n_1(x)\coloneqq B[(0,n)\to (x,1)].$ Furthermore, O'Connell \cite{OC12} demonstrated that $Y^n_1$ can be embedded into a line ensemble $Y^n=(Y^n_1,\dots,Y^n_n)$, where $Y^n$ is the collection of $n$ random continuous functions obtained by taking geometric RSK transform to Brownian motions.

The KPZ line ensemble was constructed by Corwin and Hammond \cite{CH16} as a scaling limit of $Y^n$ as follows. For $n,i\in\mathbb{N}$ and $T>0$, let
\begin{align*}
 C_1(T,n)\coloneqq & n^{1/2}T^{-1/2}+2^{-1},\ C_{3,i}(T)\coloneqq -(i-1)\log T+\log(i-1)!, \\
C_2(T,n)\coloneqq & n+2^{-1}n^{1/2}T^{1/2}- (n-1)\log( n^{1/2}T^{-1/2}) .
\end{align*}
Define $\mathcal{X}^{T,n}=\{\mathcal{X}^{T,n}_1,\mathcal{X}^{T,n}_2,\dots, \mathcal{X}^{T,n}_n\}$ by 
\begin{equation*}
\mathcal{X}^{T,n}_i(x)\coloneqq Y^n_i(n^{1/2}T^{1/2}+x)-C_1(T,n)x-C_2(T,n)-C_{3,i}(T) ,\ i\in\llbracket 1,n \rrbracket
\end{equation*}
It was proved in \cite[Theorem 3.10]{CH16} that $\{\mathcal{X}^{T,n}\}_{n\in\mathbb{N}}$ is tight. Later it was shown by Nica \cite[Corollary 1.7]{Nic21} that any subsequential distributional limit of $\mathcal{X}^{T,n}$ has the same law and hence $\mathcal{X}^{T,n}$ converges in distribution. The limiting $C(\mathbb{N}\times\mathbb{R},\mathbb{R})$-valued random variable is called the \textbf{KPZ line ensemble} and we denote it by $\mathcal{X}^T$.
\begin{proposition}[\cite{CH16,Nic21}] \label{pro:kpzlineensemble}
Fix $T>0$. When $n$ goes to infinity, $\mathcal{X}^{T,n}$ converges in distribution to the $\textup{KPZ}$ line ensemble $\mathcal{X}^T$. Here $\mathcal{X}^{T,n}$ and $\mathcal{X}^{T}$ are viewed as $C(\mathbb{N}\times\mathbb{R},\mathbb{R})$-valued random variables.
\end{proposition}
We note that O'Connell-Warren \cite{OW} constructed the exponential of the KPZ line ensemble from a different perspective, which the authors called a multi-layer extension of SHE.

The first indexed function in the KPZ line ensemble has the same distribution as the narrow wedge solutions. That is, {\color{black}$\mathcal{X}^T_1(\cdot)\overset{d}{=}\mathcal{H}(0,0;T,\cdot){=}\mathbcalboondox{h}^T(0,\cdot )$}. We further consider two-variable random functions which converge to $\mathbcalboondox{h}^T(\cdot ,\cdot)$. 
For $T>0$, $n\geq 1$, $x\in\mathbb{R} $ and $y>-n^{1/2}T^{1/2}+x$, define
\begin{align*}
\mathbcalboondox{h}^{T,n}(x,y)\coloneqq B[(x,n)\to (n^{1/2}T^{1/2}+y,1)]-C_1(T,n)(y-x)-C_2(T,n).
\end{align*}
The finite-dimensional converges of $\mathbcalboondox{h}^{T,n}(x,y)$ to $\mathbcalboondox{h}^{T}(x,y)$ was essentially proved in \cite{Nic21}. In \cite[Theorem 1.2]{Nic21}, the author proved the finite-dimensional convergence for $\mathbcalboondox{h}^{T,n}(0,\cdot)$ and the same argument applies for $\mathbcalboondox{h}^{T,n}(\cdot ,\cdot)$.  See \cite[Section 6.2]{AKQ} for a similar result for discrete polymers.
\begin{proposition} \label{pro:kpzsheet}
Fix $T>0$. When $n$ goes to infinity, the finite-dimensional marginal of $\mathbcalboondox{h}^{T,n}$ converges in distribution to the finite-dimensional marginal of the $\textup{KPZ}$ sheet $\mathbcalboondox{h}^T$. 
\end{proposition}

By utilizing Propositions~\ref{pro:kpzlineensemble} and \ref{pro:kpzsheet}, we can naturally couple the KPZ sheet and the KPZ line ensemble in one probability space. Furthermore, we will demonstrate that the Busemann functions of the KPZ sheet are linked to the polymer free energies of the KPZ line ensemble. We conjecture that these two entities are connected via the closed formula below.

\begin{conjecture}\label{conj:main}
Fix $T>0$. There exists a coupling of the KPZ line ensemble $\mathcal{X}^T(y)$ and the KPZ sheet $\mathbcalboondox{h}^T(x,y)$ such that almost surely the following holds.
\begin{enumerate}
\item $ \mathbcalboondox{h}^T(0,\cdot)= \mathcal{X}^T_1(\cdot)$.
\item For all $x>0$ and $y_1, y_2$ in $\mathbb{R}$, we have
\begin{equation}\label{equ:main}
\begin{split}
\lim_{k\to\infty}\bigg( \mathcal{X}^T[(-kT/x,k)\to (y_2,1)]-\mathcal{X}^T[(-kT/x,&k)\to (y_1,1)]\bigg)\\
&= \mathbcalboondox{h}^T(x,y_2)-\mathbcalboondox{h}^T(x,y_1).
\end{split}
\end{equation}
\end{enumerate}
\end{conjecture}

Even though we do not have a proof for Conjecture~\ref{conj:main}, we are able to reduce it to certain property about the KPZ line ensemble. 

\begin{theorem}\label{thm:main}
Suppose for any $\varepsilon>0$ and $x>0$, it holds that
\begin{align}\label{equ:wishhhh}
\sum_{k=1}^\infty\mathbb{P}\bigg( \left| \mathcal{X}^T[(0,k+1)\to (x,1)]-k\log x+\log k! \right|>\varepsilon k\bigg)<\infty.
\end{align}
Then Conjecture~\ref{conj:main} holds true.
\end{theorem}

{\color{black}After completing this paper, the author is able to prove a modulus of continuity estimate for KPZ line ensembles \cite{wu2025} which implies \eqref{equ:wishhhh} and confirms Conjecture~\ref{conj:main}. We present the proof of \eqref{equ:wishhhh} in Section~\ref{sec:proof2}.
}
\subsection{Ingredients and ideas} 
\begin{center}
\begin{figure}[ht]
\includegraphics[width=9cm]{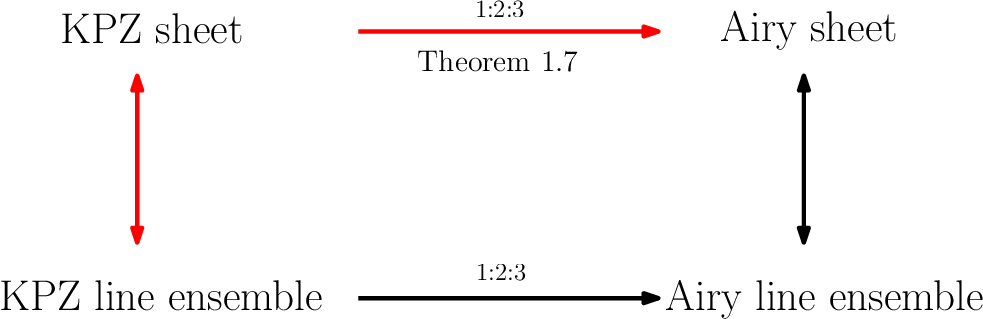}
\caption{The right and the bottom arrows are two major inputs. The top arrow is one of our main result, Theorem~\ref{thm:KPZtoAiry_sheet}. We prove it through establishing the left arrow, a relation between the KPZ sheet and the KPZ line ensemble.}\label{fig:plan}
\end{figure}
\end{center}

{\color{black}This subsection outlines the major inputs and key ideas behind the proof of our main results, Theorems~\ref{thm:KPZtoLandscape} and \ref{thm:KPZtoAiry_sheet}}. One of the inputs is the characterization of the Airy sheet, recorded as Definition~\ref{def:Airy_sheet}, proved in \cite{DOV}. The authors established a canonical coupling of the Airy sheet and the parabolic Airy line ensemble. In this coupling, Busemann functions of the Airy sheet are encoded in the parabolic Airy line ensemble with a closed formula \eqref{equ:D}. The second important ingredient is the convergence of the KPZ line ensembles to the parabolic Airy line ensemble. This convergence of line ensembles was proved by combining results from a series of works \cite{QS,V,DM,Wu21tightness,aggarwal2023}. With these two ingredients in hand, the missing link between the KPZ sheet and the Airy sheet is a connection between the KPZ sheet and the KPZ line ensemble. See Figure~\ref{fig:plan}. This connection is a key intermediate step toward the proof of main theorems. 

The relation between the KPZ sheet and the KPZ line ensemble originates in a geometric RSK invariance. The geometric RSK correspondence for continuous functions was introduced in \cite{OC12}. It maps $n$ continuous functions $f=\{f_1,f_2,\dots, f_n\}$ to another $n$ continuous functions $Wf=\{(Wf)_1,(Wf)_2,\dots, (Wf)_n\}$. {\color{black}It was proved by Noumi and Yamada \cite{MR2074600}, and more recently by Corwin \cite{Cor21}, that polymer free energies are invariant under the gRSK transformation.} Namely, for $x<y$ there holds
\begin{align}\label{equ:intro_gRSK_inv}
f[(x,n)\to (y,1)]= (Wf)[(x,n)\to (y,1)].
\end{align}
In the continuum directed random polymer model , the KPZ sheet represents polymer free energies in a white noise background. Furthermore, the KPZ line ensemble may be viewed as the output of the gRSK transform with a white noise input \cite{OW}. Therefore, \eqref{equ:intro_gRSK_inv} strongly suggests there should be a connection between the KPZ sheet and the KPZ line ensemble. Concretely, in the context of the O'Connell-Yor model, one has
\begin{align}\label{equ:intro_OY_inv} 
\mathbcalboondox{h}^{T,n} (x,y)=\mathcal{X}^{T,n}[(-n^{1/2}T^{1/2}+x ,n)\to (y,1)] )+ (n-1)\log (n^{1/2}T^{-1/2}).  
\end{align}
Considering Propositions~\ref{pro:kpzlineensemble} and \ref{pro:kpzsheet}, it can be concluded that $\mathbcalboondox{h}^{T,n}$ and $\mathcal{X}^{T,n}$ converge to the KPZ sheet and the KPZ line ensemble, respectively. One may attempt to send $n$ to infinity in \eqref{equ:intro_OY_inv}. However, a significant difficulty arises in doing so because the right-hand side of the equation involves information about $\mathcal{X}^{T,n}$ in a rapidly enlarging region $[-n^{1/2}T^{1/2}+x,y]\times \{1,2,\dots, n\}$. An emergent question is whether it is possible, and if so, how to take a meaningful limit of \eqref{equ:intro_OY_inv}. Our main contribution is to provide a (partial) affirmative answer to this question. Our resolution is inspired by \cite{DOV} in which a similar issue arose in the context of Brownian last passage percolation. However, due to the nature of the polymer model, several new ideas are necessary, which we will discuss in the following paragraphs.

Let us perform simplification and introduce notation. We fix an environment. Namely, we fix countable many continuous functions $\mathcal{X}=\{\mathcal{X}_1,\mathcal{X}_2,\dots\}$ and consider 
\begin{align}\label{equ:Xnenergy}
\mathbcalboondox{h}(x,y)= \mathcal{X}[(-n^{1/2}T^{1/2}+x ,n)\to (y,1)].
\end{align} 
Here we suppress the dependence in $n$ and $T$ for brevity. Let $G_k(z,y)$ and $F_k(x,z)$ be polymer energies from $(z,k)$ to $(y,1)$ and from $(-n^{1/2}T^{1/2}+x,n)$ to $(z,k+1)$ respectively.
\begin{align*}
G_k(z,y)=\mathcal{X}[(z ,k)\to (y,1)],\ F_k(x,z)=\mathcal{X}[(-n^{1/2}T^{1/2}+x ,n)\to (z,k+1)].
\end{align*}
It holds that (see Lemma~\ref{lem:integral})
{\color{black}
\begin{align}\label{equ:k-level}
\exp\big(\mathbcalboondox{h}(x,y)\big)=\int_{\{z< y\}} \exp\big( F_k(x,z)+G_k(z,y) \big)\, dz.
\end{align}
}
\vspace*{0.3cm}

We now discuss the key ideas and observations in our approach. Inspired by equation \eqref{equ:D}, we aim to compare the values of $\mathbcalboondox{h}(x,y_2)-\mathbcalboondox{h}(x,y_1)$ and $G_k(z,y_2)-G_k(z,y_1)$. In the case of last passage percolation, especially when geodesics merge, one can select $z$ such that $\mathbcalboondox{h}(x,y_2)-\mathbcalboondox{h}(x,y_1)$ and $G_k(z,y_2)-G_k(z,y_1)$ are equal. This equality could also be achieved for polymer models. However, determining the location of $z$ that produces this equality is difficult, as it lacks geometric significance. Instead, we make an observation: for any $z$, we can provide a bound for the difference between $\mathbcalboondox{h}(x,y_2)-\mathbcalboondox{h}(x,y_1)$ and $G_k(z,y_2)-G_k(z,y_1)$ using polymer measures. The $k$-th level marginal of the polymer measure is given by 
{\color{black}
\begin{align*}
\mu_k(x,y;dz)=\mathbbm{1}(z< y)\exp(-\mathbcalboondox{h}(x,y)+F_k(x,z)+G_k(z,y))dz.
\end{align*}
}
In Lemma~\ref{lem:GG}, we show that the difference between $\mathbcalboondox{h}(x,y_2)-\mathbcalboondox{h}(x,y_1)$ and $G_k(z,y_2)-G_k(z,y_1)$ can be controlled using the cumulative distribution functions of $\mu_k$. This is a consequence of the monotonicity of $G_k(z,y_2)-G_k(z,y_1)$, as shown in Lemma~\ref{lem:Fconvex}.

Moving forward, our focus is on the $k$-th level marginal of the polymer measure, $\mu_k(x,y;dz)$. Our crucial observation is that $F_k(x,z)$ plays a critical role in determining $\mu_k(x,y;dz)$. This is because $F_k(x,z)$ is much more sensitive to changes in $x$ compared to $\mathbcalboondox{h}(x,y)$. To illustrate this point, let us consider a simplified scenario where $\mathbcalboondox{h}(x,y)$ does not depend on $x$ at all. In this case, we claim that $\mu(x,y;dz)$ is a delta mass centered at a point $z_0$. Furthermore, $z_0$ is determined by the condition $(\partial F_k/\partial x)(x,z_0)=0$. To obtain this result, we differentiate the logarithm of \eqref{equ:k-level} with respect to $x$ and get 	
\begin{align}\label{equ:partialF}
0=\int \frac{\partial F_k}{\partial x}(x,z) \mu_k(x,y;dz).
\end{align}
Differentiating it one more time with respect to $y$, we get
\begin{align}\label{equ:partialFG}
\begin{split}
0=&\int \frac{\partial F_k}{\partial x}(x,z)\frac{\partial G_k}{\partial y}(z,y) \mu_k(x,y;dz)-\int \frac{\partial F_k}{\partial x}(x,z) \mu_k(x,y;dz)\int  \frac{\partial G_k}{\partial y}(z,y)\mu_k(x,y;dz)\\=&  2^{-1}\int\int \left( \frac{\partial F_k}{\partial x}(x,z)-\frac{\partial F_k}{\partial x}(x,z') \right)\\
&\qquad\qquad\qquad\times\left( \frac{\partial G_k}{\partial y}(z,y)-\frac{\partial G_k}{\partial y}(z',y) \right)\mu_k(x,y;dz)\mu_k(x,y;dz'). 
\end{split}
\end{align}
Note that, as shown in Lemma~\ref{lem:Fconvex}, both $ \partial F_k / \partial x $ and $ \partial G_k/\partial y $ are monotone non-decreasing in $z$. {\color{black} Assume $ \partial F_k / \partial x $ and $ \partial G_k/\partial y $ are strictly increasing in $z$.} Then \eqref{equ:partialFG} implies that $\mu_k(x,y;dz)$ is a delta measure centered at some $z_0$. From \eqref{equ:partialF}, we see that $z_0$ satisfies $(\partial F_k/\partial x)(x,z_0)=0$. In practice, the assumption that $\mathbcalboondox{h}(x,y)$ is independent of $x$ does not hold. Nevertheless, we can still use $F_{k}(x,z)$ to bound the polymer measure $\mu_k(x,y;dz)$. This is the content of Lemmas~\ref{lem:boundzz} and \ref{lem:boundz}. Note that these lemmas hold without any assumption on $F_k(x,z)$, but they are particularly useful when $F_k(x,z)$ changes significantly with $x$ for large $k$. 

We now turn to the analysis of the asymptotics of $F_k(x,z)$, which is the polymer energy from $(-n^{1/2}T^{1/2}+x, n)$ to $(z,k+1)$. It is worth noting that in this scenario, we are dealing with two parameters, $n$ and $k$, that are sent to infinity. Using a special identity \eqref{equ:Wf+WRf}, we can compute the exact distributional limit as $n$ goes to infinity. In the next paragraph, we will delve further into the details of \eqref{equ:Wf+WRf}. {\color{black}The resulting limit corresponds to the polymer energy from $(0,k+1)$ to $(x,1)$ on the {KPZ line ensemble}. This is the content of Lemma~\ref{lem:Flimit}.} Upon applying the 1:2:3 scaling, the polymer energy is transformed into the {last passage time} on the {parabolic Airy line ensemble}. The large $k$ limit in this setting has been previously studied in \cite{DOV}. In the final paragraph, we will discuss the large $k$ asymptotics without the 1:2:3 scaling.

The derivation of \eqref{equ:Wf+WRf} is inspired by \cite[Lemma 5.3]{DOV}, a similar identity in the context of the last passage percolation and RSK correspondence. However, compared to last passages, dealing with polymers is more difficult due to the fact that a polymer measure is supported on a collection of paths rather than a unique geodesic. We overcome this difficulty by exploiting a concentration phenomenon of polymer measures when the environment comes from a gRSK transform. Specifically, we consider $n$ continuous functions $f={f_1,f_2,\dots, f_n}$ and let $\mathcal{W}f$ be the gRSK transform of $f$. We show that the polymer measure on $\mathcal{Q}[(0,n)\to (y,1)]$ induced by $\mathcal{W}f$ exhibits concentration on a single geodesic path $\pi$, which is explicit and stays on the first indexed curve. It is worth noting that $\mathcal{W}f$ is not well-defined at $0$, and the polymer measure is understood through a limiting process. Consider the polymer measure on $\mathcal{Q}[(\varepsilon ,n)\to (y,1)]$ induced by $\mathcal{W}f$ for small $\varepsilon>0$. We show that as $\varepsilon$ goes to zero, this family of polymer measures becomes more and more concentrated on the first indexed curve. This concentration phenomenon makes the polymer measure much more tractable and is the starting point for proving \eqref{equ:Wf+WRf}. We establish this phenomenon using Greene's theorem and the invariance of gRSK.

Finally, we discuss the large $k$ asymptotics of the polymer energies from $(0,k+1)$ to $(x,1)$ on the KPZ line ensemble, denoted by $\mathcal{X}^T[(0,k+1)\to (x,1)]$. While we do not have a definitive answer and further investigation is necessary, we expect \eqref{equ:wishhhh} to hold based on the $\mathbf{H}$-Brownian Gibbs property of the KPZ line ensemble. The $\mathbf{H}$-Brownian Gibbs property of the KPZ line ensemble was established in \cite{CH16} and it suggests that the KPZ line ensemble should behave similarly to Brownian motions. It is well-known that the last passage time in a Brownian environment has the leading term $B[(0,k+1)\xrightarrow{\infty}  (x,0)]=2(k+1)^{1/2}x^{1/2}+O(k^{-1/6}).$ Since the volume of $\mathcal{Q}[(0,k+1) \to (x,0)]$ is $x^k/k!$, we obtain $B[(0,k+1)\to (x,0)]=k\log x-\log k!+O(k^{1/2}).$ Note that the leading term in $B[(0,k+1)\to (x,0)]$ comes from the volume of $\mathcal{Q}[(0,k+1) \to (x,0)]$. We believe the same scenario occurs for the KPZ line ensemble.


\subsection{Outline} Section~\ref{sec:polymer} contains the definition of semi-discrete polymer and some of its basic deterministic properties. In Section \ref{sec:gRSK}, we derive a crucial identity related to the geometric RSK correspondence, \eqref{equ:Wf+WRf} in Proposition~\ref{pro:Wf+WRf}. In Section \ref{sec:OY}, we introduce objects in the O’Connell-Yor polymer model and their scaled versions. We prove Theorems~\ref{thm:KPZtoAiry_sheet} and \ref{thm:KPZtoLandscape} in Sections~\ref{sec:proof1} and \ref{sec:proofnew} respectively. In Section~\ref{sec:genini}, we give an independent proof of the convergence of the KPZ equation to the KPZ fixed point in the locally uniform topology (Theorem~\ref{thm:single_sol}) and a proof for the joint convergence for multiple initial conditions (Corollary~\ref{cor:multi_sol}). In Section~\ref{sec:proof2}, we prove Theorem~\ref{thm:main} and Equation~\eqref{equ:wishhhh}, and confirm Conjecture~\ref{conj:main}. Appendix contains proofs for some results used in the paper.

\subsection{Notation}
We would like to explain some notation here. The natural numbers are defined to be $\N = \{1, 2, . . .\}$ and $\N_0 =\N\cup\{0\} $. The positive rational numbers are denoted by $\mathbb{Q}^+$. For $m\leq \ell\in\mathbb{N}_0$, we write $\llbracket m,\ell \rrbracket$ for $\{m,m+1,m+2,\dots,\ell\}$.  We use a special font $\mathsf{t}$ to denote a sequence of positive numbers $\{T_1<T_2<\dots\}$ which goes to infinity. We also denote by $\mathsf{n}$ a sequence of positive integers $\{n_1<n_2<\dots\}$ which goes to infinity.

{\color{black} We use calligraphic fonts to denote objects related to the KPZ equation without scaling.
\begin{itemize}
\item $\mathcal{H}(s,x;t,y)$ : narrow wedge solutions \eqref{equ:0309} 
\item $\mathbcalboondox{h}^T(x,y)$ : KPZ sheets \eqref{def:KPZ_sheet}
\item $\mathcal{X}$ : KPZ line ensembles
\end{itemize}

We use Fraktur fonts to denote the 1:2:3 scaled objects.

\begin{itemize}
\item $\mathfrak{H}^T(s,x;t,y)$ : scaled narrow wedge solutions \eqref{def:KPZ_4_scaled} 
\item $\mathfrak{h}^T(x,y)$ : scaled KPZ sheets \eqref{def:KPZ_sheet_rescaled0} and \eqref{def:KPZ_sheet_rescaled}
\item $\mathfrak{X}$ : scaled KPZ line ensembles \eqref{def:0310}
\end{itemize}
}
$\mathbcalboondox{h}^{T,n}$, $\mathcal{X}^{T,n}$, $\mathfrak{h}^{T,n}$, and $\mathfrak{X}^{T,n}$ are objects in the O'Connell-Yor model which converge to $\mathbcalboondox{h}^{T}$, $\mathcal{X}^{T}$, $\mathfrak{h}^{T}$, and $\mathfrak{X}^{T}$, respectively, as $n$ goes to infinity.

For a topological space $\mathcal{T}$, we equip $C(\mathcal{T},\mathbb{R})$, the collection of continuous functions on $\mathcal{T}$, with the topology of uniform convergence on compact subsets. A family of $C(\mathcal{T},\mathbb{R})$-valued random variables converges in distribution if the corresponding family of measures on $C(\mathcal{T},\mathbb{R})$ converges weakly.

\subsection{Acknowledgement}
The author thanks Promit Ghosal and Alan Hammond for the discussion related to the temporal correlation of the KPZ equation. The author thanks Ivan Corwin for advice on a draft of this paper. The authors thanks anonymous reviewers for their very careful reading of and valuable comments on the manuscript.
The author is partially supported by the NSF through NSF-2348188 and by the Simons Foundation through MPS-TSM-00007939.

\section{Semi-Discrete Polymers}\label{sec:polymer}
In this section, we introduce semi-discrete polymers with a deterministic environment and record some basic properties. The proofs for those properties can be found in Appendix. 
\vspace*{0.3cm}

A semi-discrete environment is given by finitely or countably many continuous functions defined on an interval.
\begin{definition}\label{def:CT}
For $n\in\mathbb{N}$ and an interval $I\subset\mathbb{R}$, we define
\begin{align*}
C^n(I)\coloneqq \left\{  (f_1 ,f_2 ,\dots, f_n )\, |\, f_i \in C(I,\mathbb{R})\ \textup{for}\ i\in\llbracket 1,n \rrbracket \right\}.
\end{align*}
For the special cases $I=(0,T)$ or $I=[0,T)$ for some $0<T$, we denote
\begin{align*}
C^n(T)\coloneqq C^n((0,T))\ \textup{and}\  \overline{C}^n(T)\coloneqq C^n([0,T)).
\end{align*}
\end{definition}
We define the up/right paths connecting two points as follows.
\begin{definition} 
For real numbers $x< y$ and positive integers $\ell \geq m$, we denote by $\mathcal{Q}[(x,\ell )\to (y,m)]$ the collection of non-increasing c\`adl\`ag functions $\pi:[x,y]\to\mathbb{N}$ with $\pi(x)\leq \ell$ and $\pi(y)=m$. We refer to a member of $\mathcal{Q}[(x,\ell )\to (y,m)]$ as a \textbf{path} from $(x,\ell)$ to $(y,m)$.
\end{definition}

There is an injective map from $\mathcal{Q}[(x,\ell )\to (y,m)]$ to $\mathbb{R}^{\ell-m}$, $\pi\mapsto ( t_\ell,\dots ,t_{m+1})$, given by
\begin{align}\label{equ:identify}
t_j= \inf\{ t\in [x,y]\, |\, \pi(t)\leq j-1 \}\ \textup{for}\ j\in \llbracket m+1,\ell \rrbracket. 
\end{align}
It is convenient to set
\begin{align}\label{equ:identifybd}
t_{\ell+1}=x,\ t_m=y. 
\end{align}
In particular, it holds that
\begin{align}\label{equ:pivalue}
\pi(t)=j\ \textup{for}\ t\in (t_{j+1},t_j)\ \textup{and}\ j\in\llbracket m,\ell  \rrbracket.
\end{align}
The image of $\mathcal{Q}[(x,\ell )\to (y,m)]$ is the closed convex subset $ \{  x\leq t_\ell\leq t_{\ell-1}\leq \dots\leq t_{m+1}\leq y  \}.$ 
We often abuse the notation and view $\mathcal{Q}[(x,\ell )\to (y,m)]$ as a subset of $\mathbb{R}^{\ell-m}$.

\vspace*{0.3cm}

For $f\in C^{n}([a,b])$ and $\pi\in \mathcal{Q}[(x,\ell )\to (y,m)]$ with $(x,\ell), (y,m)\in [a,b]\times\llbracket 1,n \rrbracket$, define  
\begin{align}\label{def:onecurve}
f(\pi)\coloneqq \sum_{j=m}^{\ell} f_j(t_j)-f_{j}(t_{j+1}),
\end{align}
where $t_j$ are given by \eqref{equ:identify} and \eqref{equ:identifybd}. Let $d\pi$ be the Lebesgue measure on $\mathcal{Q}[(x,\ell)\to (y,m)]$. For $\beta>0$, the $\beta$-free energy from $(x,\ell)$ to $(y,m)$ is defined by

\begin{align}\label{def:freeenergy}
f[(x,\ell)\xrightarrow{\beta} (y,m)]\coloneqq \beta^{-1} \log\int_{\mathcal{Q}[(x,\ell)\to (y,m)]} \exp\left(\beta {f(\pi)}\right)\, d\pi.
\end{align}
We also allow $\beta=\infty$ and set
\begin{align*}
 f[(x,\ell)\xrightarrow{\infty} (y,m)]=\max_{\pi\in \mathcal{Q}[(x,\ell)\to (y,m)]} f(\pi).
\end{align*}
For $\beta=1$, we denote $f[(x,\ell)\xrightarrow{1} (y,m)]$ by $f[(x,\ell)\to (y,m)]$. The following lemma follows directly from the definition above.
\begin{lemma}\label{lem:integral}
For any $k\in \llbracket  m,\ell-1  \rrbracket$, we have
\begin{align*}
\exp\left( f[(x,\ell)\to (y,m)]\right)=  \int_x^y  \exp\left(f[(x,\ell)\to (z,k+1)]+f[(z,k)\to (y,m)]\right)\, dz.
\end{align*}
\end{lemma}
{\color{black}The next three lemmas, Lemmas~\ref{lem:f-convex}-\ref{lem:changeofvairable}, are elementary and have appeared widely in the literature. For completeness, we provide their proofs in the Appendix.}
\begin{lemma}\label{lem:f-convex}
Fix $n\geq \ell\geq m\geq 1$, $a\leq x_1\leq x_2<b$ and $f\in C^n([a,b])$. Then the function $f[(x_2,\ell)\to(y,m)]-f[(x_1,\ell)\to(y,m)]$ is monotone non-decreasing for $y\in(x_2,b]$.   
\end{lemma}

\begin{lemma}\label{lem:2022last}
Fix $n\geq \ell\geq m\geq 1$, $a\leq x < y_1\leq y_2\leq b$ and $f\in C^n([a,b])$. Then
\begin{align*}
f[(x,\ell)\to(y_1,m)]\leq  f[(x,\ell)\to(y_2,m)]-f_m(y_1)+f_m(y_2).
\end{align*}
\end{lemma}

\begin{lemma}\label{lem:changeofvairable}
Fix constants $a_1,a_2>0$, $a_3,a_4\in\mathbb{R}$ and $\{a_{5,i}\}_{i\in\mathbb{N}}$. For $g$ defined by $$g_i(x)=a_1f_i(a_2x+a_3)+a_4x+a_{5,i},$$ it holds that 
\begin{align*}
g[(x,\ell)\xrightarrow{\beta} (y,k)]=a_1\cdot f [(a_2x+a_3,\ell)\xrightarrow{a_1\beta} (a_2y+a_3,k)]\\
\hfill+a_4(y-x)-\beta^{-1}(\ell-k)\log a_2. 
\end{align*}
\end{lemma} 
\vspace*{0.3cm}

Next, we consider multiple paths that do not cross each other. Let $\pi_1$ and $\pi_2$ be two paths which belong to $\mathcal{Q}[(x_1,\ell_1)\to(y_1,m_1)]$ and $\mathcal{Q}[(x_2,\ell_2)\to(y_2,m_2)]$ respectively. We write $\pi_1 \prec\pi_2$ if $\pi_1(t)<\pi_2(t)$ for all $t\in (x_1,y_1)\cap (x_2,y_2)$. In this case, we say $\pi_1$ and $\pi_2$ are \textbf{non-intersecting}. The next lemma shows that non-intersecting paths form a closed convex set. 

\begin{lemma}\label{lem:convex}
For $i\in\{1,2\}$, let $(x_i,\ell_i) $ and $(y_i,m_i)$ be pairs with $x_i< y_i$ and $\ell_i\geq m_i$. Further assume $x_1\leq x_2$ and $y_1\leq y_2$. Then the collection of $(\pi_1,\pi_2)$ in $\mathcal{Q}[(x_1,\ell_1)\to (y_1,m_1)]\times \mathcal{Q}[(x_2,\ell_2)\to (y_2,m_2)]$ with $\pi_1\prec \pi_2$ is a closed convex subset in $\mathbb{R}^{\ell_1-m_1 }\times \mathbb{R}^{\ell_2-m_2 }$. 
\end{lemma}

A pair of sequences in $\mathbb{R}\times \mathbb{N}$ which can be connected by non-intersecting paths is called an endpoint pair. Its definition is given below.
\begin{definition}
Fix $k\in\mathbb{N}$. Let $U=\{(x_i,\ell_i) \}_{i\in \llbracket 1,k \rrbracket}\ \textup{and}\ V=\{(y_i,m_i) \}_{i\in \llbracket 1,k \rrbracket}$ be two sequences in $\mathbb{R}\times \mathbb{N}$ with $x_i< y_i$ and $\ell_i\geq m_i$ for all $i$. We denote by $\mathcal{Q}[U\to V]$ the collection of paths $\pi=(\pi_1,\dots,\pi_k)$ in  $\prod_{i=1}^k \mathcal{Q}[(x_i,\ell_i)\to (y_i,m_i)]$ that satisfy $\pi_1\prec\pi_2\prec\dots\prec\pi_k$. We call $(U,V)$ an \textbf{endpoint pair} if $\mathcal{Q}[U\to V]$ is non-empty and  $x_i\leq x_{i+1}$, $y_j\leq y_{j+1}$ for $i\in\llbracket 1,k-1\rrbracket$. We may call $(U,V)$ a $k$-endpoint pair to emphasize that there are $k$ pairs of endpoints.
\end{definition}

Let $(U,V)$ be a $k$-endpoint pair and $f\in C^{n}([a,b])$ with $U,V \subset [a,b]\times\llbracket 1,n \rrbracket$. For $\pi=(\pi_1,\dots,\pi_k)\in \mathcal{Q}[U\to V]$, we define
\begin{align*}
f(\pi)\coloneqq \sum_{i=1}^k f(\pi_i),
\end{align*}
where $f(\pi_i)$ are given in \eqref{def:onecurve}. In view of Lemma~\ref{lem:convex}, $\mathcal{Q}[U\to V]$ can be identified as a closed convex set in a Euclidean space. Let $p\in\mathbb{N}_0$ be the smallest integer such that $\mathcal{Q}[U\to V]$ is contained in a $p$-dimensional subspace and let $d\pi$ be the $p$-dimensional Hausdorff measure on $\mathcal{Q}[U\to V]$. We define 
\begin{align}\label{def:freeenergy-m}
f[U\to V]\coloneqq \log\int_{\mathcal{Q} [U\to V]} \exp\left( {f(\pi)}\right)\, d\pi.
\end{align}

The following reversing map will be used in Section~\ref{sec:gRSK}. For $f\in \overline{C}^n(T)$ and $z\in (0,T)$, we define $R_zf\in {C}^n([0,z])$ by
\begin{align}\label{def:reverse}
(R_zf)_i(t):=-f_{n+1-i}(z-t)+f_{n+1-i}(z). 
\end{align}
Let $U=\{(x_i,\ell_i) \}_{i\in \llbracket 1,k \rrbracket}\ \textup{and}\ V=\{(y_i,m_i) \}_{i\in \llbracket 1,k \rrbracket}$ be an endpoint pair with $U,V\subset [0,z]\times \llbracket 1,n \rrbracket$. Let
\begin{align*}
\widetilde{V}\coloneqq &\{(z-x_{k+1-i},n+1-\ell_{k+1-i}) \}_{i\in \llbracket 1,k \rrbracket},\\ \widetilde{U}\coloneqq &\{(z-y_{k+1-i},n+1-m_{k+1-i}) \}_{i\in \llbracket 1,k \rrbracket}.
\end{align*}
\begin{lemma}\label{lem:z-reverse}
Under the setting above, it holds that
\begin{align*}
f[U\to V]=(R_zf)[\widetilde{U}\to\widetilde{V}].
\end{align*}
\end{lemma}
It is convenient to introduce certain special sequences in $\mathbb{R}\times\mathbb{N}$. We use $(x,\ell)^k$ to denote the sequence $\{\underbrace{(x,\ell),(x,\ell),\dots, (x,\ell)}_{k\ \textup{terms}}\}.$ For $1\leq k\leq n$, we set
\begin{equation}\label{def:UVk}
\begin{split}
&V_k(x)\coloneqq\{ (x,1),(x,2),\dots, (x,k) \},\\
 &V'_k(x)\coloneqq\{ (x,2),(x,3),\dots, (x,k+1) \},\\ 
&U_{n,k}(x)\coloneqq\{ (x,n-k+1),(x,n-k+2),\dots, (x,n) \}.
\end{split}
\end{equation}
Moreover, we denote by $\mathcal{V}_{n,k}(x)$ the collection
\begin{equation}\label{def:UVk2}
 \mathcal{V}_{n,k}(x)= \{ (x,\ell_1),(x,\ell_2),\dots, (x,\ell_k)\, |\,  1\leq \ell_1<\ell_2<\dots <\ell_k\leq n\}.
\end{equation}

For paths in $\mathcal{Q}[(x,n)^k\to V]$, because of the non-intersecting requirement, the starting points need to pile up. Therefore, they belong to  $\mathcal{Q}[U_{n,k}(x) \to V]$. This is the content of the next lemma.

\begin{lemma}\label{lem:UVrepeatednew}
Fix $n\geq 2$, $2\leq k\leq n$, $T>0$ and $f\in {C}^n(T)$. Let 
$$U=\{ (x_1,n)^{i_1},(x_2,n)^{i_2},\dots, (x_\ell,n)^{i_\ell} \},\ V=\{ (y_1,n)^{j_1},(y_2,n)^{j_2},\dots, (y_m,n)^{j_m} \}$$ with $0<x_1<\dots<x_\ell<T$, $0<y_1<\dots<y_m<T$  and $\sum_{p=1}^\ell i_p=\sum_{q=1}^m j_q=k$. Suppose $(U,V)$ is an endpoint pair. Then
\begin{equation}\label{equ:UVrepeatednew1}
f[U\to V ]=f[\{U_{n,i_1}(x_1),U_{n,i_2}(x_2),\dots, U_{n,i_\ell}(x_\ell)\}\to \{V_{j_1}(y_1),V_{j_2}(y_2),\dots, V_{j_m}(y_m)\}].
\end{equation}
\end{lemma}


Lastly, we introduce down/right paths which are analogous to the up/right paths.
\begin{definition} 
For real numbers $x<  y$ and positive integers $m \leq \ell $, we use the notation $\mathcal{Q}[(x,m )\searrow (y,\ell )]$ to denote the collection of non-decreasing c\`adl\`ag functions $\rho:[x,y]\to\mathbb{N}$ with $\rho(x)\geq m$ and $\rho(y)=\ell$.
\end{definition}

There is an injective map from $\mathcal{Q}[(x,m )\searrow (y,\ell)]$ to $\mathbb{R}^{\ell-m}$ given by
\begin{align}\label{equ:identify-up}
t_j= \inf\{ t\in [x,y]\, |\, \rho(t)\geq j+1 \}\ \textup{for}\ j\in \llbracket m,\ell-1 \rrbracket. 
\end{align}
The image of $\mathcal{Q} [(x, m )\searrow (y,\ell)]$ is a closed convex subset and we often view $\mathcal{Q} [(x, m )\searrow (y,\ell)]$ as the subset of $\mathbb{R}^{\ell-m}$.
 
For $f\in C^{n}([a,b])$ and $ \rho\in \mathcal{Q} [(x, m )\searrow (y,\ell)]$ with $(x,m), (y,\ell)\in [a,b]\times\llbracket 1,n \rrbracket$, we define  
\begin{align*} 
f (\rho)\coloneqq \sum_{j=m}^{\ell } f_j(t_j)-f_{j}(t_{j-1}),
\end{align*}
where $t_j,\ j\in \llbracket m, \ell-1  \rrbracket$ are given by \eqref{equ:identify-up} and $t_{m-1}=x,\ t_{\ell}=y$. Let $d\rho$ be the Lebesgue measure on $\mathcal{Q}[(x,\ell)\to (y,m)]$. We define
\begin{align}\label{def:freeenergy-up}
f[(x,\ell)\searrow (y,m)]\coloneqq -\log\int_{\mathcal{Q} [(x,\ell)\searrow (y,m)]} \exp\left(- {f (\rho)}\right)\, d\rho.
\end{align} 
We finish this section with the lemma below which shows $f[V'_{k }(x)\to V_{k }(y)]$ and $f [(x,1)\searrow (y,k+1)]$ supplement each other. Here $V'_{k}(x)$ and $ V_{k}(y)$ are given in \eqref{def:UVk}.
\begin{lemma}\label{lem:searrow}
Fix $n\geq 2$, $1\leq k\leq n-1$, $0<x<y<T$ and $f\in C^n(T)$. Then it holds that
\begin{align*}
f[V'_{k }(x)\to V_{k }(y)]+f [(x,1)\searrow (y,k+1)]=f[V_{k+1}(x)\to V_{k+1}(y)].
\end{align*}
\end{lemma}

\section{Geometric RSK correspondence}\label{sec:gRSK}
In this section we define a geometric variant of the RSK correspondence introduced in \cite{OC12}. The main goal of this section is to derive the identity \eqref{equ:Wf+WRf} in Proposition~\ref{pro:Wf+WRf}. This identity describes the polymer energy for an environment under the geometric RSK and plays a crucial role in the convergence of the scaled KPZ sheets. {\color{black}This identity is derived using the invariance of polymer energy under the geometric RSK correspondence, established by Noumi and Yamada \cite{MR2074600} and more recently by Corwin \cite{Cor21}. In the context of the last passage percolation, an analogous invariance was proved in \cite{MR2176549} and, more recently, in \cite{DOV}.}

\vspace*{0.3cm}

Fix $n\geq 2 $, $1\leq i\leq n-1$ and $f\in C^n(T)$; see Definition~\ref{def:CT}. Define $\mathcal{T}_i f\in C^n(T)$ by
\begin{align*}
(\mathcal{T}_i f)_j(t)\coloneqq\left\{ \begin{array}{cc}
f_i(t)+\left( \log\int_0^t \exp(f_{i+1}(s)-f_{i}(s))ds \right), & j=i,\\
f_{i+1}(t)-\left( \log\int_0^t \exp(f_{i+1}(s)-f_{i}(s))ds \right), & j=i+1,\\
f_j(t), & j\neq i,i+1.
\end{array} \right.  
\end{align*}
We note that $\int_0^t\exp(f_{i+1}(s)-f_{i}(s))ds$ is understood as an improper integral because $f$ is not defined at $t=0$. The reason we adopt this formulation is the following. Even if we start with functions $f\in \overline{C}^n(T)$ which are continuous up to $t=0$, $\mathcal{T}_if$ is no longer continuous at $t=0$. More precisely, $(\mathcal{T}_if)_i(t)\sim \log t$ and $(\mathcal{T}_if)_{i+1}(t)\sim -\log t$ when $t$ goes to zero.  

For $1\leq r\leq n-1$, define
\begin{align*}
\mathcal{K}_r f\coloneqq \mathcal{T}_r\mathcal{T}_{r+1}\cdots \mathcal{T}_{n-1}f.
\end{align*}
\begin{definition}\label{def:gRSK}
Given $f\in\overline{C}^n(T)$ with $f(0)=0$, we define $\mathcal{W}f\in C^n(T)$ by
\begin{align}\label{def:W}
\mathcal{W}f\coloneqq 	\mathcal{K}_{n-1}\mathcal{K}_{n-2}\cdots \mathcal{K}_{1}f.
\end{align}
\end{definition}

The following version of Greene’s theorem was proved in \cite[page 445]{OC12}.
\begin{proposition}[\cite{OC12}] \label{pro:Greene}
Fix $n\geq 2$ and $f\in\overline{C}^{n}(T)$ with $f(0)=0$. Recall that $U_{n,k}(0)$ and $ V_k(t)$ are given in \eqref{def:UVk}. Then it holds for all $t\in (0,T)$ and $1\leq k\leq n $ that
\begin{align*}
\sum_{i=1}^k (\mathcal{W}f)_i(t)=f[U_{n,k}(0) \to V_k(t)].
\end{align*}
\end{proposition}
The following invariance of the free energy was proved in  \cite[Theorem 3.4]{Cor21}. It plays an important role in the proof of Proposition~\ref{pro:Wf+WRf}.
\begin{proposition}[\cite{Cor21}]\label{pro:gRSKinv}
Fix $n\geq 2$, $f\in \overline{C}^{n}(T)$ with $f(0)=0$ and an endpoint pair  $U=\{(x_i,n)\}_{i\in \llbracket 1,k \rrbracket}$ and $V=\{(y_i,1)\}_{i\in \llbracket 1,k \rrbracket}$ with $U,V\subset (0,T)\times \llbracket 1,n \rrbracket$. Further assume that $x_1<x_2<\dots <x_n$ and $y_1<y_2<\dots <y_n$.  Then it holds that
\begin{align*}
f[U\to V]=\left(\mathcal{W}f\right)[U\to V].
\end{align*}
\end{proposition}
The condition $x_1<x_2<\dots <x_n$ and $y_1<y_2<\dots <y_n$ can be removed through approximation and we obtain the corollary below. The proof can be found in Appendix.
\begin{corollary}\label{cor:gRSKinvnew}
Proposition~\ref{pro:gRSKinv} holds true without the condition $x_1<x_2<\dots <x_n$ and $y_1<y_2<\dots <y_n$.
\end{corollary}
The next proposition relates $(\mathcal{W}f)[(x,n)\to (z,k)]$ and  $(\mathcal{W}f)_k(z).$ 
\begin{proposition}\label{pro:Wf+WRf}
Fix $n\geq 2$, $1\leq k\leq  n-1$, $0<T$ and $f\in \overline{C}^n(T)$ with $f(0)=0$. For any $0<x<z<T$, it holds that
\begin{align}\label{equ:Wf+WRf}
(\mathcal{W}f)[(x,n)\to (z,k+1)]+(\mathcal{W}R_zf)[(z-x,1)\searrow (z,k+1)]=(\mathcal{W}f)_{k+1}(z).
\end{align} 

\end{proposition}
The rest of this section is devoted to proving Proposition~\ref{pro:Wf+WRf}. We start with a direct consequence of Lemma~\ref{lem:UVrepeatednew} and Corollary~\ref{cor:gRSKinvnew}. Recall that $U_{n,k}(x)$, $V_k(x)$ and $\mathcal{V}_{n,k}(x)$ are given in \eqref{def:UVk} and \eqref{def:UVk2} respectively.
\begin{corollary}\label{cor:gRSKinv}
Fix $n\geq 2$, $1\leq k\leq n-1$, $0<x<y<z<T$ and $f\in \overline{C}^n(T)$ with $f(0)=0$. The following identities hold.
\begin{align*}
f[U_{n,k+1}(x)\to V_{k+1}(y)] &= (\mathcal{W}f)[U_{n,k+1}(x)\to V_{k+1}(y)],\\
 f[\{U_{n,k }(x), (y,n)\}\to V_{k+1}(z)] &=(\mathcal{W}f)[\{U_{n,k }(x), (y,n)\}\to V_{k+1}(z)],\\
  f[U_{n,k+1}(x) \to \{(y,1),V_{k }(z) \} ] &=(\mathcal{W}f)[U_{n,k+1}(x) \to \{(y,1),V_{k }(z) \} ].
\end{align*} 
\end{corollary}

\begin{lemma}\label{lem:fill}
Fix $n\geq 2$, $1\leq k\leq n-1$, $0<x<y<T$ and $f\in \overline{C}^n(T)$ with $f(0)=0$. It holds that
\begin{align*}
(\mathcal{W}f)[(x,n)\to (y,k+1)]+ f[U_{n,k }(0) \to V_{k }(y)]=  f[\{U_{n,k }(0), (x,n)\}\to V_{k+1}(y)].
\end{align*}
\begin{proof}
Let $g=\mathcal{W}f$. Because of the natural measure-preserving injection from $\mathcal{Q}[\{U_{n,k }(\varepsilon ), (x,n)\}\to V_{k+1}(y)]$ to $\mathcal{Q} [U_{n,k }(\varepsilon )\to V_{k }(y)]\times \mathcal{Q}[(x,n)\to (y,k+1)]$, we have
\begin{align*}
g [\{U_{n,k }(\varepsilon ), (x,n)\}\to V_{k+1}(y)]\leq g[U_{n,k }(\varepsilon )\to V_{k }(y)]+g[(x,n)\to (y,k+1)].
\end{align*}
Take $\varepsilon$ go to zero and apply Corollary~\ref{cor:gRSKinv}, we get
\begin{align*}
f[\{U_{n,k }(0), (x,n)\}\to V_{k+1}(y)] \leq   f[U_{n,k }(0) \to V_{k }(y)]+g[(x,n)\to (y,k+1)] .
\end{align*}
Because of the natural measure-preserving injection from $$\mathcal{Q} [U_{n,k }(\varepsilon )\to V_{k }(x)]\times \mathcal{Q}[V_{k }(x )\to V_{k }(y)]\times \mathcal{Q}[(x,n)\to (y,k+1)]$$  to  $\mathcal{Q}[\{U_{n,k }(\varepsilon ), (x,n)\}\to V_{k+1}(y)]$, we have
\begin{align*}
g [\{U_{n,k }(\varepsilon ), (x,n)\}\to V_{k+1}(y)]\geq &	g[U_{n,k }(\varepsilon )\to V_{k }(x)]\\
 &+g[V_{k }(x )\to V_{k }(y)]+g[(x,n)\to (y,k+1)].
\end{align*}
Take $\varepsilon$ go to zero and apply Corollary~\ref{cor:gRSKinv} and Proposition~\ref{pro:Greene}, we get
\begin{align*}
f[\{U_{n,k }(0), (x,n)\}\to V_{k+1}(y)] \geq   f[U_{n,k }(0) \to V_{k }(y)]+g[(x,n)\to (y,k+1)] .
\end{align*}
\end{proof}
\end{lemma}
\begin{lemma}\label{lem:avoidnew}
Fix $n\geq 2$, $1\leq k\leq n$, $0<x<T$ and $f\in \overline{C}^n(T)$ with $f(0)=0$. Then for any $V\in\mathcal{V}_{n,k}(x)$ with $V\neq V_k(x)$, it holds that  
\begin{align*}
\lim_{\varepsilon\to 0} \exp\bigg( (\mathcal{W}f)[U_{n,k}(\varepsilon) \to  V]\bigg)=0.  
\end{align*}
\end{lemma}
\begin{proof}
Let $g=\mathcal{W}f$ and take $y\in(x,T)$. {\color{black}
We separate $\mathcal{Q}[U_{n,k}(\varepsilon) \to  V_k(y)]$ into disjoint subsets according to $$\{(x,\pi_1(x)),(x, \pi_2(x)),\dots,(x,\pi_k(x))\}\in \mathcal{V}_{n,k}(x).$$ Then we have  } 
\begin{align*}
\exp\bigg( g[U_{n,k}(\varepsilon) \to  V_k(y)] \bigg)=\sum_{V\in\mathcal{V}_{n,k}(x)} \exp\bigg( g[U_{n,k}(\varepsilon) \to  V]+g[ V  \to  V_k(y)] \bigg). 
\end{align*}
Therefore,
\begin{align*}
&\sum_{V\in\mathcal{V}_{n,k}(x), V\neq V_k(x)} \exp\bigg( g[U_{n,k}(\varepsilon) \to  V]+g[ V  \to  V_k(y)] \bigg)\\
= &\exp\bigg( g[U_{n,k}(\varepsilon) \to  V_k(y)] \bigg)-\exp\bigg( g[U_{n,k}(\varepsilon) \to  V_k(x)]+g[ V_k(x)  \to  V_k(y)] \bigg).
\end{align*}
From Corollary~\ref{cor:gRSKinv}, the limit of the right hand side equals
\begin{align*}
\exp\bigg( f[U_{n,k}(0) \to  V_k(y)] \bigg)-\exp\bigg( f[U_{n,k}(0) \to  V_k(x)]+g[ V_k(x)  \to  V_k(y)] \bigg).
\end{align*}
From Proposition~\ref{pro:Greene}, the above vanishes. Therefore for any $V\in \mathcal{V}_{n,k}(0)$ with $V\neq V_k(x)$, we have $$\lim_{\varepsilon\to 0} \exp\big(
g[U_{n,k}(\varepsilon)\to V]\big)=0.$$
\end{proof}

\begin{lemma}\label{lem:reverse}
Fix $n\geq 2$, $1\leq k\leq n-1$, $0<x<y<T$ and $f\in\overline{C}^n(T)$ with $f(0)=0$. Then
\begin{align*}
f[U_{n,k+1}(0) \to\{ (x,1), V_{k}(y)  \}]=f[U_{n,k+1}(0) \to V_{k+1}(y)  ]-(\mathcal{W}f)[(x,1)\searrow (y,k+1)].
\end{align*}
\end{lemma}
\begin{proof}
Let $g=\mathcal{W}f$. From Corollary~\ref{cor:gRSKinv}, $f[U_{n,k+1}(0) \to\{ (x,1), V_{k }(y)  \}]$ equals
$$\lim_{\varepsilon\to 0} g[U_{n,k+1}(\varepsilon ) \to\{ (x,1), V_{k }(y)  \}].$$
From Lemma~\ref{lem:avoidnew}, the above becomes
$$\lim_{\varepsilon\to 0} g[U_{n,k+1}(\varepsilon ) \to V_{k+1}(x) ]+g[V'_{k }(x)\to V_{k}(y)].$$
From Corollary~\ref{cor:gRSKinv} and Lemma~\ref{lem:searrow}, the above equals
$$f[U_{n,k+1}(0 ) \to V_{k+1}(x) ]+g[V_{k+1}(x)\to V_{k+1}(y)]-g [(x,1)\searrow (y,k+1)].$$
In view of Proposition~\ref{pro:Greene}, the above becomes $f[U_{n,k+1}(0 ) \to V_{k+1}(y) ]-g [(x,1)\searrow (y,k+1)]$.
\end{proof}

\begin{proof}[Proof of Proposition~\ref{pro:Wf+WRf}]
From Lemma~\ref{lem:fill}, $(\mathcal{W}f)[(x,n)\to (z,k+1)]$ equals
\begin{align*}
f[\{U_{n,k }(0),(x,n)\}\to V_{k+1}(z)]-f[U_{n,k }(0)\to V_{k }(z)].
\end{align*}
From Lemma~\ref{lem:z-reverse}, the above equals
\begin{align*}
(R_zf)[U_{n,k+1}(0) \to \{(z-x,1),V_{k }(z)\} ]-f[U_{n,k }(0)\to V_{k }(z)].
\end{align*}
From Lemma~\ref{lem:reverse}, the above equals
\begin{align*}
(R_zf)[U_{n,k+1}(0)\to V_{k+1}(z)]-(\mathcal{W} R_zf)[(z-x,1)\searrow (z,k+1)]-f[U_{n,k }(0)\to V_{k }(z)].
\end{align*}
Applying Lemma~\ref{lem:z-reverse} and Proposition~\ref{pro:Greene}, it becomes
\begin{align*}
(\mathcal{W}f)_{k+1}(z)-(\mathcal{W} R_zf)[(z-x,1)\searrow (z,k+1)].
\end{align*}
\end{proof}

\section{O’Connell-Yor polymer model}\label{sec:OY}
In this section we recall the O’Connell-Yor polymer model and show that the cumulative distribution functions (c.d.f.) of the polymer measure bound the difference between Busemann functions, Lemmas~\ref{lem:FF} and \ref{lem:GG}.

\vspace{0.3cm}
 
Let $B_1,B_2,\dots$ be i.i.d. standard two-sided Brownian motions. For $n\in\mathbb{N}$, let $B^n$ be the first $n$ Brownian motions $(B_1,B_2,\dots, B_n)$ restricted on $[0,\infty)$ and define $Y^n=(Y^n_1,\dots,Y^n_n)\coloneqq \mathcal{W}B^n$. 
For $n,i\in\mathbb{N}$ and $T>0$, let
\begin{align*}
 C_1(T,n)\coloneqq & n^{1/2}T^{-1/2}+2^{-1},\ C_{3,i}(T)\coloneqq -(i-1)\log T+\log(i-1)!, \\
C_2(T,n)\coloneqq & n+2^{-1}n^{1/2}T^{1/2}- (n-1)\log( n^{1/2}T^{-1/2}) .
\end{align*}
Recall that $\mathcal{X}^{T,n}=\{\mathcal{X}^{T,n}_1,\mathcal{X}^{T,n}_2,\dots, \mathcal{X}^{T,n}_n\}$ are given  by 
\begin{equation}\label{def:kpzline_n} 
\mathcal{X}^{T,n}_i(x)\coloneqq Y^n_i(n^{1/2}T^{1/2}+x)-C_1(T,n)x-C_2(T,n)-C_{3,i}(T) ,\ i\in\llbracket 1,n \rrbracket.
\end{equation}
From Proposition~\ref{pro:Greene}, 
\begin{align*}
\mathcal{X}^{T,n}_1(x)= B^n[(0,n)\to (n^{1/2}T^{1/2}+x,1)] -C_1(T,n)x-C_2(T,n) .
\end{align*}
For $T>0$, $n\geq 1$, $x\in\mathbb{R} $ and $y>-n^{1/2}T^{1/2}+x$, recall that
\begin{align}\label{def:kpzsheet_n}
\mathbcalboondox{h}^{T,n}(x,y)\coloneqq B^n[(x,n)\to (n^{1/2}T^{1/2}+y,1)]-C_1(T,n)(y-x)-C_2(T,n).
\end{align}
Note that $\mathbcalboondox{h}^{T,n}(x,y)$ has the same distribution as $\mathcal{X}^{T,n}_1(y-x)$. In view of Proposition~\ref{pro:gRSKinv}, for $x>0$, we may rewrite $\mathbcalboondox{h}^{T,n}(x,y)$ as polymer free energies on $Y^n$.
\begin{align}\label{equ:SX}
\mathbcalboondox{h}^{T,n}(x,y)=Y^n[(x,n)\to (n^{1/2}T^{1/2}+y,1)]-C_1(T,n)(y-x)-C_2(T,n).
\end{align}


For $1\leq k\leq n-1$, $T>0$, $x>0$ and $y> z>-n^{1/2}T^{1/2}+x$, we define
\begin{align*}
F^{T,n}_k(x,z)\coloneqq &Y^n[(x,n)\to (n^{1/2}T^{1/2}+z,k+1)]-Y^n_{k+1}(n^{1/2}T^{1/2}+z)+C_1(T,n)x,
\end{align*}
and
\begin{align*}
G^{T,n}_k(z,y)\coloneqq  Y^n[ (n^{1/2}T^{1/2}+z,k )\to (n^{1/2}T^{1/2}+y,1)]&+Y^n_{k+1}(n^{1/2}T^{1/2}+z)\\
&-C_1(T,n)y-C_2(T,n).
\end{align*}
Essentially, $F^{T,n}_k(x,z)$ and $G^{T,n}_k(z,y)$ are respectively polymer energies from $(x,n)$ to $( n^{1/2}T^{1/2}+z,k+1)$ and from $(n^{1/2}T^{1/2}+z,k )$ to $(n^{1/2}T^{1/2}+y,1)$ in the environment $Y^{n}$. From Lemma~\ref{lem:integral} and \eqref{equ:SX}, we have 
\begin{equation}\label{equ:HFG}
\exp\left(\mathbcalboondox{h}^{T,n}(x,y)\right)=\int_{-n^{1/2}T^{1/2}+x}^y \exp\left(F^{T,n}_k(x,z)+G^{T,n}_k(z,y)\right)\, dz.
\end{equation}
We note that from Lemma~\ref{lem:changeofvairable} and \eqref{def:kpzline_n}, $F^{T,n}_k(x,z)$, $G^{T,n}_k(z,y)$ and $\mathbcalboondox{h}^{T,n} (x,y)$ can be expressed in terms of $\mathcal{X}^{T,n}$ as
\begin{equation*}
\begin{split}
F^{T,n}_k(x,z)=\mathcal{X}^{T,n}[(-n^{1/2}T^{1/2}+x,n )&\to (z,k+1)]-\mathcal{X}^{T,n}_{k+1}(z)\\
& + (n-1)\log (n^{1/2}T^{-1/2})-C_{3,k+1}, 
\end{split}
\end{equation*}
\begin{align}\label{equ:GH}
G^{T,n}_k(z,y)=\mathcal{X}^{T,n}[(z,k )\to (y,1)]+\mathcal{X}^{T,n}_{k+1}(z)+C_{3,k+1}, 
\end{align}
and
\begin{align*} 
\mathbcalboondox{h}^{T,n} (x,y)=\mathcal{X}^{T,n}[(-n^{1/2}T^{1/2}+x ,n)\to (y,1)]+ (n-1)\log (n^{1/2}T^{-1/2}).  
\end{align*}
The next lemma concerns the distributional limit of $F^{T,n}_k(x,z)$ when $n$ goes to infinity.
\begin{lemma}\label{lem:Flimit}
Fix $T>0$, $k\geq 1$, $x>0$ and $z\in\mathbb{R}$. Then when $n$ goes to infinity, $F^{T,n}_k(x,z)$ converges in distribution to $ \mathcal{X}^T[(0,k+1)\to(x,1)]+T^{-1}zx$.
\begin{proof}
From Proposition~\ref{pro:Wf+WRf}, and $Y^n=\mathcal{W}B^n$, $F^{T,n}_k(x,z)$ equals
\begin{align*}
-(\mathcal{W}R_{n^{1/2}T^{1/2}+z}B^n)[(n^{1/2}T^{1/2}+z-x,1)\searrow (n^{1/2}T^{1/2}+z,k+1)]+C_1(T,n)x.
\end{align*} 
Because $\mathcal{W}R_{n^{1/2}T^{1/2}+z}B^n\overset{d}{=}Y^n$, the above has the same distribution as
\begin{align*}
-Y^n[(n^{1/2}T^{1/2}+z-x,1)\searrow (n^{1/2}T^{1/2}+z,k+1)]+C_1(T,n)x.
\end{align*}
From \eqref{def:kpzline_n} and Lemma~\ref{lem:changeofvairable}, the above equals $-\mathcal{X}^{T,n}[( z-x,1)\searrow ( z,k+1)] $. By Proposition~\ref{pro:kpzlineensemble} it converges in distribution to $-\mathcal{X}^T[( z-x,1)\searrow ( z,k+1)]$. Because $\mathcal{X}^T(y)$ has the same distribution as $\mathcal{X}^{T}(-y)$,
$$-\mathcal{X}^T[( z-x,1)\searrow ( z,k+1)]\overset{d}{=}\mathcal{X}^T[(-z,k+1)\to (-z+x,1)].$$
From the stationarity of $\mathcal{X}^{T}(y)+2^{-1}T^{-1}y^2$, $\mathcal{X}^{T}(-z+y)\overset{d }{=}\mathcal{X}^T(y)+T^{-1} zy-2^{-1}T^{-1}z^2$. Therefore,
\begin{align*}
\mathcal{X}^{T}[(-z,k+1)\to (-z+x,1)]\overset{d}{=} \mathcal{X}^{T}[(0,k+1)\to (x,1)]+T^{-1}zx.
\end{align*}
\end{proof}
\end{lemma}

In the rest of the section, we drive a relation between Busemann functions and the c.d.f. of polymer measures, Lemmas~\ref{lem:FF} and \ref{lem:GG}. Those are deterministic properties and do not rely on the laws of $\mathcal{X}^{T,n}$ or $\mathbcalboondox{h}^{T,n}$.

\vspace*{0.3cm}

We start with a simple consequence of Lemma~\ref{lem:f-convex}.

\begin{lemma}\label{lem:Fconvex}
Fix $T>0$, $n\geq 2$, $1\leq k\leq n-1$, $0<x_1\leq x_2$ and $-n^{1/2}T^{1/2}<y_1\leq y_2$. Then $F^{T,n}_k(x_2,z)-F^{T,n}_k(x_1,z)$ is monotone non-decreasing in $z\in (-n^{1/2}T^{1/2}+x_2,\infty)$ and $G^{T,n}_k(z,y_2)-G^{T,n}_k(z,y_1)$ is monotone non-decreasing in $z\in (-n^{1/2}T^{1/2} ,y_1)$.  
\end{lemma}
We define the random probability measure on $\mathbb{R}$ which corresponds to \eqref{equ:HFG}. It is the marginal of the polymer measure.
\begin{definition}\label{def:AB}
Fix $T>0$, $n\geq 2$, $1\leq k\leq n-1$,  $x>0$ and $y>-n^{1/2}T^{1/2}+x$. We denote by $d\mu^{T,n}_{k,x,y}(z)$ the random probability measure with the density $$\exp\left(-\mathbcalboondox{h}^{T,n}(x,y)+F^{T,n}_k(x,z)+G^{T,n}_k(z,y)\right)\mathbbm{1}(-n^{1/2}T^{1/2}+x<z<y).$$ We also set its c.d.f.
\begin{equation*}
A^{T,n}_k(x,y;z)\coloneqq \mu^{T,n}_{k,x,y}([z,\infty)),\ B^{T,n}_k(x,y;z)\coloneqq \mu^{T,n}_{k,x,y}((-\infty,z]).
\end{equation*}
\end{definition}

\begin{lemma}\label{lem:FF}
Fix $T>0$, $n\geq 2$, $1\leq k\leq n-1$, $x_2\geq x_1>0$ and $y>-n^{1/2}T^{1/2}+x_2$. Then for $z\in (-n^{1/2}T^{1/2}+x_2,y)$, we have
\begin{align}\label{equ:FF<A}
F^{T,n}_k(x_2,z)-F^{T,n}_k(x_1,z)\leq &\mathbcalboondox{h}^{T,n}(x_2,y)-\mathbcalboondox{h}^{T,n}(x_1,y)-\log A^{T,n}_k(x_1,y;z), 
\end{align}
and
\begin{align}\label{equ:FF>B}
F^{T,n}_k(x_2,z)-F^{T,n}_k(x_1,z)\geq &\mathbcalboondox{h}^{T,n}(x_2,y)-\mathbcalboondox{h}^{T,n}(x_1,y)+\log B^{T,n}_k(x_2,y;z).
\end{align}  
\end{lemma}
\begin{proof}
We start with \eqref{equ:FF<A}.
\begin{align*}
&\exp\left(\mathbcalboondox{h}^{T,n}(x_2,y)-\mathbcalboondox{h}^{T,n}(x_1,y)\right)\\
=&\int_{-n^{1/2}T^{1/2}+x_2}^y \exp\left(F^{T,n}_k(x_2,z')-F^{T,n}_k(x_1,z')\right) d\mu^{T,n}_{k,x_1,y}(z')\\
\geq &\int_{z}^y \exp\left(F^{T,n}_k(x_2,z')-F^{T,n}_k(x_1,z')\right) d\mu^{T,n}_{k,x_1,y}(z')\\
\geq & \exp\left(F^{T,n}_k(x_2,z)-F^{T,n}_k(x_1,z)\right) A^{T,n}_k(x_1,y;z).
\end{align*}
We used Lemma~\ref{lem:Fconvex} in the last inequality. Then \eqref{equ:FF<A} follows. For \eqref{equ:FF>B}, we derive similarly, 
\begin{align*}
&\exp\left(\mathbcalboondox{h}^{T,n}(x_1,y)-\mathbcalboondox{h}^{T,n}(x_2,y)\right)\\
 \geq & \int_{-n^{1/2}T^{1/2}+x_2}^y \exp\left(F^{T,n}_k(x_1,z')-F^{T,n}_k(x_2,z')\right) d\mu^{T,n}_{k,x_2,y}(z')\\
 \geq & \int_{-n^{1/2}T^{1/2}+x_2}^z \exp\left(F^{T,n}_k(x_1,z')-F^{T,n}_k(x_2,z')\right) d\mu^{T,n}_{k,x_2,y}(z')\\
\geq &\exp\left(F^{T,n}_k(x_1,z )-F^{T,n}_k(x_2,z )\right)B^{T,n}_k(x_2,y;z).
\end{align*}
We again used Lemma~\ref{lem:Fconvex} in the last inequality. Hence \eqref{equ:FF>B} follows.
\end{proof}
The lemma below is analogous to Lemma~\ref{lem:FF} and we omit the proof.
\begin{lemma}\label{lem:GG}
Fix $T>0$, $n\geq 2$, $1\leq k\leq n-1$, $x>0$ and $y_2\geq y_1>-n^{1/2}T^{1/2}+x$. Then for $z\in (-n^{1/2}T^{1/2}+x ,y_1)$, we have
\begin{align}\label{equ:GGA}
G^{T,n}_k(z,y_2 )-G^{T,n}_k(z,y_1)\leq &\mathbcalboondox{h}^{T,n}(x,y_2)-\mathbcalboondox{h}^{T,n}(x ,y_1)-\log A^{T,n}_k(x,y_1;z), 
\end{align}
and
\begin{align}\label{equ:GGB}
G^{T,n}_k(z,y_2 )-G^{T,n}_k(z,y_1)\geq &\mathbcalboondox{h}^{T,n}(x,y_2)-\mathbcalboondox{h}^{T,n}(x ,y_1)+\log B^{T,n}_k(x,y_2;z).
\end{align}  
\end{lemma}

\vspace*{0.3cm}

The following two inequalities will be applied in Section \ref{sec:proof1} and Section \ref{sec:proof2}. {\color{black}Under the same setting as Lemma~\ref{lem:GG},}  we have
\begin{align}\label{equ:HHB0}
\begin{split} 
\mathcal{X}^{T,n}[(z,k)\to (y_2,1)]-\mathcal{X}^{T,n}[(z,k)\to (y_1,1)]  -\mathbcalboondox{h}^{T,n}&(x,y_2)+\mathbcalboondox{h}^{T,n}(x ,y_1)\\
&\leq  -\log\left(1- B^{T,n}_k(x,y_1,z)\right), 
\end{split}
\end{align}
and
\begin{align}\label{equ:HHA0} 
\begin{split} 
\mathcal{X}^{T,n}[(z,k)\to (y_2,1)]-\mathcal{X}^{T,n}[(z,k)\to (y_1,1)] - \mathbcalboondox{h}^{T,n}(x&,y_2)+\mathbcalboondox{h}^{T,n}(x ,y_1)\\
&\geq  \log\left(1- A^{T,n}_k(x,y_2,z)\right).
\end{split}  
\end{align}
{\color{black} The bounds \eqref{equ:HHB0} and \eqref{equ:HHA0} follow directly from \eqref{equ:GH}, \eqref{equ:GGA}, \eqref{equ:GGB} and $A^{T,n}_k(x,y;z)+B^{T,n}_k(x,y;z)=1$.}
\section{Proof of Theorem \ref{thm:KPZtoAiry_sheet}}\label{sec:proof1}

We present the proof for Theorem~\ref{thm:KPZtoAiry_sheet} in this section. We begin with convergence results for scaled KPZ line ensembles. 

\vspace*{0.3cm}
Recall that in \eqref{def:KPZ_sheet_rescaled}, the scaled KPZ sheet $\mathfrak{h}^T$ is obtain by performing the $1:2:3$ scaling to the KPZ sheet $\mathbcalboondox{h}^T$ as
\begin{align*}
\mathfrak{h}^T(x,y)\coloneqq 2^{1/3}T^{-1/3}\mathbcalboondox{h}^T(2^{1/3}T^{2/3}x,2^{1/3}T^{2/3}y)+ {2^{1/3}T^{2/3}}/{24}.
\end{align*}
The scaled KPZ line ensemble $\mathfrak{X}^T=\{\mathfrak{X}^T_1,\mathfrak{X}^T_2,\dots\}$ is defined accordingly as
\begin{equation}\label{def:0310}
\mathfrak{X}^T_i(x)\coloneqq 2^{1/3}T^{-1/3}\mathcal{X}_i^T(2^{1/3}T^{2/3}x)+ {2^{1/3}T^{2/3}}/{24}.
\end{equation}
We further perform the same scaling to objects in the O'Connell-Yor model as
\begin{equation}\label{def:scaled_prelimit}
\begin{split}
\mathfrak{X}^{T,n}_i(x)\coloneqq &2^{1/3}T^{-1/3}\mathcal{X}_i^{T,n}(2^{1/3}T^{2/3}x)+ {2^{1/3}T^{2/3}}/{24},\\ 
\mathfrak{h}^{T,n}(x,y)\coloneqq &2^{1/3}T^{-1/3}\mathbcalboondox{h}^{T,n}(2^{1/3}T^{2/3}x,2^{1/3}{T}^{2/3}y)+ 2^{1/3}{T^{2/3}}/{24}.
\end{split}
\end{equation}
We note that as $n$ goes to infinity, $\mathfrak{X}^{T,n}$ converges to $\mathfrak{X}^{T}$ (Proposition~\ref{pro:kpzlineensemble}) and $\mathfrak{h}^{T,n}$ converges to $\mathfrak{h}^{T}$ (Proposition~\ref{pro:kpzsheet}).

The convergence of the scaled KPZ line ensemble to the parabolic Airy line ensemble is a consequence of a series of works \cite{QS,V,DM,Wu21tightness,aggarwal2023}.
\begin{proposition}[\cite{QS,V,DM,Wu21tightness,aggarwal2023}]\label{pro:KPZtoAiry_line}
The scaled KPZ line ensemble $\mathfrak{X}^T $ converges in distribution to the parabolic Airy line ensemble $\mathcal{A}$ when $T$ goes to infinity. Here $\mathfrak{X}^T$ and $\mathcal{A}$ are considered as $C(\mathbb{N}\times\mathbb{R},\mathbb{R})$-valued random variables. 
\end{proposition} 
 
Next, we show the tightness for the scaled KPZ sheets. We rely on the following result by Dauvergne and Vir\'{a}g \cite[Lemma 3.3]{DV21}.
\begin{proposition}[\cite{DV21}]\label{pro:DV}
Let $Q = I_1\times \dots\times I_d$ be a product of bounded real intervals of length $b_1,\dots , b_d$. 
Let $c, a > 0$. Let $\mathcal{G}$ be a random continuous function from $Q$ taking values in a real vector space $V$ with Euclidean norm $|\cdot|$. Assume that for every $i\in \{1,2,\dots, d\}$, there exist $\alpha_i\in (0,1)$,  $\beta_i, r_i>0$ such that
\begin{align*}
\mathbb{P}\left( \left| \mathcal{G}(t+ue_i)-\mathcal{G}(t) \right|\geq  m u^{\alpha_i} \right)\leq ce^{-a m^{\beta_i}},
\end{align*} 
for every coordinate vector $e_i$, every $m > 0$, and every $t, t + u e_i \in  Q$ with $0\leq u < r_i$. Set $\beta=\min_i\beta_i$, $\alpha=\max_i\alpha_i$ and $r=\max_i r_i^{\alpha_i}$. Then with probability one it holds that
\begin{align*}
\left| \mathcal{G}(t+s)-\mathcal{G}(t) \right|\leq C\left( \sum_{i=1}^d |s_i|^{\alpha_i}\log^{1/\beta_i}\left( 2r^{1/\alpha_i}/|s_i| \right) \right),
\end{align*}
for every $t,t+s\in Q$  with $ |s_i|\leq r_i $ for all $i$ (here $s=(s_1,\dots, s_d)$). Here $C$ is a random constant satisfying $$\mathbb{P}(C>m)\leq \left( \prod_{i=1}^d b_i/r_i \right) cc_0 e^{-c_1m^\beta}, $$
where $c_0$ and $c_1$ are constants that depend on $\alpha_1,\dots, \alpha_d$, $\beta_1,\cdots,\beta_d $, $d$, and $a$.
\end{proposition}
\begin{proposition}\label{pro:kpzsheet_tight}
When $T$ goes to infinity, the scaled KPZ sheet $\mathfrak{h}^T$ is tight in $C(\mathbb{R}^2,\mathbb{R})$.
\end{proposition}
\begin{proof}
It suffices to prove that for all $b>0$, $\mathfrak{h}^T$ restricted on $Q=[-b,b]^2$ is tight in $C(Q,\mathbb{R})$. The tightness of $\mathfrak{h}^T(0,0)=\mathfrak{H}^T(0,0;1,0)$ follows directly from its convergence \cite{ACQ}. It remains to control the modulus of continuity. For simplicity, we denote by $D$ a constant that depends only on $b$. The value of $D$ may increase from line to line. 

From \cite[Theorem 1.3]{CGH}, there exists $T_0>0$ such that for all $T\geq T_0$, $d\in (0,1]$, $K\geq 0$ and $x,x+d\in [-2b,2b]$, we have
\begin{align}\label{equ:OYsheet_tight0}
\mathbb{P}\left( |\mathfrak{h}^{T}(0,x+d)- \mathfrak{h}^{T} (0,x)|  >K d^{1/2} \right)\leq D e^{-D^{-1}K^{3/2}}. 
\end{align}
From (3) in Proposition~\ref{thm:forcoupling}, $\mathfrak{h}^{T}(\cdot,y)\overset{d}{=}\mathfrak{h}^{T}(0,y-\cdot)$ and $\mathfrak{h}^{T}(x,\cdot )\overset{d}{=}\mathfrak{h}^{T}(0,\cdot-x )$. Therefore, for any $(x,y)\in Q$, 
\begin{align*} 
\mathbb{P}\left( |\mathfrak{h}^{T}(x+d,y)-\mathfrak{h}^{T}(x,y)|  >K d^{1/2} \right)\leq De^{-D^{-1}K^{3/2}}, 
\end{align*}
provided $x+d\in [-b,b]$. Similarly, if $y+d\in [-b,b]$, then
\begin{align*} 
\mathbb{P}\left( |\mathfrak{h}^{T}(x ,y+d)-\mathfrak{h}^{T}(x,y)|  >K d^{1/2} \right)\leq D e^{-D^{-1}K^{3/2}}, 
\end{align*}
From Proposition~\ref{pro:DV}, there exists a random constant $C^{T}$ such that almost surely for all $(x,y)$,  $(x',y')$ in $ Q$ with $|x-x'|, |y-y'|\leq 1$, we have
\begin{align*}
|\mathfrak{h}^{T}(x ,y )-\mathfrak{h}^{T}(x',y')|\leq C^{T }\left( |x-x'|^{1/2}\log^{2/3}(2/|x-x'|)+|y-y'|^{1/2}\log^{2/3}(2/|y-y'|)\right).
\end{align*}
Moreover, $\mathbb{P}(C^{T}>K)< De^{-D^{-1}K^{3/2}}$. By the Kolmogorov-Chentsov criterion (see Theorem 23.7 in \cite{Kal}) this implies the tightness of $\mathfrak{h}^{T}$ restricted on $Q$.
\end{proof}
The next proposition shows that any subsequential limit of $\mathfrak{h}^T$ and the parabolic Airy line ensemble can be coupled together with desired properties.
\begin{proposition}\label{pro:Airy_S}
Let $\mathfrak{h}$ be a distributional limit of $\mathfrak{h}^T$ along some sequence. Then there exists a coupling of $\mathfrak{h}$ and the parabolic Airy line ensemble $\mathcal{A}$ such that the following holds.
\begin{enumerate}
\item $\mathfrak{h}(0,\cdot)=\mathcal{A}_1(\cdot)$.
\item Almost surely for all $x>0$ and $y_1,y_2$ in $\mathbb{R}$, we have
\begin{equation}\label{equ:Airy_S}
\begin{split}
\lim_{k\to\infty} \mathcal{A}[(-2^{-1/2}k^{1/2}x^{-1/2},k)\xrightarrow{\infty} (y_2,1)] - \mathcal{A}[(-2^{-1/2}&k^{1/2}x^{-1/2},k)\xrightarrow{\infty} (y_1,1)]\\
& = \mathfrak{h}(x,y_2)-\mathfrak{h}(x,y_1). 
\end{split}
\end{equation}
\end{enumerate}
   
\end{proposition}
\begin{proof}[Proof of Theorem~\ref{thm:KPZtoAiry_sheet}]
Let $\mathfrak{h} $ be a distributional limit of $\mathfrak{h}^T $ along some sequence. From Proposition~\ref{pro:Airy_S}, \eqref{equ:Airy_S} holds. Because of (3) in Proposition~\ref{thm:forcoupling}, $\mathfrak{h}(\cdot +t,\cdot +t)$ has the same distribution as $\mathfrak{h}(\cdot,\cdot)$. From Definition~\ref{def:Airy_sheet}, $\mathfrak{h} $ has the same law as the Airy sheet. As a result, $\mathfrak{h}^T $ converges to the Airy sheet in distribution.
\end{proof}

It remains to prove Proposition~\ref{pro:Airy_S}. Let
\begin{equation}\label{def:frakR}
\mathfrak{R} ^{T,n}_k(x,z)\coloneqq  2^{1/3}T^{-1/3}F^{T,n}_k(2^{1/3}T^{2/3}x,2^{1/3}T^{2/3}z)-2^{3/2}k^{1/2}x^{1/2}- 2zx .
\end{equation}
This random variable is defined so that first, $F^{T,n}_k$ is suitably scaled according to \eqref{def:scaled_prelimit}; second, the major term in $F^{T,n}_k$ is subtracted such that $\mathfrak{R} ^{T,n}_k(x,z)=o(k^{1/2})$. See Proposition~\ref{pro:cpl_Airy} and its proof for details. It is convenient to define the scaled c.d.f. as
\begin{equation}\label{def:frakAB}
\begin{split}
\mathfrak{A}^{T,n}_k(x,y;z)\coloneqq &A^{T,n}_k(2^{1/3}T^{2/3}x,2^{1/3}T^{2/3}y;2^{1/3}T^{2/3} z   ),\\
\mathfrak{B}^{T,n}_k(x,y;z)\coloneqq &B^{T,n}_k(2^{1/3}T^{2/3}x,2^{1/3}T^{2/3}y;2^{1/3}T^{2/3} z   ).
\end{split}
\end{equation}
{\color{black}The next lemma can be derived from Lemma~\ref{lem:FF}.}
\begin{lemma}\label{lem:boundzz}
Fix $T\geq 2$, $n\geq 2$, and $1\leq k\leq n-1$. Then for all $x, \bar{x}>0$, and $y>-2^{-1/3}n^{1/2}T^{-1/6}+\max\{x,\bar{x}\}$, the following statements hold. Let $\bar{z}=-2^{-1/2}k^{1/2}\bar{x}^{-1/2}$. If $\bar{x} \geq x$, then
\begin{align}\label{equ:Airy_Abound} 
\begin{split} 
  \log\mathfrak{A}^{T,n}_k(x,y;\bar{z})+ 2^{1/2}k^{1/2}\bar{x} ^{1/2}&\left( 1-\bar{x}^{-1/2}x^{1/2} \right)^2 \\
&  \leq  \mathfrak{h}^{T,n}(\bar{x} , y)-\mathfrak{h}^{T,n}( x, y) -\mathfrak{R}^{T,n}_k( \bar{x} , \bar{z})+ \mathfrak{R} ^{T,n}_k( x, \bar{z}).
\end{split}
\end{align}
If $\bar{x}\leq x $, then
\begin{align}\label{equ:Airy_Bbound} 
\begin{split} 
  \log \mathfrak{B}^{T,n}_k(x,y;\bar{z} ) + 2^{1/2} k^{1/2}\bar{x} ^{1/2}&\left( 1-\bar{x}^{-1/2}x^{1/2} \right)^2 \\
&\leq  \mathfrak{h}^{T,n}(\bar{x} , y)-\mathfrak{h}^{T,n}( x, y) -\mathfrak{R}^{T,n}_k( \bar{x} , \bar{z})+ \mathfrak{R} ^{T,n}_k( x, \bar{z}).
\end{split}
\end{align}
\begin{proof}
First, we consider the case $\bar{x} \geq x$.  From \eqref{equ:FF<A}, \eqref{def:scaled_prelimit} and \eqref{def:frakAB}, we have
\begin{align*}
& 2^{-1/3}T^{1/3}\big(\mathfrak{h}^{T,n}( \bar{x} , y)- \mathfrak{h}^{T,n}( x , y)\big)-\log \mathfrak{A}^{T,n}_k(x,y;\bar{z} ) \\
 \geq & F^{T,n}_k(2^{1/3}T^{2/3}\bar{x}  ,2^{1/3}T^{2/3}\bar{z} )-F^{T,n}_k(2^{1/3}T^{2/3} x  ,2^{1/3}T^{2/3}\bar{z} ). 
\end{align*}
From \eqref{def:frakR}, the right hand side equals
\begin{align*}
2^{-1/3}T^{1/3}\bigg(2^{1/2}  k^{1/2}\bar{x}^{1/2}\left( 1-\bar{x}^{-1/2}x^{1/2} \right)^2 
 + \mathfrak{R}^{T,n}_k( \bar{x} , \bar{z})-\mathfrak{R}^{T,n}_k( x  , \bar{z})\bigg).
\end{align*}
Together with $T\geq 2$, \eqref{equ:Airy_Abound} follows by rearranging terms. The proof of \eqref{equ:Airy_Bbound} is analogous.
\end{proof}
\end{lemma}
We record the scaled version of \eqref{equ:HHB0} and \eqref{equ:HHA0} below.
\begin{lemma}
Fix $T\geq 2$, $n\geq 2$, $1\leq k\leq n-1$,  $x>0$, and $y_2\geq y_1\geq z>-2^{-1/3}n^{1/2}T^{-1/6}+x$ . Then it holds that
\begin{align}\label{equ:AirysheetB} 
\begin{split}
  \mathfrak{h}^{T,n} (x ,y_2)-\mathfrak{h}^{T,n} (x ,y_1)-  \mathfrak{X}^{T,n} [(z ,k )\xrightarrow{(T/2)^{1/3}}  (y_2,1)]+ \mathfrak{X}^{T,n} [(z ,k )&\xrightarrow{(T/2)^{1/3}}  (y_1,1)]    \\
\geq &\log\big(1-   \mathfrak{B}^{T,n}_k(x,y_1;z)  \big) .
\end{split}
\end{align}
Similarly, it holds that
\begin{align}\label{equ:AirysheetA}  
\begin{split}
  \mathfrak{h}^{T,n} (x ,y_2)-\mathfrak{h}^{T,n} (x ,y_1)- \mathfrak{X}^{T,n} [(z ,k )\xrightarrow{(T/2)^{1/3}}  (y_2,1)]+ \mathfrak{X}^{T,n} [(z &,k )\xrightarrow{(T/2)^{1/3}}  (y_1,1)]   \\
\leq & - \log\big(1-   \mathfrak{A}^{T,n}_k(x,y_2;z)  \big) .
\end{split}
\end{align}
\begin{proof}
From \eqref{def:scaled_prelimit} and Lemma~\ref{lem:changeofvairable}, $\mathfrak{X}^{T,n} [(z ,k )\xrightarrow{(T/2)^{1/3}}  (y_2,1)]- \mathfrak{X}^{T,n} [(z ,k )\xrightarrow{(T/2)^{1/3}}  (y_1,1)]$ equals
\begin{align*}
2^{1/3}T^{-1/3}\left( \mathcal{X}^{T,n} [(2^{1/3}T^{2/3}z ,k ) {\longrightarrow} (2^{1/3}T^{2/3}y_2,1)]-\mathcal{X}^{T,n} [(2^{1/3}T^{2/3}z ,k ) {\longrightarrow} (2^{1/3}T^{2/3}y_1,1)] \right).
\end{align*}
From \eqref{equ:HHB0}, it is bounded from above by
\begin{align*}
2^{1/3}T^{-1/3}\bigg( \mathbcalboondox{h}^{T,n}  (2^{1/3}T^{2/3}x ,2^{1/3}T^{2/3}y_2 )-&\mathbcalboondox{h}^{T,n}  (2^{1/3}T^{2/3}x, 2^{1/3}T^{2/3}y_1 )\\
&-\log A^{T,n}_{k}(2^{1/3}T^{2/3}x,2^{1/3}T^{2/3}y_1;2^{1/3}T^{2/3}z) \bigg).
\end{align*}
From \eqref{def:scaled_prelimit}, the above equals
\begin{align*}
\mathfrak{h}^{T,n} (x ,y_2)-\mathfrak{h}^{T,n} (x ,y_1)  -2^{1/3}T^{-1/3} \log\left(1-  \mathfrak{B}^{T,n}_k( x, y_1; z)  \right).
\end{align*}
Together with $T\geq 2$, \eqref{equ:AirysheetB} follows by rearranging terms. The proof for \eqref{equ:AirysheetA} is similar.
\end{proof}

\end{lemma}

The next proposition provides the coupling which allows us to prove Proposition~\ref{pro:Airy_S}.

\begin{proposition}\label{pro:cpl_Airy}
Fix a sequence $\mathsf{t}_0$. Then there exists a sequence $\mathsf{n}$, a subsequence $\mathsf{t}\subset \mathsf{t}_0$ and a coupling of $\{\mathfrak{X}^{T,n}, \mathfrak{h}^{T,n}\}_{(T,n)\in\mathsf{t}\times\mathsf{n}}$, $\{\mathfrak{X}^{T},\mathfrak{h}^{T}\}_{T\in\mathsf{t}}$ and the parabolic Airy line ensemble $\mathcal{A}$ such that the following statements hold.
 
\vspace*{0.3cm}
First, fix any $T\in\mathsf{t}_0$. Almost surely $\mathfrak{X}^{T,n} $ converges to $ \mathfrak{X}^{T}$ in $C(\mathbb{N}\times\mathbb{R},\mathbb{R})$, $\mathfrak{h}^{T,n}(x,y) $ converges to $ \mathfrak{h}^{T}(x,y)$ for all $x,y\in\mathbb{Q}$, and $\mathfrak{R}^{T,n}_k( \bar{x} , -2^{-1/2}k^{1/2}x^{-1/2} ) $ converges for all $k\geq 1$ and $x,\bar{x}\in\mathbb{Q}^+ $. We denote the limits by $\mathfrak{R}^{T}_k( \bar{x} , -2^{-1/2}k^{1/2}x^{-1/2} )$. 
\vspace*{0.3cm}

Second, almost surely $ \mathfrak{X}^{T}$ converges to $\mathcal{A}$ in $C(\mathbb{N}\times\mathbb{R},\mathbb{R})$, $\mathfrak{h}^{T}(x,y) $ converges in $C(\mathbb{R}^2,\mathbb{R})$, and $\mathfrak{R}^{T}_k( \bar{x} , -2^{-1/2}k^{1/2}x^{-1/2} ) $ converges for all $k\geq 1$ and $x,\bar{x}\in\mathbb{Q}^+ $. We denote the limits by $\mathfrak{h}$ and $\mathfrak{R}_k( \bar{x} , -2^{-1/2}k^{1/2}x^{-1/2} ) $ respectively.


\vspace*{0.3cm}
Lastly, $\mathfrak{h}(0,\cdot)=\mathcal{A}_1(\cdot)$. For all $x,\bar{x}\in \mathbb{Q}^+ $, it holds almost surely 
\begin{equation}\label{equ:true}
 \lim_{k\to\infty} |k^{-1/2}\mathfrak{R}_k( \bar{x} , -2^{-1/2}k^{1/2}x^{-1/2} )|=0.  
\end{equation}

\end{proposition}

\begin{proof}
Fix $T\in\mathsf{t}_0$  and an arbitrary sequence $\mathsf{n}_0$. From Proposition~\ref{pro:kpzlineensemble}, $ \{\mathfrak{X}^{T,n}\}_{n\in\mathsf{n}_0} $ is tight in $C(\mathbb{N}\times\mathbb{R},\mathbb{R})$. From Proposition~\ref{pro:kpzsheet}, the finite-dimensional distribution of $ \{\mathfrak{h}^{T,n}\}_{n\in\mathsf{n}_0} $ is tight. From \eqref{def:frakR} and Lemma~\ref{lem:Flimit}, we have the convergence in distribution of $\mathfrak{R}^{T,n}_k( \bar{x} , -2^{-1/2}k^{1/2}x^{-1/2} )$  to 
$$\mathfrak{X}^T[(0,k+1)\xrightarrow{ (T/2)^{1/3}}   ( \bar{x} ,1)]-2^{3/2} k^{1/2} \bar{x}^{1/2}+ 2^{1/3}T^{-1/3} k\log (2^{1/3}T^{2/3}).$$ 
By the Skorokhod's representation theorem \cite[Theorem 6.7]{Bil}, we may find a subsequence $\mathsf{n}'\subset\mathsf{n}_0$ and a coupling of $ \{\mathfrak{X}^{T,n},\mathfrak{h}^{T,n}\}_{n\in\mathsf{n}'} $ such that along $\mathsf{n}'$, $ \mathfrak{X}^{T,n}$, $\mathfrak{h}^{T,n}(x,y)$ and $\mathfrak{R}^{T,n}_k( \bar{x} , -2^{-1/2}k^{1/2}x^{-1/2} ) $ converge almost surely. We note that the convergences of the latter two hold at rational points. From Proposition~\ref{pro:kpzlineensemble}, the limit of $\mathfrak{X}^{T,n}$ is distributed as the scaled KPZ line ensemble and we denote it by $\mathfrak{X}^{T}$. From Proposition~\ref{pro:kpzsheet}, we may augment the probability space to accommodate a scaled KPZ sheet $\mathfrak{h}^T$ such that $\mathfrak{h}^{T,n}(x,y)$ converges almost surely to $\mathfrak{h}^{T}(x,y)$ for all $x,y\in\mathbb{Q}$. We note that since $\mathfrak{h}^{T,n}(0,\cdot)=\mathfrak{X}^{T,n}_1(\cdot)$, we may further require $\mathfrak{h}^T(0,\cdot)=\mathfrak{X}^T_1(\cdot)$. The limits of $\mathfrak{R}^{T,n}_k( \bar{x} , -2^{-1/2}k^{1/2}x^{-1/2} )$  are denoted by $\mathfrak{R}^{T}_k( \bar{x} , -2^{-1/2}k^{1/2}x^{-1/2} )$ . Moreover, 
\begin{align*}
\mathfrak{R}^{T}_k( \bar{x} , -2^{-1/2}k^{1/2}x^{-1/2} )\overset{d}{=} & \mathfrak{X}^T[(0,k+1)\xrightarrow{ (T/2)^{1/3}}   ( \bar{x} ,1)]\\
&-2^{3/2} k^{1/2} \bar{x}^{1/2}+ 2^{1/3}T^{-1/3} k\log (2^{1/3}T^{2/3}).
\end{align*}  

\vspace*{0.2cm}

By a diagonal argument, we can find a sequence $\mathsf{n}$ and couplings of $ \{\mathfrak{X}^{T,n},\mathfrak{h}^{T,n}\}_{n\in\mathsf{n}} $, $  \mathfrak{X}^{T }$ and $\mathfrak{h}^{T }  $  for each $T\in\mathsf{t}_0$ such that along $\mathsf{n}$, the convergences in the previous paragraph hold. From now on we fix such a sequence $\mathsf{n}$. From Proposition~\ref{pro:KPZtoAiry_line}, $ \{\mathfrak{X}^{T}\}_{T\in\mathsf{t}_0} $ is tight in $C(\mathbb{N}\times\mathbb{R},\mathbb{R})$. From Proposition~\ref{pro:kpzsheet_tight}, $ \{\mathfrak{h}^{T }\}_{T\in\mathsf{t}_0} $ is tight in $C(\mathbb{R}^2,\mathbb{R})$. Similarly, $\{\mathfrak{R}^{T}_k( \bar{x} , -2^{-1/2}k^{1/2}x^{-1/2} ) \}_{T\in\mathsf{t}_0} $ is tight. By the Skorokhod's representation theorem, we can find a subsequence $\mathsf{t} \subset\mathsf{t}_0$ and a coupling such that along $\mathsf{t}$, $ \mathfrak{X}^{T }$, $\mathfrak{h}^{T }$ and $\mathfrak{R}^{T }_k( \bar{x} , -2^{-1/2}k^{1/2}x^{-1/2} ) $ converge almost surely. From Proposition~\ref{pro:KPZtoAiry_line}, the limit of $\mathfrak{X}^T$ is distributed as a parabolic Airy line ensemble and we denote it by $\mathcal{A}$. Denote by $\mathfrak{h}$ and $\mathfrak{R}_k( \bar{x} , -2^{-1/2}k^{1/2}x^{-1/2} )$ the limits of $\mathfrak{h}^T$ and $\mathfrak{R}^T_k( \bar{x} , -2^{-1/2}k^{1/2}x^{-1/2} )$ respectively. From $\mathfrak{h}^T(0,\cdot)=\mathfrak{X}^T_1(\cdot)$, we have $\mathfrak{h}(0,\cdot)=\mathcal{A}_1(\cdot)$. Moreover,
\begin{align*}
\mathfrak{R}_k( \bar{x} , -2^{-1/2}k^{1/2}x^{-1/2} )\overset{d}{=}  \mathcal{A} [(0,k+1)\xrightarrow{\infty}   ( \bar{x} ,1)]-2^{3/2} k^{1/2} \bar{x}^{1/2}.
\end{align*}
From \cite[Theorem 6.3]{DOV}, for all $\varepsilon>0$,
\begin{align*}
\sum_{k=1}^\infty \mathbb{P}\left( |\mathfrak{R}_k( \bar{x} , -k^{1/2}x^{-1/2} )|>\varepsilon k^{1/2} \right)<\infty.
\end{align*} 
Then \eqref{equ:true} follows from the  Borel-Cantelli lemma. 
\end{proof}

\begin{proof}[Proof of Proposition~\ref{pro:Airy_S}]
Let $\mathfrak{h}$ be the distributional limit of $\mathfrak{h}^T$ along some sequence $\mathsf{t}_0$. From Proposition~\ref{pro:cpl_Airy}, we can find a sequence $\mathsf{n}$, a subsequence $\mathsf{t}$ of $\mathsf{t}_0$ and a coupling of $\mathfrak{h}$ and the parabolic Airy line ensemble $\mathcal{A}$ and $\{ \mathfrak{X}^{T,n},\mathfrak{h}^{T,n} \}_{(T,n)\in \mathsf{t}\times\mathsf{n}}$ such that the assertions in Proposition~\ref{pro:cpl_Airy} hold. In particular, $\mathfrak{h}(0,\cdot)=\mathcal{A}_1(\cdot)$. From Definition~\ref{def:Airy_sheet}, we may further augment the probability space to accommodate an Airy sheet $\mathcal{S}$ such that on an event with probability one,  
\begin{equation}\label{equ:Airydiff}
\begin{split}
\lim_{k\to\infty} \mathcal{A}[(-2^{-1/2}k^{1/2}x^{-1/2},k)\xrightarrow{\infty} (y_2,1)] - \mathcal{A}[(-2^{-1/2}&k^{1/2}x^{-1/2},k)\xrightarrow{\infty} (y_1,1)]\\
&= \mathcal{S}(x,y_2)-\mathcal{S}(x,y_1), 
\end{split}
\end{equation}
for all $x>0$ and $y_1,y_2$ in $\mathbb{R}$. From now on, we fix an event $\Omega_0$ with probability one such that for each element in $\Omega_0$, all assertions in Proposition~\ref{pro:cpl_Airy} and \eqref{equ:Airydiff} hold. Our goal is to prove that when this event $\Omega_0$ occurs, 
\begin{align}\label{equ:Airy_middle3}
\mathfrak{h}(x,y_2)-\mathfrak{h}(x,y_1)=\mathcal{S}(x,y_2)-\mathcal{S}(x,y_1),
\end{align} 
for all $x>0$ and $y_1,y_2$ in $\mathbb{R}$.

\vspace*{0.3cm}
Fix $x_-<x_0$ in $\mathbb{Q}^+$ and $y_1\leq y_2 $ in $\mathbb{Q}$. We want to show that 
\begin{equation}\label{equ:B_middle3}
\mathfrak{h}(x_0,y_2)-\mathfrak{h}(x_0,y_1)\geq  \mathcal{S}(x_-,y_2)-\mathcal{S}(x_-,y_1). 
\end{equation}
Let $z_k=-2^{-1/2}k^{1/2}x_-^{-1/2}$. From \eqref{equ:AirysheetB}, we have
\begin{align}\label{equ:Airy_middle}
\begin{split} 
 \mathfrak{h}^{T,n}(x_0,y_2)-\mathfrak{h}^{T,n}(x_0,y_1)-\mathfrak{X}^{T,n} [(z_k   ,k )\xrightarrow{(T/2)^{1/3}} (y_2,1)]+ &\mathfrak{X}^{T,n} [(z_k ,k )\xrightarrow{(T/2)^{1/3}} (y_1,1)]\\
 &  \geq   \log\big(1-   \mathfrak{B}^{T,n}_k( x_0, y_1;  z_k )  \big).
\end{split}
\end{align}
From our arrangement, 
\begin{align*}
\lim_{k\to\infty}\lim_{\substack{ T\in\mathsf{t} \\ T \to \infty}}\lim_{\substack{ n\in\mathsf{n} \\ n \to \infty}} \bigg( \textup{LHS\ of\ \eqref{equ:Airy_middle}}\bigg)= \mathfrak{h}(x_0,y_2)-\mathfrak{h}(x_0,y_1)-\mathcal{S}(x_-,y_2)+\mathcal{S}(x_-,y_1).  
\end{align*}
Therefore, to prove \eqref{equ:B_middle3}, it suffices to show
\begin{align}\label{equ:B_middle2}
\liminf_{k\to\infty}\liminf_{\substack{ T\in\mathsf{t} \\ T \to \infty}}\liminf_{\substack{ n\in\mathsf{n} \\ n \to \infty}}\left(\log \mathfrak{B}^{T,n}_k( x_0, y_1;  z_k )\right)=-\infty.
\end{align}
Applying \eqref{equ:Airy_Bbound} with $x=x_0$ and $\bar{x}=x_-$, $\log \mathfrak{B}^{T,n}_k( x_0, y_1;  z_k )$ is bounded from above by
\begin{align*} 
-2^{1/2}k^{1/2}x_-^{1/2}\left( 1-x_-^{-1/2}x_0^{1/2} \right)^2 + \mathfrak{h}^{T,n}( x_-, y_1)-\mathfrak{h}^{T,n}( x_0, y_1) -\mathfrak{R}^{T,n}_k( x_-, z_k)+ \mathfrak{R} ^{T,n}_k( x_0, z_k ). 
\end{align*}
Because of \eqref{equ:true}, the above goes to $-\infty$. Therefore \eqref{equ:B_middle2} holds. A similar argument yields
\begin{equation}\label{equ:B_middle4}
\mathfrak{h}(x_0,y_2)-\mathfrak{h}(x_0,y_1)\leq \mathcal{S}(x_+,y_2)-\mathcal{S}(x_+,y_1),
\end{equation}
for all $x_0<x_+$ in $\mathbb{Q}^+$ and $y_1\leq y_2$ in $\mathbb{Q}$. As a result, \eqref{equ:Airy_middle3} holds for all $x\in\mathbb{Q}^+$ and $y_1,y_2\in\mathbb{Q}$. By the continuity, \eqref{equ:Airy_middle3} holds for all $x>0$ and $y_1,y_2\in\mathbb{R}$.

\end{proof}

\section{Proof of Theorem~\ref{thm:KPZtoLandscape}}\label{sec:proofnew}
In this section, we prove Theorem~\ref{thm:KPZtoLandscape} based on Theorem~\ref{thm:KPZtoAiry_sheet}. For $T>0$, recall that the scaled narrow wedge solutions are given by
\begin{align*}
\mathfrak{H}^T( s,x ; t,y )=2^{1/3}T^{-1/3}\mathcal{H}(Ts , 2^{1/3}T^{2/3}x ; Tt, 2^{1/3}T^{2/3}y )+(t-s)2^{1/3}T^{2/3}/24.
\end{align*}
From (3) in Proposition~\ref{thm:forcoupling} and \eqref{def:KPZ_sheet_rescaled0}, for fixed $s<t$ it holds that	
\begin{align}\label{equ:KPZ_scape_rescaled}
\mathfrak{H}^T( s,x ;t,y )\overset{d}{=}   (t-s)^{1/3}\mathfrak{h}^{(t-s)T}((t-s)^{-2/3}x,(t-s)^{-2/3}y).
\end{align}
The both sides of \eqref{equ:KPZ_scape_rescaled} are viewed as $C(\mathbb{R}^2,\mathbb{R})$-valued random variables. The linearity \eqref{equ:KPZ_linear} can be rewritten as
\begin{align}
\begin{split}
\mathfrak{H}^T( s,x ; t ,y )=2^{1/3}T^{-1/3}\log \int_{-\infty}^\infty \exp\big[2^{-1/3}T^{1/3}\big( \mathfrak{H}^T( s,x ;&\tau ,z )+\mathfrak{H}^T(\tau ,z ;t,y )\big) \big] dz\\
& + 2^{1/3}T^{-1/3}\log (2^{1/3}T^{2/3}).
\end{split}
\end{align}
We first show the finite-dimensional convergence of the scaled narrow wedge solutions $\mathfrak{H}^T$ to the directed landscape $\mathcal{L}$.
\begin{lemma}\label{lem:KPZtoLandscape}
Fix a finite set $\Lambda=\{ t_1<t_2<\dots<t_m\}$. Then $\{\mathfrak{H}^T(t_i,x ;t_j,y )\} $ converges in distribution to the directed landscape $\{\mathcal{L}((t_i,x);(t_j,y))\}$ as $T$ goes to infinity. Here $\mathfrak{H}^T( t_i,x ; t_j,y )$ and $\mathcal{L}( t_i,x ; t_j,y )$ are viewed as $C(\Lambda^2_+\times\mathbb{R}^2,\mathbb{R})$-valued random variables with $\Lambda_+^2=\{(s,t)\in\Lambda^2\, |\, s<t\}$. 
\end{lemma}
\begin{proof}
Fix a finite set $\{t_1<t_2<\dots<t_m\}$. From \eqref{equ:KPZ_scape_rescaled} and Proposition~\ref{pro:kpzsheet_tight}, $\{\mathfrak{H}^T( t_i,\cdot ; t_j,\cdot ) \}_{i<j}$ is tight. Denote by $\{\mathfrak{H}( t_i,\cdot ; t_j,\cdot ) \}_{i<j}$ a subsequential limit. By the Skorokhod's representation theorem \cite[Theorem 6.7]{Bil}, we may take a coupling such that $\mathfrak{H}^T( t_i,\cdot  ; t_j ,\cdot )$ {\color{black}jointly} converges to $\mathfrak{H}( t_i,\cdot  ; t_j ,\cdot )$ almost surely in $C(\mathbb{R}^2,\mathbb{R})$ {\color{black}for $1\leq i<j\leq m$}. From (4) in Proposition~\ref{thm:forcoupling}, $\{\mathfrak{H}( t_i,\cdot ; t_{i+1},\cdot )\}_{i=1}^{m-1}$ are independent. Moreover, from \eqref{equ:KPZ_scape_rescaled} and Theorem~\ref{thm:KPZtoAiry_sheet}, $\mathfrak{H}( t_i,\cdot  ; t_j ,\cdot )$ is distributed as an Airy sheet of scale $(t_j-t_i)^{1/3}$. In view of Corollary~\ref{cor:landscape_finite}, it remains to prove that for any $t_i<t_j<t_k$, it holds almost surely 
\begin{align}\label{equ:landscape1}
\mathfrak{H}( t_i,x  ; t_k ,y )=\max_{z\in \mathbb{R}}\left( \mathfrak{H}( t_i,x  ; t_j ,z )+\mathfrak{H}( t_j,z  ; t_k ,y ) \right).
\end{align} 
From \cite[Proposition 9.2]{DOV}, the right hand side of \eqref{equ:landscape1} is well-defined as a random variable on $C(\mathbb{R}^2,\mathbb{R})$. Moreover, it is distributed as an Airy sheet of scale $(t_k-t_i)^{1/3}$. Therefore, it suffices to show that almost surely for all $x,y\in\mathbb{R}$,
\begin{align}\label{equ:landscape2}
\mathfrak{H}( t_i,x ; t_k ,y )\geq  \max_{z\in \mathbb{R}}\left( \mathfrak{H}( t_i,x  ; t_j ,z )+\mathfrak{H}( t_j,z ; t_k ,y ) \right).
\end{align} 

Let $\Omega_0$ be the event on which the following holds. First, $\mathfrak{H}^T( t_i,\cdot  ; t_j ,\cdot )$ converges to $\mathfrak{H}( t_i,\cdot ; t_j ,\cdot )$ in $C(\mathbb{R}^2,\mathbb{R})$ for all $t_i<t_j$. Second, the right hand side of \eqref{equ:landscape1} defines a continuous function in $x$ and $y$. We will show that \eqref{equ:landscape2} holds on $\Omega_0$.   

Fix $t_i<t_j$ and $x,y\in\mathbb{R}$. Denote by $Z_j( t_i,x ; t_k,y )$ the collection of maximum points of $\mathfrak{H}( t_i,x  ; t_j ,z )+\mathfrak{H}( t_j,z  ; t_k ,y) $. Note that when $\Omega_0$ occurs, $Z_j( t_i,x ; t_k,y )\neq \varnothing$. For $M>0$, consider the event $\Omega_0\cap \{ Z_j( t_i,x ; t_k,y )\cap [-M,M]\neq \varnothing  \}$. When such an event occurs, we have 
\begin{align*}
&\max_{z\in\mathbb{R}}\left( \mathfrak{H}( t_i,x  ; t_j ,z )+\mathfrak{H}( t_j,z  ; t_k ,y ) \right)\\
=&\max_{z\in [-M,M]}\left( \mathfrak{H}( t_i,x ; t_j ,z )+\mathfrak{H}( t_j,z  ; t_k ,y ) \right)\\
{\color{black}=}&{\color{black}\lim_{T\to\infty} 2^{1/3}T^{-1/3}\log \int_{-M}^M \exp\big[2^{-1/3}T^{1/3}\big( \mathfrak{H}( t_i,x ; t_j,z )+\mathfrak{H} ( t_j,z ;  t_k,y )\big) \big]\, dz}\\
=&\lim_{T\to\infty} 2^{1/3}T^{-1/3}\log \int_{-M}^M \exp\big[2^{-1/3}T^{1/3}\big( \mathfrak{H}^T( t_i,x ; t_j,z )+\mathfrak{H}^T( t_j,z ;  t_k,y )\big) \big]\, dz\\
\leq & \lim_{T\to\infty}  \mathfrak{H}^T( t_i,x ; t_k,y )-2^{1/3} T^{-1/3}\log (2^{1/3}T^{2/3})=\mathfrak{H}( t_i,x ; t_k,y ).
\end{align*}
{\color{black}In the third equality, we used 
$$\max_{z\in [-M,M]}|(\mathfrak{H}^T( t_i,x ; t_j,z )+\mathfrak{H}^T( t_j,z ;  t_k,y ))-(\mathfrak{H}( t_i,x ; t_j,z )+\mathfrak{H}( t_j,z ;  t_k,y ))|\to 0.$$
}
Then \eqref{equ:landscape2} is shown and the proof is finished.
\end{proof}
\begin{proposition}\label{pro:kpzfull_tight}
When $T$ goes to infinity, $\mathfrak{H}^T$ is tight in $C(\mathbb{R}^4_+,\mathbb{R})$.
\end{proposition}
We assume for a moment Proposition~\ref{pro:kpzfull_tight} is valid and give the proofs of Theorem~\ref{thm:KPZtoLandscape} and Corollary~\ref{cor:multi_sol}.

\begin{proof}[Proof of Theorem~\ref{thm:KPZtoLandscape}]
In view of Proposition~\ref{pro:kpzfull_tight}, $\{\mathfrak{H}^T(s,x;t,y)\}$ is tight as $C(\mathbb{R}^4_+,\mathbb{R})$-valued random variables. From Lemma~\ref{lem:KPZtoLandscape}, any subsequential limit has the same law as the directed landscape $\mathcal{L}(s,x;t,y)$. Therefore $\mathfrak{H}^T(s,x;t,y)$ converges to $\mathcal{L}(s,x;t,y)$ in distribution.
\end{proof}

 In the rest of the section, we prove Proposition~\ref{pro:kpzfull_tight}. In the following proposition we record the spatial and temporal modulus of continuity estimates in \cite[Theorem 1.3]{CGH} and \cite[Propositions 5.1 and 5.2]{DG}.
\begin{proposition}[\cite{CGH,DG}]\label{pro:spacial_continuity}
There exist universal constants $C_0$ and $T_0$ such that the following holds. For any $T\geq T_0$, $d\in (0,1] $ and $K\geq 0$, we have
\begin{align}\label{equ:continuity1}
\mathbb{P}\left(\left| \mathfrak{H}^T(0,0;1,d)-\mathfrak{H}^T(0,0;1,0)\right|\geq Kd^{1/2} \right)\leq C_0 e^{-C_0^{-1}K^{3/2}}.
\end{align}
Moreover, for any $T\geq T_0$, $\beta\in (0,1]$ and $K\geq 0$, we have
\begin{align}\label{equ:continuity2}
\mathbb{P}\left(\left| \mathfrak{H}^T(0,0;1+\beta,0)-\mathfrak{H}^T(0,0;1,0)\right|\geq K\beta ^{1/8} \right)\leq C_0 e^{-C_0^{-1}K^{3/2}}.
\end{align}
\end{proposition}
Combining Proposition~\ref{pro:spacial_continuity} with the symmetries of $\mathfrak{H}^T$, we may replace $(0,0;1,0)$ with a general point $(s,x;t,y)$ in $\mathbb{R}^4_+$. This is the content of the corollary below. Its proof can be found in Appendix.
\begin{corollary}\label{cor:modulus}
Fix $b\geq 1$. There exist constants $r_0=r_0(b)$, $D=D(b)$ and $T_1=T_1(b)$ such that the following statements hold. Let $(s,x,t,y)\in [-b,b]^4 $ with $t-s\geq b^{-1}$ and $T\geq T_1$. Then for $d\in (0,r_0]$, we have
\begin{align}
\mathbb{P}\left(\left| \mathfrak{H}^T(s,x;t,y+d)-\mathfrak{H}^T(s,x;t,y )\right|\geq Kd^{1/2} \right)\leq D e^{-D^{-1}K^{3/2}},  \label{equ:modulus1}\\
\mathbb{P}\left(\left| \mathfrak{H}^T(s,x+d;t,y)-\mathfrak{H}^T(s,x;t,y )\right|\geq Kd^{1/2} \right)\leq D e^{-D^{-1}K^{3/2}}. \label{equ:modulus2}
\end{align}
Also, for $\beta\in (0,r_0]$, we have
\begin{align}
\mathbb{P}\left(\left| \mathfrak{H}^T(s,x;t+\beta ,y)-\mathfrak{H}^T(s,x;t,y )\right|\geq K\beta ^{1/8} \right)\leq D e^{-D^{-1}K^{3/2}}, \label{equ:modulus3} \\
\mathbb{P}\left(\left| \mathfrak{H}^T(s+\beta ,x;t,y)-\mathfrak{H}^T(s,x;t,y )\right|\geq K\beta ^{1/8} \right)\leq D e^{-D^{-1}K^{3/2}}. \label{equ:modulus4}
\end{align}
\end{corollary}

\begin{proof}[Proof of Proposition~\ref{pro:kpzfull_tight}]
For $b\in\mathbb{N}$, let $r_{b,0}=r_0(b)$ be the constant in Corollary~\ref{cor:modulus}. We may further assume {\color{black}$r_{b,0}<1$ and} $r_{b,0}$ decreases to zero as $b$ goes to infinity. Consider a sequence of compact sets $Q_1\subset Q_2\subset \dots \mathbb{R}^4_+$ with the following properties.
\begin{enumerate}
\item $Q_b\subset [-b,b]^4 \cap\{ t-s\geq b^{-1} \}$. 
\item $Q_b$ is the union of finitely many small hypercubes $Q_{b,j}$ of the form 
\begin{align*}
Q_{b,j}=\{(s,x;t,y)\in\mathbb{R}^4\, |\,  (s ,x ; t , y )  \in (s_{b,j},x_{b,j}; t_{b,j}, y_{b,j})+ [0,  r_{b,0}]^4  \}.
\end{align*} 
\item Any compact subset of $\mathbb{R}^4_+$ is contained in $Q_b$ for $b$ large enough.
\end{enumerate}
From the requirement (3) above, it suffices to show that $\mathfrak{H}^T$ restricted in $Q_b$ is tight for any $b\in\mathbb{N}$. From the requirement (2) above, it suffices to show that $\mathfrak{H}^T$ restricted in $Q_{b,j}$ is tight. From now on we fix $Q_{b,j}$. We will denote $r_{b,0}$ by $r_0$ for simplicity. Note that from Lemma~\ref{lem:KPZtoLandscape}, $\mathfrak{H}^T(s_{b,j},x_{b,j};t_{b,j},y_{b,j})$ is tight. It remains to control the modulus of continuity.

\vspace*{0.2cm}

From the estimates \eqref{equ:modulus1}-\eqref{equ:modulus4}, Proposition~\ref{pro:DV} implies that there exists a random constant $C^T$ such that almost surely for any $(s_1,x_1;t_1,y_1),(s_2,x_2;t_2,y_2)\in Q_{b,j}$
\begin{align*}
&\left| \mathfrak{H}^T(s_1,x_1;t_1,y_1)-\mathfrak{H}^T(s_2,x_2;t_2,y_2)\right|\\
\leq &C^T\bigg( |x_2-x_1|^{1/2} \log^{2/3}\left( 2r^2/|x_2-x_1| \right)+|y_2-y_1|^{1/2} \log^{2/3}\left( 2r^2/|y_2-y_1| \right)\\
&+|s_2-s_1|^{1/8} \log^{2/3}\left( 2r^8/|s_2-s_1| \right)+|t_2-t_1|^{1/8} \log^{2/3}\left( 2r^8/|t_2-t_1|\right)\bigg),
\end{align*}
where {\color{black}$r=  r_0^{1/8}$}. Moreover, there is a constant $D$ depending only on $b$ such that 
$$\mathbb{P}(C^T>K)\leq De^{-D^{-1}K^{3/2}}.$$
By the Kolmogorov-Chentsov criterion (see Corollary 16.9 in \cite{Kal}) this implies the tightness of $\mathfrak{H}^{T}$ in $Q_{b,j}$.
\end{proof}

\section{General Initial Condition}\label{sec:genini}
In this section, we provide an independent proof of Theorem~\ref{thm:single_sol}, the convergence of KPZ equation to the KPZ fixed point for general initial conditions, see also \cite{QS, V}. We then combine Theorem~\ref{thm:single_sol} and Theorem~\ref{thm:KPZtoLandscape} to prove Corollary~\ref{cor:multi_sol}, the joint convergence for multiple initial conditions.

This section is organized as follows. We prove finite dimensional convergence to the KPZ fixed point for general initial conditions in Proposition~\ref{pro:genini_onept}. We establish modulus of continuity estimates in space and in time respectively in Proposition~\ref{pro:spaceHoler} and Proposition~\ref{pro:timeHoler}. We then prove Theorem~\ref{thm:single_sol} and Corollary~\ref{cor:multi_sol}.

In the next lemma, we give an upper bound for $f\otimes_T\mathfrak{H}^T$ (defined in \eqref{equ:T_convolution}) with $f(x)=2a|x|$. A similar estimate has appeared in \cite{V}.

\begin{lemma}\label{lem:lineargrowth}
Fix $a>0$, $t>0$ and $T\geq \max \{1,t^{-1}\}$. There exists a random variable $Z$ such that 
\begin{equation}
(2a|\cdot|\otimes_T\mathfrak{H}^T)(t,y)\leq 4a|y|+Z 
\end{equation}
for all $y\in\mathbb{R}$. Moreover, $\mathbb{E}[Z^2]\leq C$ for some $C=C(a,t)$.
\begin{proof}
Let $B(x)$ be a two-sided Brownian motion which is independent to the narrow wedge solutions $\mathcal{H}(s,x;t,y)$. It is known that $B(x)$ is a stationary solution to the KPZ equation \eqref{eq:KPZ} modulo additive constants \cite{MR1462228}. For any $t>0$, there exists a random variable $Z_0(t)$ such that
\begin{align*}
\log\int\exp\left( B(x)+\mathcal{H}(0,x;t,y) \right) dx\overset{d}{=} B(y)+Z_0(t).
\end{align*} 
The random variable $Z_0(t)$ has mean $\mathbb{E}[Z_0(t)]=-\frac{t}{4!}$ \cite[(3.4)]{MR2784327}. It was established \cite[Theorem 1.3]{MR2784327} that for $t\geq 1$, $\textup{Var}(Z_0(t))\leq C t^{3/2}$ for a uniform constant $C$. Through a direct computation, we deduce for any $\nu\in\mathbb{R}$ and $T>0$
\begin{align*}
((B(2\cdot)+\nu(\cdot))\otimes_T\mathfrak{H}^T)(t,y)\overset{d}{=}B(2y+\nu t)+\nu y +Z_1(T,\nu,t),
\end{align*}
where 
$$Z_1(T,\nu,t)\overset{d}{=}2^{1/3}T^{-1/3}(Z_0(Tt)-\mathbb{E}[Z_0(Tt)])-2^{1/3}T^{-1/3}\log(2^{1/3}T^{2/3})+4^{-1}\nu^2 t.$$
For $a>0$, let $U^+_a=\sup_{x\in\mathbb{R}} (B(2x)-a|x|)$ and $U^-_a=\sup_{x\in\mathbb{R}} (-B(2x)-a|x|)$. From
\begin{align*}
 -a|x|+\nu x\leq   B(2x)+\nu x+U^-_a,\ \textup{and}\ B(2y+\nu t)+\nu y\leq  a|y|+\nu y+2^{-1}a|\nu| t+U^+_a,
\end{align*}
we deduce
\begin{align*}
((-a|\cdot|+\nu\cdot)\otimes_T\mathfrak{H}^T)(t,y)\leq  a|y|+\nu y+Z_2(T,\nu,a,t),
\end{align*}
where $Z_2(T,\nu,a,t)\overset{d}{=}\big(Z_1(T,\nu,t)\big)_++U^+_a+U^-_a+2^{-1}a|\nu|t$. Applying the above inequality to $\nu=\pm 3 a$ and using 
\begin{align*}
(\max\{f,g\} \otimes_T\mathfrak{H}^T)\leq \max\{ f\otimes_T\mathfrak{H}^T,g\otimes_T\mathfrak{H}^T \}+2^{1/3}T^{-1/3}\log 2,
\end{align*}
we conclude
\begin{align*}
(2a|\cdot| \otimes_T\mathfrak{H}^T)(t,y)\leq  4a|y|+Z,
\end{align*}
where $Z\overset{d}{=}2^{1/3}T^{-1/3}\log 2+Z_2(T,3a,a,t)+Z_2(T,-3a,a,t)$. It is direct to check that the second moment of $Z$ is bounded by a constant depending only on $a$ and $t$.
\end{proof}
\end{lemma}

In the following proposition, we prove the finite dimensional convergence for the KPZ equation to the KPZ fixed point for general initial conditions.
\begin{proposition}\label{pro:genini_onept}
Let $f(x)$ be a continuous function that satisfies $f(x)\leq C(1+|x|)$ for some $C>0$. Then the finite-dimensional marginal of $(f\otimes_T\mathfrak{H}^T)(t,y)$ converges in distribution to the one of $(f\,\overline{\otimes}\,\mathcal{L})(t,y)$.
\begin{proof}
Let $a$ be a positive number such that $f(x)\leq a|x|$ for all $|x|\geq 1$. Fix finitely many times $0<t_1<t_2<\dots< t_m$. For $1\leq i\leq m$ and $T\geq t_1^{-1}$, let $Z_{i,T}$ be the random variable given in Lemma~\ref{lem:lineargrowth} such that $\mathbb{E}[Z^2_{i,T}]\leq C_i$ and $$(2a|\cdot|\otimes_T\mathfrak{H}^T)(t_i,y)\leq 4a|y|+Z_{i,T}\ \textup{for all}\ y\in\mathbb{R}.$$ In particular, for each $i$, $\{Z_{i,T}\}_{T\geq t_1^{-1}}$ is tight. There exists a subsequence such that $Z_{i,T}$ converges in distribution to $Z_{i,\infty}$, where $\mathbb{E}[Z_{i,\infty}^2]\leq C_i$. Together with Theorem~\ref{thm:KPZtoLandscape}, we may apply the Skorokhod's representation theorem \cite[Theorem 6.7]{Bil} to take a coupling such that the following holds almost surely for $1\leq i\leq m$:
\begin{equation}\label{equ:Z_T}
Z_{i,T}\ \textup{converges to}\ Z_{i,\infty},
\end{equation}
and
\begin{equation}\label{equ:H_TT}
\mathfrak{H}^T(0,\cdot;t_i,\cdot)\ \textup{converges to}\ \mathcal{L}(0,\cdot ; t_i,\cdot)\ \textup{in}\ C(\mathbb{R}^2,\mathbb{R}).
\end{equation}
From now on, we fix a realization for which the convergences \eqref{equ:Z_T} and \eqref{equ:H_TT} hold true. We will show that for this realization, and for each $1\leq i\leq m$ and $y\in\mathbb{R}$,
\begin{equation}\label{equ:L<H}
(f\,\overline{\otimes}\,\mathcal{L})(t_i,y)\leq\liminf_{T\to\infty}(f\otimes_T\mathfrak{H}^T)(t_i,y),
\end{equation}
and
\begin{equation}\label{equ:L>H}
(f\,\overline{\otimes}\,\mathcal{L})(t_i,y)\geq\limsup_{T\to\infty}(f\otimes_T\mathfrak{H}^T)(t_i,y).
\end{equation}
\eqref{equ:L<H} and \eqref{equ:L>H} imply $(f\,\overline{\otimes}\,\mathcal{L})(t_i,y)$ converges to $(f\,\overline{\otimes}\,\mathcal{L})(t_i,y)$, which proves the desired assertion.

To prove \eqref{equ:L<H}, we compute
\begin{align*}
\max_{x\in [-M,M]} (f(x)+\mathcal{L}(0,x,t_i,y))&=\lim_{T\to\infty} 2^{1/3}T^{-1/3}\log\int_{-M}^M \exp\left( 2^{-1/3}T^{1/3}(f(x)+\mathcal{L}(0,x;t_i,y)) \right) dx\\
&=\lim_{T\to\infty} 2^{1/3}T^{-1/3}\log\int_{-M}^M \exp\left( 2^{-1/3}T^{1/3}(f(x)+\mathfrak{H}^T(0,x;t_i,y)) \right) dx\\
&\leq \liminf_{T\to\infty} (f\otimes_T\mathfrak{H}^T)(t_i,y).
\end{align*}
In the second equality, we used
\begin{align*}
\max_{x\in [-M,M]}\left|\mathfrak{H}^T(0,x;t_i,y)-\mathcal{L}(0,x;t_i,y)\right|\to 0.
\end{align*}
Sending $M$ to infinity in the above estimate yields \eqref{equ:L<H}. 

Next, we turn to \eqref{equ:L>H}. For any $r\geq 1$, we have
\begin{align*}
\exp\left( 2^{-1/3}T^{1/3}(f\otimes_T\mathfrak{H}^T)(t_i,y) \right)=\int_{|x|\leq r}\exp\left( 2^{-1/3}T^{1/3}(f(x)+\mathfrak{H}^T(0,x;t_i,y)) \right) dx\\
+\int_{|x|> r}\exp\left( 2^{-1/3}T^{1/3}(f(x)+\mathfrak{H}^T(0,x;t_i,y)) \right) dx.
\end{align*}
Using $f(x)\leq a|x|$ for $x\geq 1$, the second term is bounded from above as
\begin{align*}
&\int_{|x|> r}\exp\left( 2^{-1/3}T^{1/3}(f(x)+\mathfrak{H}^T(0,x;t_i,y)) \right) dx \\
\leq &\int_{|x|> r}\exp\left( 2^{-1/3}T^{1/3}(-ar+2a|x| +\mathfrak{H}^T(0,x;t_i,y)) \right) dx\\
\leq &\exp\left(2^{-1/3}T^{1/3}(-ar+4a|y|+Z_{i,T})\right) .
\end{align*}
Set 
\begin{align*}
Z'_{i,T}=2^{1/3}T^{-1/3}\log  \int_{|x|\leq 1}\exp\left( 2^{-1/3}T^{1/3}(f(x)+\mathfrak{H}^T(0,x;t_i,y)) \right) dx.
\end{align*} 
Using $\log(b+c)\leq \log b+b^{-1}c$ for $b,c>0$, we have
\begin{align*}
(f\otimes_T\mathfrak{H}^T)(t_i,y)\leq 2^{1/3}T^{-1/3}\log  \int_{|x|\leq r}\exp\left( 2^{-1/3}T^{1/3}(f(x)+\mathfrak{H}^T(0,x;t_i,y)) \right) dx\\
+2^{1/3}T^{-1/3}\exp\left( 2^{-1/3}T^{1/3}(-ar+4a|y|+Z_{i,T}-Z'_{i,T}) \right).
\end{align*}
Since $Z_{i,T}$ and $Z'_{i,T}$ are bounded by \eqref{equ:Z_T}, \eqref{equ:H_TT}, we may choose $r$ large (depending on the realization and $y$) such that 
\begin{align*}
(f\otimes_T\mathfrak{H}^T)(t_i,y)\leq 2^{1/3}T^{-1/3}\log  \int_{|x|\leq r}\exp\left( 2^{-1/3}T^{1/3}(f(x)+\mathfrak{H}^T(0,x;t_i,y)) \right) dx +2^{1/3}T^{-1/3}.
\end{align*}
Sending $T$ to infinity in the above estimate yields
\begin{align*}
\limsup_{T\to\infty} (f\otimes_T\mathfrak{H}^T)(t_i,y)\leq \max_{x\in [-r,r]} (f(x)+\mathcal{L}(0,x;t_i,y)\leq (f\overline{\otimes}\mathcal{L})(t_i,y). 
\end{align*}
This completes the derivation for \eqref{equ:L>H}. The proof is  finished.
\end{proof}
\end{proposition}

The following two propositions concern the spatial and temporal modulus of continuity estimates for $f\otimes_{T}\mathfrak{H}^T$.
\begin{proposition}\label{pro:spaceHoler}
Let $f(x)$ be a continuous function that satisfies $f(x)\leq C(1+|x|)$ for some $C>0$. Then there exist $T_0, K_0\geq 1$, $r_0>0$ and $D>0$ such that the following holds. For any $T\geq T_0$, $K\geq K_0$, $t_0\in [1,2]$ and $y_1,y_2\in [-2^{-1},2^{-1}]$ with $|y_1-y_2|\leq r_0$, we have
\begin{align}\label{equ:spaceHolder}
\mathbb{P}\left(\left| (f\otimes_T\mathfrak{H}^T)(t_0,y_1)-(f\otimes_T\mathfrak{H}^T)(t_0,y_2) \right|\geq K |y_1-y_2|^{1/8} \right)\leq De^{-D^{-1}K^{3/2}}.
\end{align}
\end{proposition}

\begin{proposition}\label{pro:timeHoler}
Let $f(x)$ be a continuous function that satisfies $f(x)\leq C(1+|x|)$ for some $C>0$. Then there exist $T_0, K_0\geq 1$, $r_0>0$ and $D>0$ such that the following holds. For any $T\geq T_0$, $K\geq K_0$, $t_1, t_2\in [1,2]$ and $y_0\in [-2^{-1},2^{-1}]$ with $|t_1-t_2|\leq r_0$, we have
\begin{align}\label{equ:timeHolder}
\mathbb{P}\left(\left| (f\otimes_T\mathfrak{H}^T)(t_1,y_0)-(f\otimes_T\mathfrak{H}^T)(t_2,y_0) \right|\geq K |t_1-t_2|^{1/16} \right)\leq De^{-D^{-1}K^{3/2}}.
\end{align}
\end{proposition}

We now provide the proof for Theorem~\ref{thm:single_sol} and Corollary~\ref{cor:multi_sol}, after which we finish this section with the proof for Proposition~\ref{pro:spaceHoler} and Proposition~\ref{pro:timeHoler}.

\begin{proof}[Proof of Theorem~\ref{thm:single_sol}] Proposition~\ref{pro:genini_onept} provides the finite-dimensional convergence of $f\otimes_T\mathfrak{H}^T$ to $f\,\overline{\otimes}\, \mathcal{L}$. Hence it suffices to show that $f\otimes_T\mathfrak{H}^T$ is tight in $C((0,\infty)\times\mathbb{R},\mathbb{R})$.  Combining Propositions~\ref{pro:spaceHoler}, \ref{pro:timeHoler}, and Proposition \ref{pro:DV}, we may apply the Kolmogorov-Chentsov criterion (see Theorem 23.7 in \cite{Kal}) to obtain the tightness $f\otimes_T\mathfrak{H}^T$ restricted on $[1,2]\times[-2^{-1},2^{-1}]$. Applying symmetries \eqref{equ:translations} and \eqref{equ:scale}, we conclude the tightness of $f\otimes_T\mathfrak{H}^T$ restricted on any compact subsets of $(0,\infty)\times \mathbb{R}$. The proof is finished. 
\end{proof}

\begin{proof}[Proof of Corollary~\ref{cor:multi_sol}]
Let $s_1,\dots, s_N$ and $f_1(x),\dots, f_N(x)$ be given as in Corollary~\ref{cor:multi_sol}. From Theorem~\ref{thm:single_sol}, $f_i\otimes_T \mathfrak{H}^T(s_i; t,y)$ converges in distribution to $f_i\,\overline{\otimes}\, \mathcal{L}(s_i; t,y)$. This implies $$f_i\otimes_T \mathfrak{H}^T(s_i; t,y),\ 1\leq i\leq N $$ are jointly tight. Let $(g_1(t,y),\dots, g_N(t,y))$ be the distributional limit along some sequence $\mathsf{t}$. We aim to show it has the same distribution as $$(f_1\,\overline{\otimes}\, \mathcal{L}(s_i; t,y) ,\dots, f_N\,\overline{\otimes}\, \mathcal{L}(s_i; t,y) ).$$

From Theorem~\ref{thm:KPZtoLandscape}, we may couple $\{\mathfrak{H}^T(s,x;t,y)\}_{T\in\mathsf{t}}$, $\mathcal{L}(s,x;t,y)$ and $g_i(t,y)$ in one probability space such that almost surely 
\begin{enumerate}
\item $\mathfrak{H}^T(s,x;t,y)$ converges to $\mathcal{L}(s,x;t,y)$ in $C(\mathbb{R}^4_+,\mathbb{R})$.
\item $f_i\otimes_T \mathfrak{H}^T(s_i; t,y)$ converges to $g_i(t,y)$ in $C((s_i,\infty)\times\mathbb{R},\mathbb{R})$ for $1\leq i\leq N$.
\end{enumerate}
Let $\Omega_0$ be an event with probability one such that both (1) and (2) hold. When $\Omega_0$ occurs, we derive
\begin{align*}
&(f_i\,\overline{\otimes}\, \mathcal{L})(s_i;t,y)\\
=&\sup_{x\in \mathbb{R}}(f_i(x)+\mathcal{L}(s_i,x;t,y))\\
=&\lim_{M\to\infty}\sup_{|x|\leq M}(f_i(x)+\mathcal{L}(s_i,x;t,y))\\
{\color{black}}=&{\color{black}\lim_{M\to\infty}\lim_{T\to\infty} 2^{1/3}T^{-1/3}\log\int_{-M}^{M}\exp\left( 2^{-1/3}T^{1/3}\left( f_i(x)+\mathcal{L}(s_i,x;t,y) \right) \right) dx}\\
=&\lim_{M\to\infty}\lim_{T\to\infty} 2^{1/3}T^{-1/3}\log\int_{-M}^{M}\exp\left( 2^{-1/3}T^{1/3}\left( f_i(x)+\mathfrak{H}^T(s_i,x;t,y) \right) \right) dx\\
 \leq &\lim_{T\to\infty} (f_i\otimes_T\mathfrak{H}^T)(s_i;t,y)=g_i(t,y).
\end{align*}
{\color{black}In the fourth equality, we used 
$$\max_{x\in [-M,M]}|\mathfrak{H}^T(s_i,x;t,y)-\mathcal{L}(s_i,x;t,y) |\to 0.$$
}
Since each $g_i(t,y)$ has the same law as $(f_i\,\overline{\otimes}\, \mathcal{L})(s_i;t,y)$, we get $(f_i\,\overline{\otimes}\, \mathcal{L})(s_i;t,y)= g_i(t,y)$ almost surely. This finishes the proof.
\end{proof}

To prove Propositions~\ref{pro:spaceHoler} and \ref{pro:timeHoler}, we need the following modulus of continuity estimates for $\mathfrak{H}^T$. 

\begin{lemma}\label{lem:7.5}
There exist $T_1\geq 1$, $r_1>0$ and $D_1>0$ such that the following holds. For any $T\geq T_1$, there is a random constant $C^T$ such that for all $(t_1,x_1,y_1), (t_2,x_2,y_2)\in [1,2]\times [-2,2]\times [-2,2]$ with $|t_1-t_2|\leq r_1$, $|x_1-x_2|\leq r_1$, and $|y_1-y_2|\leq r_1$, we have
\begin{align*}
|\mathfrak{H}^T(0,x_1;t_1,y_1)-\mathfrak{H}^T(0,x_2;t_2,y_2)|
\leq &C^T\left( |x_2-x_1|^{1/4} +|y_2-y_1|^{1/4}+  |t_2-t_1|^{1/16}\right).
\end{align*}
Moreover, for any $K\geq 0$ we have $\mathbb{P}(C^T>K)\leq D_1e^{-D_1^{-1}K^{3/2}}.$
\begin{proof}
We apply Corollary~\ref{cor:modulus} with $b=2$ to $\mathfrak{H}^T(0,x;t,y)$ on $(t,x,y)\in [1,2]\times[-2,2]\times[-2,2]$. Set $T_1=T_1(2)$, $r_1=r_0(2)$ and $D=D(2)$ in Corollary~\ref{cor:modulus}. Then we have
\begin{align*}
\mathbb{P}\left( \left| \mathfrak{H}^T(0,x;t,y+d)-\mathfrak{H}^T(0,x;t,y) \right|\geq Kd^{1/2} \right)\leq De^{-D^{-1}K^{3/2}},\\
\mathbb{P}\left( \left| \mathfrak{H}^T(0,x+d;t,y)-\mathfrak{H}^T(0,x;t,y) \right|\geq Kd^{1/2} \right)\leq De^{-D^{-1}K^{3/2}},\\
\mathbb{P}\left( \left| \mathfrak{H}^T(0,x;t+\beta,y)-\mathfrak{H}^T(0,x;t,y) \right|\geq K\beta^{1/8} \right)\leq De^{-D^{-1}K^{3/2}}
\end{align*} 
for all $T\geq T_1$, $(t,x,y)\in [1,2]\times [-2,2]\times [-2,2]$ and $\beta, d\in [0,r_1]$. Applying Proposition~\ref{pro:DV}, we have a random constant $C^T$ such that with probability one,
\begin{align*}
 &\left| \mathfrak{H}^T(0,x_1;t_1,y_1)-\mathfrak{H}^T(0,x_2;t_2,y_2) \right|\\
 \leq &C^T\bigg( |x_2-x_1|^{1/2}\log^{2/3}\left( 2r_1^2/|x_1-x_2| \right)+|y_2-y_1|^{1/2}\log^{2/3}\left( 2r_1^2/|y_1-y_2| \right)\\
 &+|t_2-t_1|^{1/8}\log^{2/3}\left( 2r_1^2/|t_1-t_2| \right) \bigg)
\end{align*} 
provided $(t_1,x_1,y_1), (t_2,x_2,y_2)\in [1,2]\times [-2,2]\times [-2,2]$, $|t_1-t_2|\leq r_1$, $|x_1-x_2|\leq r_1$, and $|y_1-y_2|\leq r_1$. Moreover, $\mathbb{P}(C^T\geq K)\leq D'e^{-(D')^{-1}K^{2/3}}$ for some constant $D'$. Then the assertion follows.
\end{proof}   
\end{lemma}

\begin{proof}[Proof of Proposition~\ref{pro:spaceHoler}]
Let $r_1$ and $T_1$ be the constants given in Lemma~\ref{lem:7.5}. We set $r_0=r_1$ and fix $t_0\in [1,2]$, and $y_1,y_2\in [-2^{-1},2^{-1}]$ with $0<|y_1-y_2|\leq r_0$. Let $T_0\geq T_1$ and $K_0$ be large constants to be determined and assume $T\geq T_0$ and $K\geq K_0$. We use $D$ to denote a constant that depends only on $f$, $T_0$ and $K_0$ (in particular, not on $t_0, y_1$ and $y_2$). The exact value of $D$ may increase from line to line.

Set $r=2^{-3}K^{1/2}|y_1-y_2|^{-1/8}$. To prove \eqref{equ:spaceHolder}, it is sufficient to show 
\begin{equation}\label{equ:0624}
\mathbb{P}(\mathsf{A}^{\textup{c}})\leq De^{-D^{-1}K^{3/2}},
\end{equation}
where $\mathsf{A}=\{ (f\otimes_T \mathfrak{H}^T)(t_0,y_1)\leq (f\otimes_T \mathfrak{H}^T)(t_0,y_2)+K|y_1-y_2|^{1/8} \}$. Define the events:
\begin{align*}
\mathsf{A}_1&=\left\{ \sup_{x\in\mathbb{R}}\left( \mathfrak{H}^T(0,x;t_0,y_1)+2^{-2}(x-y_1)^2 \right)\leq 2^{-4}r^2 \right\},\\
\mathsf{A}_2&=\left\{ f\otimes_T \mathfrak{H}^T(t_0,y_2)\geq -2^{-4}r^2 \right\},\\
\mathsf{A}_3&=\left\{ \sup_{|x-y_1|\leq r, |y|,|y'|\leq 2^{-1}, |y-y'|\leq r_0  }\frac{|\mathfrak{H}^T(0,x;t_0,y)-\mathfrak{H}^T(0,x;t_0,y')| }{|y-y'|^{1/4}}\leq 2^{-1} K |y_1-y_2|^{-1/8}  \right\}. 
\end{align*}
Then \eqref{equ:0624} holds true if we prove both 
\begin{equation}\label{equ:7.6}
\mathsf{A}_1\cap\mathsf{A}_2\cap \mathsf{A}_3\subset\mathsf{A}.
\end{equation}
and
\begin{equation}\label{equ:7.7}
\mathbb{P}(\mathsf{A}_1^{\textup{c}})+\mathbb{P}(\mathsf{A}_2^{\textup{c}})+\mathbb{P}(\mathsf{A}_3^{\textup{c}})\leq De^{-D^{-1}K^{3/2}}. 
\end{equation}
We begin with \eqref{equ:7.6}. We write $\exp\left( 2^{-1/3}T^{1/3}(f\otimes_T \mathfrak{H}^T)(t_0,y_1)  \right)=\mathrm{I}+\mathrm{II},$ where
\begin{align*}
\mathrm{I}=\int_{|x-y_1|\leq r} \exp\left( 2^{-1/3}T^{1/3}(f(x)+\mathfrak{H}^T(0,x;t_0,y_1)) \right)dx,\\
\mathrm{II}=\int_{|x-y_1|> r} \exp\left( 2^{-1/3}T^{1/3}(f(x)+\mathfrak{H}^T(0,x;t_0,y_1)) \right)dx,
\end{align*}
When the event $\mathsf{A}_3$ occurs, we have $\mathfrak{H}^T(0,x;t_0,y_1)\leq \mathfrak{H}^T(0,x;t_0,y_2)+2^{-1}K|y_1-y_2|^{1/8}$ for all $x$ with $|x-y_1|\leq r$. Hence
\begin{align*}
\mathrm{I}\leq \exp\left(2^{-1/3}T^{1/3}\left( (f\otimes_T\mathfrak{H}^T)(t_0,y_2) +2^{-1}K|y_1-y_2|^{1/8} \right)\right).
\end{align*}
When the event $\mathsf{A}_1$ occurs, we have $f(x)+\mathfrak{H}^T(0,x;t_0,y_1)\leq f(x)-2^{-2}(x-y_1)^2+2^{-4}r^2.$ By increasing $K_0$ if necessary, we have $f(x)\leq 2^{-4}(x-y_1)^2$ for all $|x-y_1|\geq r$. Then 
\begin{align*}
f(x)+\mathfrak{H}^T(0,x;t_0,y_1)\leq -2^{-3}(x-y_1)^2
\end{align*}
for $|x-y_1|\geq r$. Therefore,
\begin{align*}
\mathrm{II}\leq \int_{|x-y_1|> r} \exp\left( -2^{-1/3}T^{1/3}\cdot 2^{-3}(x-y_1)^2 \right)dx .
\end{align*}
By taking $T_0$ large enough, we have
$$\mathrm{II}\leq \exp\left(-2^{-1/3}T^{1/3}\cdot 2^{-3}r^2\right)$$
from a Gaussian integral bound. Using $\log(\mathrm{I}+\mathrm{II})\leq \log\mathrm{I}+\mathrm{II}/\mathrm{I}$, we obtain
\begin{align*}
(f\otimes_T\mathfrak{H}^T)(t_0,y_1)\leq &(f\otimes_T\mathfrak{H}^T)(t_0,y_2)+2^{-1}K|y_1-y_2|^{1/8}\\
&+2^{1/3}T^{-1/3}\exp\left( 2^{-1/3}T^{1/3}(-2^{-3}r^2-(f\otimes_T\mathfrak{H}^T)(t_0,y_2)) \right).
\end{align*}
When the event $\mathsf{A}_2$ occurs, the quantity above is bounded from above by
\begin{align*}
(f\otimes_T\mathfrak{H}^T)(t_0,y_2)+2^{-1}K|y_1-y_2|^{1/8}+2^{1/3}T^{-1/3}\exp\left(- 2^{-1/3}T^{1/3}\cdot  2^{-4}r^2 \right).
\end{align*}
Recall that $r=2^{-3}K^{1/2}|y_1-y_2|^{-1/8}$. By increasing $K_0$ if necessary, we have
\begin{align*}
2^{1/3}T^{-1/3}\exp\left(- 2^{-1/3}T^{1/3}\cdot  2^{-4}r^2 \right) \leq 2^{-1}K|y_1-y_2|^{1/8}.
\end{align*}
In short, when the events $\mathsf{A}_1$, $\mathsf{A}_2$, and $\mathsf{A}_3$ all occur, we have $(f\otimes_T\mathfrak{H}^T)(t_0,y_1)\leq (f\otimes_T\mathfrak{H}^T)(t_0,y_2)+K|y_1-y_2|^{1/8}$. We have proved \eqref{equ:7.6}.

Next, we bound the probabilities $\mathbb{P}(\mathsf{A}^{\textup{c}}_1),\mathbb{P}(\mathsf{A}^{\textup{c}}_2)$ and $ \mathbb{P}(\mathsf{A}^{\textup{c}}_3)$. Applying \cite[Proposition 4.2]{CGH}, we have 
\begin{equation}\label{equ:A1}
\mathbb{P}(\mathsf{A}^{\textup{c}}_1)\leq De^{-D^{-1}r^3}\leq De^{-D^{-1}K^{3/2}}. 
\end{equation}
To bound $\mathbb{P}(\mathsf{A}^{\textup{c}}_2)$, we estimate that
\begin{align*}
(f\otimes_T\mathfrak{H}^T)(t_0,y_2)= &2^{1/3}T^{-1/3}\log\int \exp\left( 2^{-1/3}T^{1/3}(f(x)+\mathfrak{H}^T(0,x;t_0,y_2)) \right) dx \\
\geq &\inf_{|x-y_2|\leq 2^{-1}} (f(x)+\mathfrak{H}^T(0,x;t_0,y_2)). 
\end{align*}
This implies 
$$\mathsf{A}_2^{\textup{c}}\subset\left\{\inf_{|x-y_2|\leq 2^{-1}} (f(x)+\mathfrak{H}^T(0,x;t_0,y_2))\leq -2^{-4}r^2\right\}.$$
Applying the tail bound of $\mathfrak{H}^T(0,x;t_0,y_2)$ in \cite{CG1,CG2}, we have 
\begin{equation}\label{equ:A2}
\mathbb{P}(\mathsf{A}^{\textup{c}}_2)\leq De^{-D^{-1}r^3}\leq De^{-D^{-1}K^{3/2}}.
\end{equation}
To bound $\mathbb{P}(\mathsf{A}_3^{\textup{c}})$, we define for $i\in\mathbb{Z}$,
\begin{align*}
Z_i=\sup_{|x-y_1-i|, |y|,|y'|\leq 2^{-1}, |y-y'|\leq r_0}\frac{|\mathfrak{H}^T(0,x;t_0,y)-\mathfrak{H}^T(0,x;t_0,y')-2it_0^{-1}(y-y')| }{|y-y'|^{1/4}}.
\end{align*}  
From \eqref{equ:translations} and \eqref{equ:skew}, it follows that every $Z_i$ has the same distribution. For $|i|\leq \lceil r\rceil$, we have
\begin{align*}
Z_i\geq &\sup_{|x-y_1-i|, |y|,|y'|\leq 2^{-1}, |y-y'|\leq r_0}\frac{|\mathfrak{H}^T(0,x;t_0,y)-\mathfrak{H}^T(0,x;t_0,y')}{|y-y'|^{1/4}}-2\lceil r\rceil\\
\geq  &\sup_{|x-y_1-i|, |y|,|y'|\leq 2^{-1}, |y-y'|\leq r_0}\frac{|\mathfrak{H}^T(0,x;t_0,y)-\mathfrak{H}^T(0,x;t_0,y')|}{|y-y'|^{1/4}}- 3^{-1}K|y_1-y_2|^{-1/8}.
\end{align*} 
Therefore, $\mathsf{A}_3^{\textup{c}}\subset\cup_{|i|\leq \lceil r\rceil} \{Z_i\geq 6^{-1}K|y_1-y_2|^{-1/8} \}$ and $$\mathbb{P}(\mathsf{A}_3^{\textup{c}})\leq (2\lceil r \rceil+1)\mathbb{P}(Z_0\geq 6^{-1}K|y_1-y_2|^{-1/8}).$$
From Lemma~\ref{lem:7.5}, we have $$\mathbb{P}(Z_0\geq 6^{-1}K|y_1-y_2|^{-1/8})\leq De^{-D^{-1}|y_1-y_2|^{-3/16}K^{3/2}}.$$
Hence,
\begin{equation}\label{equ:A3}
\begin{split}
\mathbb{P}(\mathsf{A}_3^{\textup{c}})\leq & (2\lceil r \rceil+1)  \mathbb{P}(Z_0\geq 6^{-1}K|y_1-y_2|^{-1/8})\\
\leq & D|y_1-y_2|^{-1/8}K^{1/2} e^{-D^{-1}|y_1-y_2|^{-3/16}K^{3/2}}\leq D e^{-D^{-1}K^{3/2}}. 
\end{split} 
\end{equation}
Combining \eqref{equ:A1}, \eqref{equ:A2}, and \eqref{equ:A3} yields \eqref{equ:7.7}. The proof is finished.  
\end{proof}

\begin{proof}[Proof of Proposition~\ref{pro:timeHoler}]
Let $r_1$ and $T_1$ be the constants given in Lemma~\ref{lem:7.5}. We set $r_0=r_1$ and fix $t_1, t_2\in [1,2]$ and $y_0\in [-2^{-1},2^{-1}]$ with $0<|t_1-t_2|\leq r_0$. Let $T_0\geq T_1$ and $K_0$ be large constants to be determined and assume $T\geq T_0$ and $K\geq K_0$. We use $D$ to denote a constant that depends only on $f$, $T_0$ and $K_0$ (in particular, not on $y_0, t_1$ and $t_2$). The exact value of $D$ may increase from line to line.
For $C_0\geq 1$, consider the event
\begin{align*}
\mathsf{B}_3=\left\{ \sup_{|x-y_0|\leq C_0^{-1}K^{1/2}|t_1-t_2|^{-1/8}}  |\mathfrak{H}^T(0,x;t_1,y_0)-\mathfrak{H}^T(0,x;t_2,y_0)|   \leq 2^{-1} K|t_1-t_2|^{1/16} \right\}  
\end{align*} 
\begin{claim}\label{clm:1}
For $K_0$ and $C_0$ large enough, we have 
\begin{equation}\label{equ:B3}
\mathbb{P}(\mathsf{B}_3^{\textup{c}})\leq  De^{-D^{-1}K^{3/2}}.
\end{equation}
\end{claim}
We postpone the proof of Claim~\ref{clm:1} to the end of this section. From now on, we fix $C_0$ given in Claim~\ref{clm:1} and assume $K_0$ is large enough such that \eqref{equ:B3} holds. Set $r=C_0^{-1}K^{1/2}|t_1-t_2|^{-1/8}$. To prove \eqref{equ:timeHolder}, it is sufficient to show 
\begin{equation}\label{equ:0934}
\mathbb{P}(\mathsf{B}^{\textup{c}})\leq De^{-D^{-1}K^{3/2}},
\end{equation}
where $\mathsf{B}=\{ (f\otimes_T \mathfrak{H}^T)(t_1,y_0)\leq (f\otimes_T \mathfrak{H}^T)(t_2,y_0)+K|t_1-t_2|^{1/16} \}$. Define the events:
\begin{align*}
\mathsf{B}_1&=\left\{ \sup_{x\in\mathbb{R}}\left( \mathfrak{H}^T(0,x;t_1,y_0)+2^{-2}(x-y_0)^2 \right)\leq 2^{-4}r^2 \right\},\\
\mathsf{B}_2&=\left\{ f\otimes_T \mathfrak{H}^T(t_2,y_0)\geq -2^{-4}r^2 \right\}. 
\end{align*}
Then \eqref{equ:0934} holds true if we prove both 
\begin{equation}\label{equ:7.13}
\mathsf{B}_1\cap\mathsf{B}_2\cap \mathsf{B}_3\subset\mathsf{B}.
\end{equation}
and
\begin{equation}\label{equ:7.14}
\mathbb{P}(\mathsf{B}_1^{\textup{c}})+\mathbb{P}(\mathsf{B}_2^{\textup{c}})+\mathbb{P}(\mathsf{B}_3^{\textup{c}})\leq De^{-D^{-1}K^{3/2}}. 
\end{equation}
We begin with \eqref{equ:7.13}. We write $\exp\left( 2^{-1/3}T^{1/3}(f\otimes_T \mathfrak{H}^T)(t_1,y_0)  \right)=\mathrm{I}+\mathrm{II},$ where
\begin{align*}
\mathrm{I}=\int_{|x-y_0 |\leq r} \exp\left( 2^{-1/3}T^{1/3}(f(x)+\mathfrak{H}^T(0,x;t_1,y_0 )) \right)dx,\\
\mathrm{II}=\int_{|x-y_0 |> r} \exp\left( 2^{-1/3}T^{1/3}(f(x)+\mathfrak{H}^T(0,x;t_1,y_0 )) \right)dx.
\end{align*}
When the event $\mathsf{B}_3$ occurs, we have $\mathfrak{H}^T(0,x;t_1,y_0)\leq \mathfrak{H}^T(0,x;t_2,y_0)+2^{-1}K|t_1-t_2|^{1/16}$ for $|x-y_0|\leq r$. Hence
\begin{align*}
\mathrm{I}\leq \exp\left(2^{-1/3}T^{1/3}\left( (f\otimes_T\mathfrak{H}^T)(t_2,y_0 ) +2^{-1}K|t_1-t_2|^{1/16} \right)\right).`
\end{align*}
When the event $\mathsf{B}_1$ occurs, we have $f(x)+\mathfrak{H}^T(0,x;t_1,y_0)\leq f(x)-2^{-2}(x-y_0)^2+2^{-4}r^2.$ By increasing $K_0$ if necessary, we may assume $f(x)\leq 2^{-4}(x-y_0)^2$ for all $|x-y_0|\geq r$. Then
\begin{align*}
f(x)+\mathfrak{H}^T(0,x;t_1,y_0)\leq-2^{-3}(x-y_0)^2
\end{align*}
for $|x-y_0|\geq r$. Therefore,
\begin{align*}
\mathrm{II}\leq \int_{|x-y_0|> r} \exp\left( -2^{-1/3}T^{1/3}\cdot 2^{-3}(x-y_0)^2 \right)dx .
\end{align*}
By taking $T_0$ large enough, we have
$$\mathrm{II}\leq \exp\left(-2^{-1/3}T^{1/3}\cdot 2^{-3}r^2\right)$$
from a Gaussian integral bound. Using $\log(\mathrm{I}+\mathrm{II})\leq \log\mathrm{I}+\mathrm{II}/\mathrm{I}$, we obtain
\begin{align*}
(f\otimes_T\mathfrak{H}^T)(t_1,y_0)\leq &(f\otimes_T\mathfrak{H}^T)(t_2,y_0)+2^{-1}K|t_1-t_2|^{1/16}\\
&+2^{1/3}T^{-1/3}\exp\left( 2^{-1/3}T^{1/3}(-2^{-3}r^2-(f\otimes_T\mathfrak{H}^T)(t_2,y_0)) \right).
\end{align*}
When the event $\mathsf{B}_2$ occurs, the quantity above is bounded from above by
\begin{align*}
(f\otimes_T\mathfrak{H}^T)(t_2,y_0)+2^{-1}K|t_1-t_2|^{1/16}+2^{1/3}T^{-1/3}\exp\left(- 2^{-1/3}T^{1/3}\cdot  2^{-4}r^2 \right).
\end{align*}
By increasing $K_0$ if necessary, we have
\begin{align*}
2^{1/3}T^{-1/3}\exp\left(- 2^{-1/3}T^{1/3}\cdot  2^{-4}r^2 \right)& \leq 2^{-1}K|t_1-t_2|^{1/16}.
\end{align*}
In short, when the events $\mathsf{B}_1$, $\mathsf{B}_2$, and $\mathsf{B}_3$ all occur, we have $(f\otimes_T\mathfrak{H}^T)(t_1,y_0)\leq (f\otimes_T\mathfrak{H}^T)(t_2,y_0)+K|t_1-t_2|^{1/16}$.
We have proved \eqref{equ:7.13}.

Next, we bound the probabilities $\mathbb{P}(\mathsf{B}^{\textup{c}}_1)$ and $\mathbb{P}(\mathsf{B}^{\textup{c}}_2)$. Applying \cite[Proposition 4.2]{CGH}, we have 
\begin{equation}\label{equ:B1}
\mathbb{P}(\mathsf{B}^{\textup{c}}_1)\leq De^{-D^{-1}r^3}\leq De^{-D^{-1}K^{3/2}}. 
\end{equation}
To bound $\mathbb{P}(\mathsf{B}^{\textup{c}}_2)$, we estimate that
\begin{align*}
(f\otimes_T\mathfrak{H}^T)(t_2,y_0)= &2^{1/3}T^{-1/3}\log\int \exp\left( 2^{-1/3}T^{1/3}(f(x)+\mathfrak{H}^T(0,x;t_2,y_0)) \right) dx \\
\geq &\inf_{|x-y_0|\leq 2^{-1}} (f(x)+\mathfrak{H}^T(0,x;t_2,y_0)). 
\end{align*}
This implies $$\mathsf{B}_2^{\textup{c}}\subset\left\{\inf_{|x-y_0|\leq 2^{-1}} (f(x)+\mathfrak{H}^T(0,x;t_2,y_0))\leq -2^{-4}r^2\right\}.$$
Applying the tail bound of $\mathfrak{H}^T(0,x;t_2,y_0)$ in \cite{CG1,CG2}, we have 
\begin{equation}\label{equ:B2}
\mathbb{P}(\mathsf{B}^{\textup{c}}_2)\leq De^{-D^{-1}r^3}\leq De^{-D^{-1}K^{3/2}}.
\end{equation}
Combining \eqref{equ:B1}, \eqref{equ:B2}, and \eqref{equ:B3} yields \eqref{equ:7.14}.  

\end{proof}
\begin{proof}[Proof of Claim~\ref{clm:1}]
For $i\in\mathbb{Z}$, set $x_i=i|t_1-t_2|^{1/2}+y_0$ and $I_i=[x_i,x_{i+1}]$. Let $N= \lceil C_0^{-1}K^{1/2}|t_1-t_2|^{-5/8} \rceil$. Then $$[y_0-C_0^{-1}|t_1-t_2|^{-1/8}K^{1/2},y_0+C_0^{-1}|t_1-t_2|^{-1/8}K^{1/2}]\subset\cup_{|i|\leq N}I_i.$$
Define
\begin{align*}
\mathsf{B}_{31}=&\cap_{|i|\leq N}\left\{\sup_{x\in I_i}\left| \mathfrak{H}^T(0,x;t_1,y_0)-\mathfrak{H}^T(0,x_i;t_1,y_0) \right|\leq 2^{-3}K|t_1-t_2|^{1/16}  \right\},\\
\mathsf{B}_{32}=&\cap_{|i|\leq N}\left\{\sup_{x\in I_i}\left| \mathfrak{H}^T(0,x;t_2,y_0)-\mathfrak{H}^T(0,x_i;t_2,y_0) \right|\leq 2^{-3}K|t_1-t_2|^{1/16}  \right\},\\
\mathsf{B}_{33}=&\cap_{|i|\leq N}\left\{ \left| \mathfrak{H}^T(0,x_i;t_1,y_0)-\mathfrak{H}^T(0,x_i;t_2,y_0) \right|\leq 2^{-2}K|t_1-t_2|^{1/16}  \right\}.
\end{align*}
From the triangle inequality, we have $\mathsf{B}_{31}\cap\mathsf{B}_{32}\cap\mathsf{B}_{33}\subset \mathsf{B}_{3}.$ Hence it suffices to show that for $C_0$ and $K_0$ large enough, we have $\mathbb{P}(\mathsf{B}^{\textup{c}}_{31})+\mathbb{P}(\mathsf{B}^{\textup{c}}_{32})+\mathbb{P}(\mathsf{B}^{\textup{c}}_{33})\leq De^{-D^{-1}K^{3/2}}.$

Define
\begin{align*}
Z_i=\sup_{x\in I_i}\frac{\left| \mathfrak{H}^T(0,x;t_1,y_0)-\mathfrak{H}^T(0,x_i;t_1,y_0)+2it_1^{-1}|t_1-t_2|^{1/2}(x-x_i) \right|}{|x-x_i|^{1/4}}.
\end{align*}
From \eqref{equ:skew}, it follows that every $Z_i$ has the same distribution. For $|i|\leq N$, 
\begin{align*}
Z_i\geq &\sup_{x\in I_i}\frac{\left| \mathfrak{H}^T(0,x;t_1,y_0)-\mathfrak{H}^T(0,x_i;t_1,y_0) \right|}{|x-x_i|^{1/4}}-2Nt_1^{-1}|t_1-t_2|^{1/2}|x-x_i|^{3/4}
\end{align*}
From $N= \lceil C_0^{-1}K^{1/2}|t_1-t_2|^{-5/8}  \rceil$, $|x-x_i|\leq |t_1-t_2|^{1/2}$ and $K\geq 1$, by taking $C_0$ large, we have
\begin{align*}
Z_i\geq \sup_{x\in I_i}\frac{\left| \mathfrak{H}^T(0,x;t_1,y_0)-\mathfrak{H}^T(0,x_i;t_2,y_0) \right|}{|x-x_i|^{1/4}}-2^{-4}K|t_1-t_2|^{-1/16}.
\end{align*}
When $\mathsf{B}_{31}^{\textup{c}}$ occurs, there exists $|i|\leq N$ such that 
$$\sup_{x\in I_i}\left| \mathfrak{H}^T(0,x;t_1,y_0)-\mathfrak{H}^T(0,x_i;t_2,y_0) \right|> 2^{-3}K|t_1-t_2|^{1/16}.$$
Since $|x-x_i|^{1/4}\leq |t_1-t_2|^{1/8}$ for $x\in I_i$, we have
$$\sup_{x\in I_i}\frac{\left| \mathfrak{H}^T(0,x;t_1,y_0)-\mathfrak{H}^T(0,x_i;t_2,y_0) \right|}{|x-x_i|^{1/4}}> 2^{-3}K|t_1-t_2|^{-1/16}.$$
Therefore, $\mathsf{B}_{31}^{\textup{c}}\subset \cup_{|i|\leq N}\{Z_i\geq 2^{-4}K|t_1-t_2|^{-1/16}\}$ and $$\mathbb{P}(\mathsf{B}_{31}^{\textup{c}})\leq (2N+1)\mathbb{P}(Z_0\geq 2^{-4}K|t_1-t_2|^{-1/16}).$$
From Lemma~\ref{lem:7.5}, we have $$\mathbb{P}(Z_0\geq 2^{-4}K|t_1-t_2|^{-1/16})\leq De^{-D^{-1}|t_1-t_2|^{-3/32}K^{3/2}}.$$
Therefore,
\begin{align}\label{equ:B31}
\mathbb{P}(\mathsf{B}_{31}^{\textup{c}})\leq  D|t_1-t_2|^{-5/8}K^{1/2} e^{-D^{-1}|t_1-t_2|^{-3/32}K^{3/2}}\leq De^{-D^{-1}K^{3/2}}.
\end{align}
A similar argument yields
\begin{align}\label{equ:B32}
\mathbb{P}(\mathsf{B}_{32}^{\textup{c}}) \leq De^{-D^{-1}K^{3/2}}.
\end{align}
Lastly, from \cite[Propositions 5.1 and 5.2]{DG} we have,
\begin{align*}
\mathbb{P}\left( \frac{\left| \mathfrak{H}^T(0,x;t_1,y_0)-\mathfrak{H}^T(0,x;t_2,y_0) \right|}{|t_1-t_2|^{1/8}}> 2^{-2}K|t_1-t_2|^{-1/16}\right)\leq De^{-D^{-1}|t_1-t_2|^{-3/32}K^{3/2}}
\end{align*}
for all $x$ with $|x-y_0|\leq C_0^{-1}K^{1/2}|t_1-t_2|^{-1/8}$.  Hence
\begin{align}\label{equ:B33}
\mathbb{P}(\mathsf{B}_{33}^{\textup{c}}) \leq (2N+1) De^{-D^{-1}|t_1-t_2|^{-3/32}K^{3/2}}\leq De^{-D^{-1}K^{3/2}}.
\end{align}
Combining \eqref{equ:B31}, \eqref{equ:B32} and \eqref{equ:B33} yields \eqref{equ:B3}.
\end{proof}

\section{Proof of Theorem~\ref{thm:main}}\label{sec:proof2}
In this section, we prove Theorem~\ref{thm:main} and Equation~\eqref{equ:wishhhh}, and confirm Conjecture~\ref{conj:main}. 

\vspace*{0.3cm}

We begin by giving upper bounds for the c.d.f. $A^{T,n}_k(x,y;z)$ and $B^{T,n}_k(x,y;z)$ (see Definition~\ref{def:AB}). Let
\begin{align}\label{def:Rnk}
R^{T,n}_k(x,z)\coloneqq F^{T,n}_k(x,z)-k\log x-T^{-1}zx+\log k!.
\end{align} 
The random variable $R^{T,n}_k(x,z)$ is defined so that, if \eqref{equ:wishhhh} holds true, $R^{T,n}_k(x,z)=o(k)$. See Proposition~\ref{pro:cpl_KPZ} and its proof for details.

We need the following elementary inequality. Fix $\varepsilon>0$ and $k\geq 1$. There exists a constant $D=D( \varepsilon)>0$ such that for all $T>0$ and $x,\bar{x} \in [\varepsilon, \varepsilon^{-1}]$, we have
\begin{align}\label{equ:kconcave}
k\log  \bar{x} + T^{-1}\bar{z}\bar{x} - k\log x-T^{-1}\bar{z} x \geq D^{-1}k |x-\bar{x}|^2,
\end{align} 
where $\bar{z}=-kT/\bar{x}$. 
  
\begin{lemma}\label{lem:boundz}
Fix $T, \varepsilon>0$, $n\geq 2$, and $1\leq k\leq n-1$. There exists $D=D(\varepsilon)>0$ such that for all $x,\bar{x}\in [\varepsilon, \varepsilon^{-1}]$, and $y>-n^{1/2}T^{1/2}+\max\{x, \bar{x}\}$, the following statements hold. Let $\bar{z}=-kT/\bar{x}$. If $\bar{x}\geq x$, then
\begin{align}\label{equ:z<A}
\log A^{T,n}_k(x,y;\bar{z})\leq -D^{-1}k |x-\bar{x}|^2 +\mathbcalboondox{h}^{T,n}(\bar{x},y)-\mathbcalboondox{h}^{T,n}(x,y)-R^{T,n}_k(\bar{x},\bar{z})+ R^{T,n}_k(x,\bar{z}).
\end{align}
If $\bar{x}\leq x$, then
\begin{align}\label{equ:z>B}
\log B^{T,n}_k(x,y;\bar{z})\leq -D^{-1}k |x -\bar{x}|^2  +\mathbcalboondox{h}^{T,n}(\bar{x},y)-\mathbcalboondox{h}^{T,n}(x,y)-R^{T,n}_k(\bar{x},\bar{z})+ R^{T,n}_k(x,\bar{z}).
\end{align}
\end{lemma}
\begin{proof}
First, we consider the case $\bar{x} \geq   x$. From \eqref{equ:FF<A}, we have
\begin{align*}
\mathbcalboondox{h}^{T,n}(\bar{x} ,y)-\mathbcalboondox{h}^{T,n}(x ,y)-\log A^{T,n}_k(x,y;\bar{z} ) \geq  &F^{T,n}_k(\bar{x} ,\bar{z} )-F^{T,n}_k(x,\bar{z}).
\end{align*}
Using  \eqref{def:Rnk}, the right hand side of the above equals
\begin{align*}
(k\log \bar{x} +T^{-1}\bar{z}  \bar{x} )-(k\log x+T^{-1}\bar{z} x)+R^{T,n}_k(\bar{x} ,\bar{z})-R^{T,n}_k(x,\bar{z}).
\end{align*}
Using \eqref{equ:kconcave}, it is bounded from below by
\begin{align*}
D^{-1}k |x -\bar{x}|^2+R^{T,n}_k(\bar{x} ,\bar{z} )-R^{T,n}_k(x,\bar{z}).
\end{align*}
Hence \eqref{equ:z<A} follows by rearranging terms. The proof of \eqref{equ:z>B} is similar.

\end{proof}

The next proposition provides us the coupling to prove Theorem~\ref{thm:main}.
\begin{proposition}\label{pro:cpl_KPZ}
Fix $T>0$. There exists a sequence $\mathsf{n}$, and a coupling of $\{\mathcal{X}^{T,n}, \mathbcalboondox{h}^{T,n}\}_{n \in\mathsf{n}}$, the KPZ line ensemble $\mathcal{X}^{T}$ and the KPZ sheet $\mathbcalboondox{h}^{T}$ such that the following statements hold.
 
\vspace*{0.3cm}
Almost surely, $ \mathcal{X}^{T,n}$ converges to $ \mathcal{X}^{T}$ in $C(\mathbb{N}\times\mathbb{R},\mathbb{R})$, $\mathbcalboondox{h}^{T,n}(x,y)$ converges to $\mathbcalboondox{h}^{T}(x,y)$ for all $x,y\in\mathbb{Q}$ and $ {R}^{T,n}_k( x , -Tk/\bar{x}  ) $ converge for all $k\geq 1$ and $x,\bar{x}\in\mathbb{Q}^+ $.  The limits of $ {R}^{T,n}_k( x , -Tk/\bar{x}  ) $ are denoted by $ {R}^{T}_k( x , -Tk/\bar{x}  ) $. It holds that $ \mathbcalboondox{h}^T(0,\cdot)= \mathcal{X}^T_1(\cdot)$.

Moreover, suppose \eqref{equ:wishhhh} holds. Then for all $x,\bar{x}\in \mathbb{Q}^+ $, it holds almost surely 
\begin{equation}\label{equ:wish}
\lim_{k\to\infty}  |R^{T}_k(x,-Tk/\bar{x})/k|=0.  
\end{equation}
\end{proposition}
\begin{proof}
From Proposition~\ref{pro:kpzlineensemble}, $ \{\mathcal{X}^{T,n}\}_{n\in\mathbb{N} } $ is tight in $C(\mathbb{N}\times\mathbb{R},\mathbb{R})$. From Proposition~\ref{pro:kpzsheet}, the finite-dimensional distribution of $\mathbcalboondox{h}^{T,n}$ is tight. From Lemma~\ref{lem:Flimit}, $ {R}^{T,n}_k( x , -Tk/\bar{x}  )$  has the distributional limit 
$$\mathcal{X}^T[(0,k+1) \to    ( x ,1)]- k\log x+\log k!.$$ 
By the Skorokhod's representation theorem \cite[Theorem 6.7]{Bil}, we may find a sequence $ \mathsf{n}$ and a coupling of $ \{\mathcal{X}^{T,n},\mathbcalboondox{h}^{T,n}\}_{n\in\mathsf{n}} $ such that along $\mathsf{n}$, $ \mathcal{X}^{T,n}$, $\mathbcalboondox{h}^{T,n}(x,y)$ and ${R}^{T,n}_k( x , -Tk/\bar{x}  )$ converge almost surely. We note that the convergences of the latter two hold at rational points. From Proposition~\ref{pro:kpzlineensemble}, the limit of $\mathcal{X}^{T,n}$ is distributed as the KPZ line ensemble and we denote it by $\mathcal{X}^{T}$. From Proposition~\ref{pro:kpzsheet}, we may augment the probability space to accommodate the KPZ sheet $\mathbcalboondox{h}^T$ such that $\mathbcalboondox{h}^{T,n}(x,y)$ converges to $\mathbcalboondox{h}^T(x,y)$ for all $x,y\in\mathbb{Q}$.  From $ \mathbcalboondox{h}^{T,n}(0,\cdot)=\mathcal{X}^{T,n}_1(\cdot)$, we may further require $ \mathbcalboondox{h}^{T}(0,\cdot)=\mathcal{X}^{T}_1(\cdot)$. Denote the limits of ${R}^{T,n}_k( x , -Tk/\bar{x}  )$ by ${R}^{T}_k( x , -Tk/\bar{x}  )$. From Proposition~\ref{lem:Flimit},
\begin{align*}
R^{T}_k(x,-Tk/\bar{x})\overset{d}{=}  \mathcal{X}^T[(0,k+1) \to    ( x ,1)]- k\log x+\log k!.
\end{align*}  

Suppose \eqref{equ:wishhhh} holds. This implies for all $\varepsilon>0$,
\begin{align*}
\sum_{k=1}^\infty \mathbb{P}\left( |R^{T}_k(x,-Tk/\bar{x})|>\varepsilon k  \right)<\infty.
\end{align*} 
Then \eqref{equ:wish} follows from the  Borel-Cantelli lemma. 
\end{proof}

\begin{proof}[Proof of Theorem~\ref{thm:main}]
Throughout this proof, we assume \eqref{equ:wishhhh} is valid. Fix $T>0$. From Proposition~\ref{pro:cpl_KPZ}, we can find a sequence $\mathsf{n}$ and a coupling of $\{ \mathcal{X}^{T,n},\mathbcalboondox{h}^{T,n} \}_{n\in  \mathsf{n}}$ with the following property. There exists an event $\Omega_0$ with probability one on which the statements below hold.
\begin{enumerate}
\item $\mathcal{X}^{T,n}$ converges to the KPZ line ensemble $\mathcal{X}^T$ in $C(\mathbb{N}\times\mathbb{R},\mathbb{R})$.
\item $\mathbcalboondox{h}^{T,n}(x,y)$ converges to the KPZ sheet $\mathbcalboondox{h}^T(x,y)$ for all $x,y\in\mathbb{Q}$.
\item $R^{T,n}_k(x,-Tk/\bar{x})$ converges to $R^{T}_k(x,-Tk/\bar{x})$ for all $x,\bar{x}\in\mathbb{Q}^+$ and $k\in\mathbb{N}$.
\item \eqref{equ:wish} holds.
\end{enumerate}
Our goal is to show that \eqref{equ:main} holds on $\Omega_0$. 

\vspace*{0.3cm}
Fix arbitrary $x_0<x_+$ in $\mathbb{Q}^+$ and $y_1\leq y_2$ in $\mathbb{Q}$. We claim that
\begin{equation}\label{equ:KPZ_middle}
\begin{split}
\limsup_{k\to\infty}\bigg( \mathcal{X}^T [(-Tk/x_0,k )\to (y_2,1)]- \mathcal{X}^T [&(-Tk/x_0,k )\to  (y_1,1)]\bigg)\\
&\leq  \mathbcalboondox{h}^T (x_+ ,y_2)-\mathbcalboondox{h}^T (x_+ ,y_1).
\end{split} 
\end{equation}
Let $z_k=-kT/x_0$. From \eqref{equ:HHB0}, we have
\begin{equation*} 
\begin{split} 
\mathcal{X}^{T,n}[(z_k,k)\to (y_2,1)]-\mathcal{X}^{T,n}[(z_k,k)\to (y_1,1)] & -\mathbcalboondox{h}^{T,n}(x_+,y_2)+\mathbcalboondox{h}^{T,n}(x_+ ,y_1)\\
&\leq  -\log\left(1- B^{T,n}_k(x_+,y_1,z_k)\right).
\end{split}
\end{equation*}
Let $n$ and $k$ go to infinity, we have
\begin{align*}
\limsup_{k\to\infty}\bigg( \mathcal{X}^T [(z_k,k )\to (y_2,1)]- \mathcal{X}^T &[(z_k,k )\to (y_1,1)]\bigg)-\mathbcalboondox{h}^T (x_+ ,y_2)+\mathbcalboondox{h}^T (x_+ ,y_1)\\
\leq &-\log\bigg(1- \limsup_{k\to\infty}\limsup_{\substack{ n\in\mathsf{n} \\ n \to \infty}} B^{T,n}_k(x_+,y_1;z_k)  \bigg).
\end{align*} 
To obtain \eqref{equ:KPZ_middle}, it suffice to show that the limit of $B^{T,n}_k(x_+,y_1;z_k)$ is zero. Equivalently, 
\begin{equation}\label{equ:KPZ_middle2}
\limsup_{k\to\infty}\limsup_{\substack{ n\in\mathsf{n}\\ n\to\infty }}\log B^{T,n}_k(x_+,y_1;z_k )=-\infty.
\end{equation}
Applying \eqref{equ:z>B} with $\bar{x}=x_0$ and $x=x_+$, we have
\begin{align*}
\log B^{T,n}_k(x_+,y_1;z_k )  \leq &-D^{-1}k |x_+ -x_0|^2  +\mathbcalboondox{h}^{T,n}(x_0,y_1)-\mathbcalboondox{h}^{T,n}(x_+,y_1)\\
&-R^{T,n}_k(x_0,z_k)+ R^{T,n}_k(x_+,z_k).
\end{align*}
Because of \eqref{equ:wish}, the limit of the right hand side is $-\infty$. Therefore we proved \eqref{equ:KPZ_middle2} and \eqref{equ:KPZ_middle}. For any $x_-<x_0$ in $\mathbb{Q}^+$, a similar argument yields
\begin{align}\label{equ:KPZ_middle3}
\begin{split}
\liminf_{k\to\infty}\bigg( \mathcal{X}^T [(-Tk/x_0,k )\to (y_2,1)]- \mathcal{X}^T [(-Tk /x_0,k )&\to (y_1,1)]\bigg)\\
& \geq     \mathbcalboondox{h}^T (x_- ,y_2)-\mathbcalboondox{h}^T (x_- ,y_1).  
\end{split}
\end{align}
Combining \eqref{equ:KPZ_middle} and \eqref{equ:KPZ_middle3}, we obtain \eqref{equ:main} for $x,y_1,y_2\in\mathbb{Q}$.

Next, we show that \eqref{equ:main} holds for $x\in\mathbb{Q}^+$ and $y_1,y_2\in\mathbb{R}$. Let $y_{1,j}$ and $y_{2,j}$ be a sequence of rational numbers that converge to $y_1$ and $y_2$ respectively. We further require $y_{1,j}\leq y_1$ and $y_{2,j}\geq y_2$, from Lemma~\ref{lem:2022last}, we have
\begin{align*}
\mathcal{X}^T [(-Tk/x,k )\to (y_2,1)]\leq &\mathcal{X}^T [(-Tk/x,k )\to (y_{2,j},1)]-\mathcal{X}^T_1(y_{2,j})+\mathcal{X}^T_1(y_{2}),\\
  \mathcal{X}^T [(-Tk/x,k )\to (y_1,1)]\geq &\mathcal{X}^T [(-Tk/x,k )\to (y_{1,j},1)]-\mathcal{X}^T_1(y_{1,j})+\mathcal{X}^T_1(y_{1}).
\end{align*}  
Therefore,
\begin{equation*} 
\begin{split}
&\limsup_{k\to\infty}\bigg( \mathcal{X}^T [(-Tk/x ,k )\to (y_2,1)]- \mathcal{X}^T [(-Tk /x ,k )\to (y_1,1)]\bigg)\\
 \leq  &   \mathbcalboondox{h}^T (x ,y_{2,j})-\mathbcalboondox{h}^T (x  ,y_{1,j})-\mathcal{X}^T_1(y_{2,j})+\mathcal{X}^T_1(y_{2})+\mathcal{X}^T_1(y_{1,j})-\mathcal{X}^T_1(y_{1}).  
\end{split}
\end{equation*}
Let $j$ go to infinity, we get
\begin{equation*} 
\begin{split}
\limsup_{k\to\infty}\bigg( \mathcal{X}^T [(-Tk/x ,k )\to (y_2,1)]- \mathcal{X}^T [(-Tk /x ,k )&\to (y_1,1)]\bigg)\\
 & \leq    \mathbcalboondox{h}^T (x ,y_{2})-\mathbcalboondox{h}^T (x  ,y_1).  
\end{split}
\end{equation*}
The other direction can be proved similarly.

Lastly, the condition $x\in\mathbb{Q}_+$ can be replaced by $x>0$ through noting that from Lemma~\ref{lem:f-convex}, $\mathcal{X}^T [(-Tk/x ,k )\to (y_2,1)]- \mathcal{X}^T [(-Tk /x ,k )\to (y_1,1)]$ is monotone non-decreasing in $x$.
\end{proof}

{\color{black}
\begin{proof}[Proof of \eqref{equ:wishhhh}]
We begin by recalling \eqref{equ:wishhhh}, which claims that for all $\varepsilon>0$ and $x>0$, it holds that
\begin{align*}
\sum_{k=1}^\infty\mathbb{P}\bigg( \left| \mathcal{X}^T[(0,k+1)\to (x,1)]-k\log x+\log k! \right|>\varepsilon k\bigg)<\infty.
\end{align*}
For convenience, we prove \eqref{equ:wishhhh} with $\varepsilon k$ replaced by $\varepsilon(k+1)$. Because $k\log x-\log k!$ is the volume of $\mathcal{Q}[(0,k+1)\to (x,1)]$, we have
\begin{align*}
 \{ &\left| \mathcal{X}^T[(0,k+1)\to (x,1)]-k\log x+\log k! \right|>\varepsilon (k+1) \}\\
&\subset \left\{ \max_{\pi\in\mathcal{Q}[(0,k+1)\to (x,1)]} \mathcal{X}^T(\pi)> (k+1)\varepsilon\right\}\cup\left\{ \min_{\pi\in\mathcal{Q}[(0,k+1)\to (x,1)]} \mathcal{X}^T(\pi)< -(k+1)\varepsilon\right\}.
\end{align*}
Set $\alpha=1/4$, $M_k = \varepsilon x^{-\alpha}(k+1)^\alpha $,  and
\begin{align*}
\|\mathcal{X}^T_i\|_{\alpha,[0,x]}=\sup_{y_1,y_2\in [0,x], y_1\neq y_2} \frac{|\mathcal{X}^T_i(y_1)-\mathcal{X}^T_i(y_2)|}{|y_1-y_2|^\alpha}.
\end{align*}
It can be checked that 
\begin{align*}
\left\{ \max_{\pi\in\mathcal{Q}[(0,k+1)\to (x,1)]} \mathcal{X}^T(\pi)> (k+1)\varepsilon\right\}\cup\left\{ \min_{\pi\in\mathcal{Q}[(0,k+1)\to (x,1)]} \mathcal{X}^T(\pi)< -(k+1)\varepsilon\right\}\\
\subset\left\{ \max_{1\leq i\leq k+1}\|\mathcal{X}^T_i\|_{\alpha,[0,x]}> M_k \right\}.
\end{align*}
Let $\hat{\mathcal{X}}^T_i(x)\coloneqq \mathcal{X}^T_i(x)+2^{-1}x^2.$ For $k$ large enough, we have
\begin{align*}
 \left\{ \max_{1\leq i\leq k+1}\|\mathcal{X}^T_i\|_{\alpha,[0,x]}> M_k \right\}\subset\left\{ \max_{1\leq i\leq k+1}\|\hat{\mathcal{X}}^T_i\|_{\alpha,[0,x]}> 2^{-1}M_k \right\}.
\end{align*}
In \cite[Corollary 1.5]{wu2025}, it is proved that for all $p\geq 2$, there exists $C=C(\alpha,p)>0$ such that
\begin{align*}
\mathbb{E}\|\hat{\mathcal{X}}^T_i\|_{\alpha,[0,x]}^p\leq Cx^{p/2-\alpha p} 
\end{align*}
for all $i\in\mathbb{N}$. Therefore,
\begin{align*}
\mathbb{P}\left(\max_{1\leq i\leq k+1}\|\hat{\mathcal{X}}^T_i\|_{\alpha,[0,x]}> 2^{-1}M_k \right)\leq & (k+1)\cdot 2^pM_k^{-p}\cdot C x^{p/2-\alpha p}\\
=&2^p\varepsilon^{-p}x^{p/2}C(k+1)^{1-\alpha p}.
\end{align*}
By taking $p=12$, the above is summable in $k$. This finishes the proof of \eqref{equ:wishhhh}.
\end{proof}}

\begin{appendix}
\section{}\label{sec:appendix}
In the appendix we provide proofs for basic results used in the paper.
\begin{proof}[Proof of Lemma~\ref{lem:f-convex}]
We use an induction argument on $\ell-m$. The assertion holds when $\ell=m$ because $f[(x,m)\to(y,m)]=f_m(y)-f_m(x)$. From Lemma~\ref{lem:integral}, we have
\begin{align*}
\exp\big( f[(x,\ell)\to (y,m)]\big)=  \int_x^y  \exp\big(f[(x,\ell)\to (z,m+1)]+f_m(y)-f_m(z) \big)\, dz.
\end{align*}
Hence $ e^{f_m(y)} \displaystyle \frac{d}{dy}\bigg(f[(x_2,\ell)\to(y,m)]-f[(x_1,\ell)\to(y,m)]\bigg)$ equals
\begin{align}\label{equ:12171132}
\begin{split}
&\left( \int_{x_2}^y  \exp\big(f[(x_2,\ell)\to (z,m+1)] -f_m(z) \big)\, dz \right)^{-1}\exp\bigg( f[(x_2,\ell)\to (y,m+1)]  \bigg)\\
-&\left( \int_{x_1}^y  \exp\big(f[(x_1,\ell)\to (z,m+1)] -f_m(z) \big)\, dz \right)^{-1}\exp\bigg( f[(x_1,\ell)\to (y,m+1)]  \bigg).
\end{split}
\end{align}
From the induction hypothesis, $f[(x_2,\ell)\to (z,m+1)]-f[(x_1,\ell)\to (z,m+1)]$ is non-decreasing in $z$. Therefore,
\begin{align*} 
&\int_{x_2}^y  \exp\big(f[(x_2,\ell)\to (z,m+1)] -f_m(z) \big)\, dz\\
\leq &\exp\big(f[(x_2,\ell)\to (y,m+1)] -f[(x_1,\ell)\to (y,m+1)] \big)\\
&\qquad\qquad\qquad\qquad\qquad\qquad\qquad\qquad \times\int_{x_2}^y  \exp\big(f[(x_1,\ell)\to (z,m+1)] -f_m(z) \big)\, dz\\
\leq & \exp\big(f[(x_2,\ell)\to (y,m+1)] -f[(x_1,\ell)\to (y,m+1)] \big)\\
&\qquad\qquad\qquad\qquad\qquad\qquad\qquad\qquad \times\int_{x_1}^y  \exp\big(f[(x_1,\ell)\to (z,m+1)] -f_m(z) \big)\, dz.
\end{align*}
Apply the above inequality to \eqref{equ:12171132}, we obtain $$\displaystyle \frac{d}{dy}\bigg(f[(x_2,\ell)\to(y,m)]-f[(x_1,\ell)\to(y,m)]\bigg)\geq 0.$$ 
\end{proof}
\begin{proof}[Proof of Lemma~\ref{lem:2022last}]
Consider the following measure-preserving injection from $\mathcal{Q}[(x,\ell)\to(y_1,m)]$ to $Q[(x,\ell)\to(y_2,m)]$. Given $\pi\in \mathcal{Q}[(x,\ell)\to(y_1,m)]$, let 
\begin{align*}
\bar{\pi}(t)=\left\{ \begin{array}{cc}
\pi(t), & t\in [x,y_1],\\
m, & t\in (y_1,y_2].
\end{array}  \right.
\end{align*} 
Then the assertion follows from $f(\pi) =f(\bar{\pi})-f_m(y_2)+f_m(y_1)$.
\end{proof}
 \begin{proof}[Proof of Lemma~\ref{lem:changeofvairable}]
From $\mathcal{Q}[(x,\ell)\to (y,m)]$ to $\mathcal{Q}[(a_2x+a_3,\ell)\to (a_2y+a_3,m)]$, there is a natural map given by $\pi(t)\mapsto \pi'(t)=\pi(a_2^{-1}t-a_2^{-1}a_3)$. Moreover, $d\pi=a_2^{-(\ell-m)}d\pi'$. Together with $g(\pi)=a_1f(\pi')+a_4(y-x)$, we derive
\begin{align*}
&g[(x,\ell)\xrightarrow{\beta} (y,m) ]\\
=&\beta^{-1}\log \int_{\mathcal{Q}[(x,\ell)\to (y,m)]} \exp(\beta g(\pi))d\pi\\
=&\beta^{-1}\log \int_{\mathcal{Q}[(a_2x+a_3,\ell)\to (a_2y+a_3,m)]} a_2^{-(\ell-m)} \exp(a_1 \beta  f(\pi')+a_4\beta (y-x))d\pi'\\
=&a_1\cdot f[(a_2x+a_3,\ell)\xrightarrow{a_1\beta} (a_2y+a_3,m)]+a_4(y-x)-\beta^{-1}(\ell-m)\log a_2.
\end{align*}
\end{proof} 
\begin{proof}[Proof of Lemma~\ref{lem:convex}]
Fix $\pi_i\in \mathcal{Q}[(x_i,\ell_i)\to (y_i,m_i)]$ and let $(t_{i,j})_{j\in\llbracket m_i+1,\ell_i \rrbracket}$ be the coordinates of $\pi_i$ under the identification \eqref{equ:identify}. We again follow the convention \eqref{equ:identifybd} and set $t_{i,\ell_i+1}=x_i$ and $t_{i,m_i}=y_i.$ It suffices to show that $\pi_1\prec\pi_2$ if and only if for all $j_1\in  \llbracket m_1,\ell_1  \rrbracket$ and $j_2 \in\llbracket m_2 , \ell_2 \rrbracket$ with $j_1\geq j_2$, it holds that
\begin{align}\label{equ:t1t2}
t_{1,j_1}\leq t_{2,j_2+1}.
\end{align}

Suppose $\pi_1\prec \pi_2$ fails. There exists $t_0\in (x_1,y_1)\cap (x_2,y_2)$ such that $\pi_1(t_0)\geq \pi_2(t_0)$. Set $j_i=\pi_i(t_0)$. Because $\pi_i$ are c\`adl\`ag and integer-valued, there exists $\varepsilon>0$ such that $\pi_i(t)=j_i$ for $t\in [t_0,t_0+\varepsilon)$. In view of \eqref{equ:pivalue}, this implies $t_{1,j_1}>t_{2,j_2+1}$ and $\eqref{equ:t1t2}$ is violated.  

Suppose $t_{1,j_1}> t_{2,j_2+1}$  for some $j_1\in  \llbracket m_1,\ell_1  \rrbracket$ and $j_2 \in\llbracket m_2 , \ell_2 \rrbracket$ with $j_1\geq j_2$. Because $x_1\leq x_2$ and $y_1\leq y_2$, we may assume $(t_{1,j_1+1},t_{1,j_1})$ and $(t_{2,j_2+1},t_{2,j_2})$ are non-empty by increasing $j_1$ or decreasing $j_2$ if necessary. Moreover, by further increasing $j_1$, we may assume $(t_{1,j_1+1},t_{1,j_1})\cap (t_{2,j_2+1},t_{2,j_2})$ is non-empty. In view of \eqref{equ:pivalue}, this implies there exists $t\in (x_1,y_1)\cap (x_2,y_2)$ such that $\pi_1(t)=j_1\geq \pi_2(t)=j_2$ and hence $\pi_1\prec \pi_2$ fails.
\end{proof}

\begin{proof}[Proof of Lemma~\ref{lem:z-reverse}]
For $i\in\llbracket 1,k \rrbracket$, set
$$(\tilde{x}_i,\tilde{\ell}_i)=(z-y_{k+1-i},n+1-m_{k+1-i}),\ (\tilde{y}_i,\tilde{m}_i)=(z-x_{k+1-i},n+1-\ell_{k+1-i}).$$
There is a natural measure-preserving bijection between $\mathcal{Q}[(x_{k+1-i},\ell_{k+1-i})\to (y_{k+1-i},m_{k+1-i})]$ and $\mathcal{Q}[(\tilde{x}_i,\tilde{\ell}_i)\to (\tilde{y}_i,\tilde{m}_i)]$. Given $\pi_{k+1-i}\in \mathcal{Q}[(x_{k+1-i},\ell_{k+1-i})\to (y_{k+1-i},m_{k+1-i})]$, set $\tilde{\pi}_i\in \mathcal{Q}[(\tilde{x}_i,\tilde{\ell}_i)\to (\tilde{y}_i,\tilde{m}_i)]$ as follows. Let $\tilde{\pi}_i(\tilde{y}_i)=\tilde{m}_i$ and $\tilde{\pi}_i(t)=n+1-\lim_{t'\to t+} \pi_{k+1-i}(z-t)$. It can be checked that $f(\pi_{k+1-i})=(R_zf)(\tilde{\pi}_i)$ and then $f(\pi)=(R_zf)(\tilde{\pi})$. Hence the assertion follows. 
 
\end{proof}

{\color{black}The following lemma is used to prove Lemma~\ref{lem:UVrepeatednew}.}
 
\begin{lemma}\label{lem:appendix}
Fix $n\geq 2$, $2\leq k\leq n$. Let $(U,V)$ be an endpoint pair with $U=(x,n)^{k}$. Then there exists a measure-preserving bijection from $\mathcal{Q}[U\to V]$ to $\mathcal{Q}[U_{n,k}(x)\to V]$. Moreover, for any continuous environment $f$, $f(\pi)$ is preserved under this bijection.

Similarly, let $(U,V)$ be an endpoint pair with $V=(y,1)^{k}$. Then there exists a measure-preserving bijection from $\mathcal{Q}[U\to V]$ to $\mathcal{Q}[U\to V_k(y)]$ which preserves the value of $f(\pi)$. 
\begin{proof}
We give the proof for the case $V=(y,1)^k$. The argument for the case $U=(x,n)^k$ is similar and simpler. Consider a map $\mathbf{G}:\mathcal{Q}[U\to (y,1)^{k }]\to \mathcal{Q}[U\to V_{k}(y)]$ given by the following. For $\pi=(\pi_1,\dots,\pi_{k})\in \mathcal{Q}[U\to (y,1)^{k}]$, we define $\mathbf{G}(\pi)=\bar{\pi}=(\bar{\pi}_1,\dots,\bar{\pi}_{k})$ through
\begin{align*}
\bar{\pi}_j(t)=\left\{ \begin{array}{cc}
\pi_j(t), & t<y,\\
j, & t=y.
\end{array} \right.
\end{align*} 
From the above arrangement, we have $f(\pi)=f(\bar{\pi})$. It can be checked that $\bar{\pi}\in \mathcal{Q}[U\to V_{k}(y)]$ and that $\mathbf{G}$ is a bijection. Moreover, $\mathbf{G}$ is the restriction of a projection map between Euclidean spaces. This implies $\mathbf{G}$ is measure-preserving. 
\end{proof}
\end{lemma}

\begin{proof}[Proof of Lemma~\ref{lem:UVrepeatednew}]
The assertion holds by applying Lemma~\ref{lem:appendix} repeatedly.
\end{proof}

\begin{proof}[Proof of Lemma~\ref{lem:searrow}]
Consider the map from $\mathcal{Q}[V'_{k }(x)\to V_{k }(y)]$ to $\mathcal{Q} [(x,1)\searrow (y,k+1)]$ given by the following. For $\pi=(\pi_1,\dots,\pi_{k })\in \mathcal{Q}[V'_{k }(x)\to V_{k }(y)]$, let $\rho$ be defined in the way such that for all $t\in [x,y]$, 
$\{\pi_1(t),\pi_2(t),\dots, \pi_{k }(t),\rho(t) \}=\llbracket 1,k+1 \rrbracket$. It is straightforward to check that $\rho$ belongs to $ \mathcal{Q} [(x,1)\searrow (y,k+1)]$. Moreover, this map is a measure-preserving bijection. Together with $$f(\pi)+f (\rho) =\sum_{i=1}^{k+1} f_i(y)-f_i(x),$$ we have 
\begin{align*}
f[V'_{k-1}(x)\to V_{k-1}(y)]&=\log\int_{\mathcal{Q}[V'_{k }(x)\to V_{k }(y)] }\exp(f(\pi))\, d\pi\\
&=\sum_{i=1}^{k+1} f_i(y)-f_i(x)+\log \int_{\mathcal{Q} [(x,1)\searrow (y,k+1)]} \exp(-f (\rho ))\, d\rho \\
&=f[V_{k+1}(x)\to V_{k+1}(y)] -f [(x,1)\searrow (y,k+1)]. 
\end{align*}
\end{proof}

\begin{proof}[Proof of Corollary~\ref{cor:gRSKinvnew}]
Consider $U_0=\{(x_1,n)^{i_1} ,\dots, (x_\ell,n)^{i_\ell} \}$ and $V_0=\{(y_1,1)^{j_1},\dots, (y_m,1)^{j_m} \}$ with $\sum_{p=1}^\ell i_p=\sum_{q=1}^m j_q=k$, $x_1<x_2<\dots<x_\ell$ and $y_1<y_2<\dots<y_m$. We aim to show that
\begin{align*}
f[U_0\to V_0]=(\mathcal{W}f)[U_0\to V_0].
\end{align*} 
Let $$U=\{U_{n,i_1}(x_1), U_{n,i_2}(x_2),\dots, U_{n,i_\ell}(x_\ell)\},\ V=\{V_{j_1}(y_1),V_{j_2}(y_2),\dots, V_{j_m}(y_m)\}.$$ From Lemma~\ref{lem:UVrepeatednew}, it suffices to show
\begin{align*}
f[U\to V]=(\mathcal{W}f)[U\to V].
\end{align*} 
For $\varepsilon>0$, we consider
\begin{align*}
U^\varepsilon_p\coloneqq &\{ (x_p-i_p\varepsilon, n),(x_p-(i_p-1)\varepsilon, n),\dots, (x_p-\varepsilon, n) \},\ U^\varepsilon\coloneqq \{ U_1^\varepsilon, U_2^\varepsilon,\dots, U^\varepsilon_\ell \}\\
V^\varepsilon_q\coloneqq &\{ (y_q+\varepsilon, 1),(y_q+2\varepsilon, 1),\dots, (y_q+j_q \varepsilon, 1) \},\ V^\varepsilon\coloneqq \{ V_1^\varepsilon,V_2^\varepsilon,\dots, V^\varepsilon_m \}.
\end{align*}
For $\varepsilon$ small enough, $(U^\varepsilon,V^\varepsilon)$ is an endpoint pair. From Proposition~\ref{pro:gRSKinv},
\begin{align*}
f[U^\varepsilon\to V^\varepsilon]=(\mathcal{W}f)[U^\varepsilon\to V^\varepsilon].
\end{align*} 
Let $N=2^{-1}\sum_{p=1}^\ell i_p(i_p-1)+2^{-1}\sum_{q=1}^m j_q(j_q-1)$. It is then sufficient to show that for some constant $c_0>0$
\begin{align}
f[U \to V ]=&\lim_{\varepsilon\to 0}  f[U^\varepsilon\to V^\varepsilon]-\log(c_0\varepsilon^{N}), \label{equ:middlegRSK1} \\ 
(\mathcal{W}f)[U \to V ]=&\lim_{\varepsilon\to 0} (\mathcal{W}f)[U^\varepsilon\to V^\varepsilon]-\log(c_0\varepsilon^{N}). \label{equ:middlegRSK2}
\end{align}
We derive \eqref{equ:middlegRSK1} below and the proof for \eqref{equ:middlegRSK2} is similar.
\begin{align*}
\exp\big( f[U^\varepsilon\to V^\varepsilon]\big)=&\exp\big( f[U^\varepsilon\to U ]+f[U\to V ]+f[V\to V^\varepsilon ]\big) \\
&+\sum_{(U',V')\neq (U,V)}\exp\big( f[U^\varepsilon\to U' ]+f[U'\to V' ]+f[V'\to V^\varepsilon ]\big).
\end{align*}
Here the summation is taken over $U'=\{U_1',U_2',\dots, U_\ell'\}$, $U_p'\in\mathcal{V}_{n,i_p}(x_p)$ and $V'=\{V_1',V_2',\dots, V_m'\}$, $V_q'\in\mathcal{V}_{n,j_q}(y_q)$ such that $(U',V')$ is an endpoint pair. It is clear that $\mathcal{Q}[U^\varepsilon\to U]\times \mathcal{Q}[V\to V^\varepsilon ]$, as a subset in a Euclidean space, has dimension $N$ with volume $c_0\varepsilon^{N}$ for some constant $c_0>0$. Also, for $(U',V')\neq (U,V)$, $\mathcal{Q}[U^\varepsilon\to U]\times \mathcal{Q}[V\to V^\varepsilon ]$ has dimension at least $N+1$ with volume $O(\varepsilon^{N+1})$. Moreover, $f(\pi)=o(1)$ for $\pi$ in $\mathcal{Q}[U^\varepsilon\to U]$, $\mathcal{Q}[U^\varepsilon\to U']$, $\mathcal{Q}[V \to V^\varepsilon]$, or $\mathcal{Q}[V'\to V^\varepsilon]$. Therefore, it holds that
\begin{align*}
\exp\big( f[U^\varepsilon\to V^\varepsilon]\big)=c_0\varepsilon^{N}\exp\big( f[U \to V ]\big)+o(\varepsilon^{N}).
\end{align*}
Taking a logarithm yields \eqref{equ:middlegRSK1}.

 \end{proof}
 
To prove Corollary~\ref{cor:modulus}, we first translate the symmetries of $\mathcal{H}$ in (3) of Proposition~\ref{thm:forcoupling} to the symmetries for $\mathfrak{H}^T$. They are simple consequences of \eqref{def:KPZ_4_scaled}. 
\begin{lemma}\label{lem:symmetry}
As a $C(\mathbb{R}^4_+,\mathbb{R})$-valued random variable, $\mathfrak{H}^T$ has the following symmetries.
\begin{equation}\label{equ:flip}
\mathfrak{H}^T(-t,y;-s,x) \overset{d}{=}\mathfrak{H}^T(s,x;t,y).
\end{equation}
For any $(r,u)\in\mathbb{R}^2$,
\begin{equation}\label{equ:translations}
\mathfrak{H}^T(s+r,x+u;t+r,y+u)\overset{d}{=}\mathfrak{H}^T(s,x;t,y).
\end{equation}
For any $\tau>0$, 
\begin{equation}\label{equ:scale}
\tau^{1/3} \mathfrak{H}^{\tau T}(\tau^{-1} s, \tau^{-2/3}x ;\tau^{-1}t,\tau^{-2/3} y) \overset{d}{=}\mathfrak{H}^T(s,x;t,y) .
\end{equation}
For any $\nu\in\mathbb{R}$,
\begin{equation}\label{equ:skew}
\mathfrak{H}^T(s,x+\nu s;t,y+\nu t) \overset{d}{=}\mathfrak{H}^T(s,x;t,y)-2\nu (y-x)-\nu^2(t-s).
\end{equation}

\end{lemma}  
 
 \begin{proof}[Proof of Corollary~\ref{cor:modulus}]
It is sufficient to prove \eqref{equ:modulus1}-\eqref{equ:modulus4} for $K\geq 1$ and hence we assume $K\geq 1$ throughout this proof. For simplicity, we write $D$ for a constant that depends only on $b$. The exact value of $D$ may increase from line to line. We always assume $d,\beta\in (0,r_0]$ and $T\geq T_1$ for some small constant $r_0$ and large constant $T_1$ to be determined. 

We begin with \eqref{equ:modulus1}. Let $z=y-x$ and $\tau=t-s$. From \eqref{equ:translations},
\begin{align*}
\mathfrak{H}^T(s,x;t,y+d)-\mathfrak{H}^T(s,x;t,y )\overset{d}{=} \mathfrak{H}^T(0,0;\tau ,z+d)-\mathfrak{H}^T(0,0;\tau,z ).
\end{align*}
From \eqref{equ:skew} with $\nu=\tau^{-1}z$, the above has the same distribution as
\begin{align*}
\mathfrak{H}^T(0,0;\tau ,d)-\mathfrak{H}^T(0,0;\tau,0 )-2\tau^{-1}zd.
\end{align*}
From \eqref{equ:scale}, the above has the same distribution as
\begin{align*}
\tau^{1/3}\left( \mathfrak{H}^{\tau T}(0,0;1 ,\tau^{-2/3} d)-\mathfrak{H}^{\tau T}(0,0;1,0 )\right)-2\tau^{-1}zd.
\end{align*}
Therefore,
\begin{align*}
&\mathbb{P}\left(\left| \mathfrak{H}^T(s,x;t,y+d)-\mathfrak{H}^T(s,x;t,y )\right|\geq Kd^{1/2} \right) \\
\leq &\mathbb{P}\left(\left| \mathfrak{H}^{\tau T}(0,0;1,\tau^{-2/3}d )-\mathfrak{H}^{\tau T}(0,0;1,0 )\right|\geq \tau^{-1/3} Kd^{1/2}-2\tau^{-4/3}|z|d \right).
\end{align*}
{\color{black}Recall that $d\in (0,r_0]$.} By taking $r_0$ small enough, we have $\tau^{-1/3} Kd^{1/2}-2\tau^{-4/3}|z|d\geq 2^{-1}\tau^{-1/3} Kd^{1/2}$. Therefore, the above is bounded by
\begin{align*}
\mathbb{P}\left(\left| \mathfrak{H}^{\tau T}(0,0;1,\tau^{-2/3}d )-\mathfrak{H}^{\tau T}(0,0;1,0 )\right|\geq 2^{-1}\tau^{-1/3} Kd^{1/2}\right).
\end{align*}
By taking $r_0$ small enough and $T_1$ large enough, we have $\tau^{2/3} d\leq 1$ and $\tau T\geq T_0$. Here $T_0$ is the constant in Proposition~\ref{pro:spacial_continuity}. Hence we can apply \eqref{equ:continuity1} to obtain
\begin{align*}
\mathbb{P}\left(\left| \mathfrak{H}^{\tau T}(0,0;1,\tau^{-2/3}d )-\mathfrak{H}^{\tau T}(0,0;1,0 )\right|\geq 2^{-1}\tau^{-1/3} Kd^{1/2}\right)\leq De^{-D^{-1}K^{3/2}}.
\end{align*}
This implies \eqref{equ:modulus1}. In view of \eqref{equ:flip}, \eqref{equ:modulus2} also holds true.

Next, we turn to \eqref{equ:modulus3}. From \eqref{equ:translations},
\begin{align*}
\mathfrak{H}^T(s,x;t+\beta ,y)-\mathfrak{H}^T(s,x;t,y )\overset{d}{=} \mathfrak{H}^T(0,0;\tau +\beta ,z)-\mathfrak{H}^T(0,0;\tau,z ).
\end{align*}
Hence
\begin{equation}\label{equ:modulus3_1}
\begin{split}
&\mathbb{P}\left(\left| \mathfrak{H}^T(s,x;t+\beta ,y)-\mathfrak{H}^T(s,x;t,y )\right|\geq K\beta ^{1/8} \right)\\
\leq &\mathbb{P}\left(\left| \mathfrak{H}^T(0,0;\tau +\beta ,z)-\mathfrak{H}^T(0,0;\tau  , (\tau+\beta )^{-1} \tau z )\right|\geq 2^{-1}K\beta ^{1/8} \right)\\
&+ \mathbb{P}\left(\left| \mathfrak{H}^T(0,0;\tau  , (\tau+\beta )^{-1} \tau z )-\mathfrak{H}^T(0,0;\tau  ,z)\right|\geq 2^{-1}K\beta ^{1/8} \right).
\end{split}
\end{equation}
Note that in the second term on the right hand side of \eqref{equ:modulus3_1}, the temporal variables are the same and the spatial variables are different. Arguing similar to the proof of \eqref{equ:modulus1}, we have for $r_0$ small enough and $T_1$ large enough,
\begin{align}\label{equ:modulus3_2}
\mathbb{P}\left(\left| \mathfrak{H}^T(0,0;\tau  , (\tau+\beta )^{-1} \tau z )-\mathfrak{H}^T(0,0;\tau  ,z)\right|\geq 2^{-1}K\beta ^{1/8} \right)\leq D e^{-D^{-1}K^{3/2}}.
\end{align}
It remains to deal with the first term on the right hand side of \eqref{equ:modulus3_1}. From \eqref{equ:skew} with $\nu=(\tau+\beta)^{-1}z$, we have
\begin{align*}
\mathfrak{H}^T(0,0;\tau +\beta ,z)-\mathfrak{H}^T(0,0;\tau  , (\tau+\beta )^{-1} \tau z )\overset{d}{=}\mathfrak{H}^T(0,0;\tau +\beta ,0)&-\mathfrak{H}^T(0,0;\tau  , 0 )\\
&-(\tau+\beta)^{-2}z^2\beta. 
\end{align*}
From \eqref{equ:scale}, it has the same distribution as
\begin{align*}
\tau^{1/3}\left( \mathfrak{H}^{\tau T}(0,0;1 +\tau^{-1} \beta ,0)-\mathfrak{H}^{\tau T}(0,0;1, 0 )\right)-(\tau+\beta)^{-2}z^2\beta. 
\end{align*}
Hence
\begin{align*}
&\mathbb{P}\left(\left| \mathfrak{H}^T(0,0;\tau +\beta ,z)-\mathfrak{H}^T(0,0;\tau  , (\tau+\beta )^{-1} \tau z )\right|\geq 2^{-1}K\beta ^{1/8} \right)\\
\leq &\mathbb{P}\left(\left| \mathfrak{H}^{\tau T}(0,0;1 +\tau^{-1} \beta ,0)-\mathfrak{H}^{\tau T}(0,0;1,0)\right|\geq 2^{-1} \tau^{-1/3} K\beta ^{1/8}- \tau^{-1/3} (\tau+\beta)^{-2}z^2\beta \right).
\end{align*}
From \eqref{equ:continuity2}, we have for $r_0$ small enough and $T_1$ large enough,
\begin{align}\label{equ:modulus3_3}
\mathbb{P}\left(\left| \mathfrak{H}^T(0,0;\tau +\beta ,z)-\mathfrak{H}^T(0,0;\tau  , (\tau+\beta )^{-1} \tau z ) \right|\geq 2^{-1}K\beta ^{1/8} \right)\leq D e^{-D^{-1}K^{3/2}}.
\end{align}
Combining \eqref{equ:modulus3_1}, \eqref{equ:modulus3_2} and \eqref{equ:modulus3_3}, we obtain \eqref{equ:modulus3}. In view of \eqref{equ:flip}, \eqref{equ:modulus4} also holds true.
\end{proof}
\end{appendix}
\printbibliography
\end{document}